\documentclass{gtart_sp}
\pdfoutput=1

\usepackage[all]{xy}
\usepackage{graphicx}

\title[Geometry of contact transformations and domains]{Geometry of contact 
transformations and domains:\\orderability versus squeezing}

\author{Yakov Eliashberg}
\givenname{Yakov}
\surname{Eliashberg}
\address{Department of Mathematics\\Stanford University\\\newline
Stanford, CA 94305-2125\\USA}
\email{eliash@math.stanford.edu}
\urladdr{}

\author{Sang Seon Kim}
\givenname{Sang Seon}
\surname{Kim}
\address{Departamento de Matematica\\Instituto Superior
Tecnico\\\newline
Av Roviso Pais\\1049-001 Lisboa\\Portugal}
\email{sskim@math.ist.utl.pt}
\urladdr{}

\author{Leonid Polterovich}
\givenname{Leonid}
\surname{Polterovich}
\address{School of Mathematical Sciences\\The Raymond and
Beverly Sackler Faculty of Exact Sciences\\Tel Aviv University\\
69978 Tel Aviv\\Israel}
\email{polterov@math.tau.ac.il}
\urladdr{}

\volumenumber{10}
\issuenumber{}
\publicationyear{2006}
\papernumber{41}
\lognumber{0714}
\startpage{1635}
\endpage{1747}

\doi{}
\MR{}
\Zbl{}

\keyword{contact manifolds}
\keyword{contact squeezing and orderability}
\keyword{Floer homology}
\keyword{holomorphic curves}
\subject{primary}{msc2000}{53D10}
\subject{primary}{msc2000}{53D40}
\subject{secondary}{msc2000}{53D35}
\subject{secondary}{msc2000}{53D50}

\received{12 February 2006}
\revised{30 September 2006}
\published{28 October 2006}
\publishedonline{28 October 2006}
\accepted{4 September 2006}
\proposed{Eleny Ionel}
\seconded{Tomasz Mrowka, Peter Ozsv\'ath}
\corresponding{Eliashberg}
\editor{CPR}
\version{}

\arxivreference{math.SG/0511658}

\let\xysavmatrix\xymatrix
\def\xymatrix{\disablesubscriptcorrection\xysavmatrix}
\AtBeginDocument{\let\bar\wbar\let\tilde\wtilde\let\hat\what}

\makeatletter
\def\cnewtheorem#1[#2]#3{\newtheorem{#1}{#3}[section]
\expandafter\let\csname c@#1\endcsname\c@thm}

\newtheorem{thm}{Theorem}[section]
\cnewtheorem{lemma}[thm]{Lemma}
\cnewtheorem{prop}[thm]{Proposition}
\newtheorem{proposition}[thm]{Proposition}
\cnewtheorem{cor}[thm]{Corollary}
\newtheorem{conjecture}[thm]{Conjecture}
\cnewtheorem{lem}[thm]{Lemma}
\cnewtheorem{sublem}[thm]{Sublemma}

\theoremstyle{definition}
\cnewtheorem{rem}[thm]{Remark}
\cnewtheorem{exam}[thm]{Example}
\cnewtheorem{example}[thm]{Example}
\cnewtheorem{defin}[thm]{Definition}

\makeatother  

\newcommand{\N}{{\mathbb{N}}}
\newcommand{\PP}{{\mathbb{P}}}
\newcommand{\II}{{\mathbb{I}}}
\newcommand{\JJ}{{\mathbb{J}}}

\newcommand{\CH}{{\mathsf{CH}}}
\newcommand{\SH}{{\mathsf{SH}}}

\newcommand{\FH}{{\mathsf{FH}}}
\newcommand{\GFH}{{\mathsf{GFH}}}

\newcommand{\sC}{{\mathtt{C}}}
\newcommand{\sK}{{\mathtt{K}}}
\newcommand{\sP}{{\mathtt{P}}}
\newcommand{\sQ}{{\mathtt{Q}}}

\newcommand{\bV}{{\overline{V}}}
\newcommand{\oc}{{\overline{c}}}
\newcommand{\uc}{{\underline{c}}}

\newcommand{\cA}{{\mathcal{A}}}

\newcommand{\cF}{{\mathcal{F}}}
\newcommand{\cG}{{\mathcal{G}}}
\newcommand{\cH}{{\mathcal{H}}}

\newcommand{\cP}{{\mathcal{P}}}
\newcommand{\cS}{{\mathcal{S}}}
\newcommand{\cT}{{\mathcal{T}}}
\newcommand{\cU}{{\mathcal{U}}}

\newcommand{\cX}{{\mathcal{X}}}

\newcommand{\cL}{{\mathcal{L}}}
\newcommand{\cR}{{\mathcal{R}}}

\newcommand{\ocH}{{\overrightarrow{\mathcal{H}}}}
\newcommand{\ocG}{{\overrightarrow{G}}}

\newcommand{\mon}{{\mathsf{mon}}}
\newcommand{\incl}{{\mathsf{in}}}

\newcommand{\Ker}{{\mathrm {Ker\,}}}
\newcommand{\Mc}{\mathcal{M}}
\newcommand{\Pc}{\mathcal{P}}
\newcommand{\p}{{\partial}}
\newcommand{\wt}{\wtilde}
\newcommand{\wh}{\what}
\newcommand{\e}{\varepsilon}
\newcommand{\CZ}{{\mathrm{CZ}}}

\newcommand{\ppp}{\frac{\partial}{\partial p}}
\newcommand{\ppq}{\frac{\partial}{\partial q}}
\newcommand{\pps}{\frac{\partial}{\partial s}}
\newcommand{\ppt}{\frac{\partial}{\partial t}}

\def\g{\gamma}

\newcommand{\codim}{\operatorname{codim}}

\newcommand{\image}{{\mathrm{Image}}}
\newcommand{\spec}{{\mathrm{spec}}}
\newcommand{\Spec}{\mathrm{spec}}
\newcommand{\st}{\mathrm{st}}
\newcommand{\id}{{\text{{\bf 1}}}}

\newcommand{\Cl}{{\mathrm{Closure}}}
\newcommand{\Supp}{\mathrm{Supp}}

\newcommand{\Span}{\mathrm{Span}}

\newcommand{\Core}{{\mathrm{Core}}}
\newcommand{\Cont}{{\hbox{\it Cont\,}}}
\newcommand{\tCont}{\widetilde{\hbox{\it Cont}\, }}
\newcommand{\sgrad}{{\mathrm{sgrad}}}
\newcommand{\Graph}{{\mathrm{graph}}}

\newcommand{\Con}{{\mathrm{Con}}}

\newcommand{\psucceq}{\mathop{\succeq}\limits_+}

\newcommand{\gothic}{\mathfrak}

\begin{document}

\begin{abstract} Gromov's famous non-squeezing theorem (1985) states
that the standard symplectic ball cannot be symplectically
squeezed into any cylinder of smaller radius. Does there exist an
analogue of this result in contact geometry? Our main finding is
that the answer depends on the sizes of the domains in question:
We establish contact non-squeezing on large scales, and show that
it disappears on small scales. The algebraic counterpart of the
(non)-squeezing problem for contact domains is the question of
existence of a natural partial order on the universal cover of the
contactomorphisms group of a contact manifold. In contrast to our
earlier beliefs, we show that the answer to this question is very
sensitive to the topology of the manifold. For instance, we prove
that the standard contact sphere is non-orderable while the real
projective space is known to be orderable. Our methods include a
new embedding technique in contact geometry as well as a
generalized Floer homology theory which contains both cylindrical
contact homology and Hamiltonian Floer homology. We discuss links
to a number of miscellaneous topics such as topology of free loops
spaces, quantum mechanics and semigroups.
\end{abstract}

\maketitle

\cl{\small\it Dedicated to Dusa McDuff on the occasion of 
her $\hbox{60}^{\it th}$ birthday}

{\leftskip 25pt\small\hyperlink{Err}{Erratum attached}\par}

\section{Introduction and main results} \label{sec-intro}

\subsection{Contact (non)-squeezing} \label{subsec-sq}

Consider the standard symplectic vector space $\R^{2n}$ endowed
with the symplectic form $\omega = dp \wedge dq=\sum\limits_1^n
dp_i\wedge dq_i$. We often identify $\R^{2n}$ with $\C^n$ and
write $z = p+iq$ for the complex coordinate. Symplectic embeddings
preserve the volume, and hence the Euclidean ball
$$B^{2n}(R_1): = \{\pi|z|^2 < R_1\}$$ cannot be symplectically embedded into
$B^{2n}(R_2)$ if $R_2<R_1$. Gromov's famous non-squeezing theorem
states that there are much more subtle obstructions for symplectic
embeddings and, in particular, $B^{2n}(R_1)$ cannot be
symplectically embedded into the cylinder
$$C^{2n}(R_2) : = B^{2}(R_2)\times \R^{2n-2}$$
when $R_2 < R_1$, see \cite{Gro}. This result led to the first
non-trivial invariants of symplectic domains in dimension $2n \geq
4$.

 In the present paper we address the question whether there
are any analogues of non-squeezing results in contact geometry.
Consider {\it the prequantization space} of $\R^{2n}$, that is the
contact manifold $V = \R^{2n} \times S^1, \; S^1 = \R/\Z$, with
contact structure $\xi = \mathrm{Ker}(dt-\alpha)$ where $\alpha$ is
the Liouville form $\frac{1}{2}(pdq-qdp)$. Given a subset $D
\subset \R^{2n}$, write $\wh{D} = D \times S^1$ for its
prequantization. The naive attempt to extend the non-squeezing
from $D$ to $\wh{D}$ fails. It is is easy to show (see \fullref{thm-all-sq} and \fullref{subseq-loopc}) that for {\it
any} $R_1,R_2 > 0$ there exists a contact embedding of
$\wh{B}(R_1)$ into $\wh{B}(R_2)$ which, for $n>1$, is isotopic to
the inclusion through smooth embeddings into $V$. Furthermore, due
to the conformal character of the contact structure, the domain
$\wh{B}(R)$ can be contactly embedded into an arbitrarily small
neighborhood of a point in $V$ (see \fullref{cor-small-nbhd}
below).

 However, the situation becomes more sophisticated if one
considers only those contact embeddings which come from globally
defined compactly supported contactomorphisms of $(V,\xi)$. We
write $\cG = \Cont (V,\xi)$ for the group of all such
contactomorphisms.

  Given two open subsets $U_1$ and $U_2$ of a contact
manifold $V$, we say that $U_1$ {\it can be squeezed} into $U_2$
if there exists a contact isotopy $\Psi_t\co  \Cl (U_1) \to V,\; t
\in [0,1],$ such that $\Psi_0 = \id$ and
$$\Psi_1 (\Cl (U_1)) \subset U_2.$$ The isotopy $\{\Psi_t\}$ is called
a {\it contact squeezing} of $U_1$ into $U_2$. If, in addition, $W
\subset V$ is an open subset such that $\Cl (U_2) \subset W$ and
$\Psi_t(\Cl (U_1)) \subset W$ for all $t$, we say that $U_1$ can be
squeezed into $U_2$ {\it inside} $W$. If the closure of $U_1$ is
compact, the ambient isotopy theorem (see, for instance, Geiges
\cite{Ge}) guarantees that any squeezing of $U_1$ into $U_2$ inside
$W$ extends to a contactomorphism from $\cG$ whose support lies in
$W$.\footnote{If the group $\cG$ is not connected than the possibility
to squeeze by an isotopy is stronger than by a global
contactomorphism. All squeezing and non-squeezing results in this
paper are proven in the strongest sense, ie,  squeezing is always
done by a contact isotopy while in our non-squeezing results we prove
non-existence of the corresponding global contactomorphism.}

\begin{rem}\label{remconvexity} {
\rm If a domain $U$ has a {\it convex contact boundary} then it admits
a {\it contact squeezing inside itself}. Let us recall that a
hypersurface $\Sigma$ in a contact manifold is called {\it convex}
(see Eliashberg and Gromov \cite{Eliashberg??}) if there exists a
contact vector field $X$ which is transversal to $\Sigma$. Note that
the vector field $-X$ is also contact, and hence one cannot assign to
a convex hypersurface any canonical co-orientation. Giroux showed (see
\cite{Giroux-convex}) that in a $3$--dimensional contact manifold any
co-orientable surface can be made convex by a generic $C^\infty$--small
perturbation. On the other hand, it is easy to check that the boundary
of a domain $\wh D\subset\R^{2n}\times S^1$ is never convex.}
\end{rem}
 Our main results concerning the contact squeezing problem
are given in the next theorems.

\begin{thm}[Non-Squeezing]\label{thm-nonsq}
Assume that $R_2 \leq m \leq R_1$ for some positive integer $m$.
Then the closure of $\wh{B}^{2n}(R_1)$ cannot be mapped into
$\wh{C}^{2n}(R_2)$ by a contactomorphism from $\cG$. In
particular, $\wh{B}^{2n}(R_1)$ cannot be squeezed into
$\wh{C}^{2n}(R_2)$.
\end{thm}

 As a counterpoint to this result, we prove

\begin{thm}[Squeezing] \label{thm-sq} Assume that $2n \geq
4$. Then $\wh{B}^{2n}(R_1)$ can be squeezed into
$\wh{B}^{2n}(R_2)$ for all $R_1,R_2 < 1$.
\end{thm}
\begin{rem}\label{remn=on}
{\rm The restriction $n>1$ is essential: it was shown by Eliashberg in
\cite{Eliash-shape} that $\wh{B}^{2}(R_1)$ cannot be squeezed into
$\wh{B}^{2}(R_2)$ for any $R_1>R_2$.}
\end{rem}
 We do not know whether $\wh{B}^{2n}(R_1)$ can be squeezed
into $\wh{B}^{2n}(R_2)$ or $\wh{C}^{2n}(R_2)$ when
$$m+1 > R_1
>R_2>m$$ for an integer $m \geq 1$.

\begin{thm} \label{thm-nosmallsq} Assume that
$$R_2 \leq \frac{m}{k} \leq R_1 < R_3 < \frac{m}{k-1}$$ for some integers
$k,m \geq 1$. Then the closure of $\wh{B}^{2n}(R_1)$ cannot be
mapped into $\wh{B}^{2n}(R_2)$ by any contactomorphism $\Phi \in
\cG$ with $\Phi\big{(}\wh{B}^{2n}(R_3)\big{)} = \wh{B}^{2n}(R_3)$.
In particular, $\wh{B}^{2n}(R_1)$ cannot be squeezed into
$\wh{B}^{2n}(R_2)$ inside $\wh{B}^{2n}(\frac{m}{k-1})$.
\end{thm}

 In the case $m=1, k>1$ this result imposes a restriction on a
squeezing of $\wh{B}^{2n}(\frac{1}{k})$ into itself guaranteed by
\fullref{thm-sq}. Roughly speaking, such a squeezing requires
some extra room. As we will see in \fullref{rem-sharp} below,
this restriction is sharp: 

\cl{\sl $\wh{B}^{2n}(\frac{1}{k})$ can be
squeezed into itself inside $\wh{B}^{2n}(\rho)$ for any $\rho >
1/(k-1) $.}

The proofs of Theorems \ref{thm-nonsq} and
\ref{thm-nosmallsq} are based on cylindrical contact homology
theory (see Eliashberg, Givental and Hofer \cite{SFT}, Ustilovsky
\cite{Ustilovsky}, Bourgeois \cite{B} and Yau \cite{Yau})
which is discussed in Sections \ref{subsec-ch} and
\ref{secflavors} below. The (non)-squeezing phenomenon described
above is closely related to the geometry of the group of
contactomorphisms of the standard sphere $S^{2n-1}$, see Sections
\ref{subsec-cont-dom} and \ref{subsec-semigroup} below. This
resembles the link between symplectic non-squeezing and the
geometry of the group of symplectomorphisms which was explored
by Lalonde and McDuff in \cite{LM1,LM2}.

\subsection{Negligible domains and symplectic capacities} \label{subsec-capac}

We say that a domain $X \subset V$ is {\it negligible} if every
bounded open subset $U$ with $\mathrm{Closure}(U) \subset X$ can be
contactly squeezed into $\wh{B}^{2n}(r)$ for any $r>0$. We start
this section with the following generalization of \fullref{thm-sq} above.

\begin{thm} \label{thm-sq-strong} The cylinder $\wh{C}^{2n}(1)$
is negligible when $2n \geq 4$.
\end{thm}

 This result is sharp: indeed, for $R>1$ the cylinder
$\wh{C}^{2n}(R)$ contains $\wh{B}^{2n}(R')$ with $R' >1$ and hence
it is not negligible due to \fullref{thm-nonsq}.

 The proof is given in \fullref{sec-sq-strong} below.  With
this result at hand, we can present the transition from non-squeezing
to squeezing in terms of symplectic capacities. For a bounded domain
$U \subset \R^{2n}$ define $\uc (U)$ as the supremum of $R$ such that
the ball $B^{2n}(R)$ can be symplectically squeezed into $U$, and $\oc
(U)$ as the infimum of $R$ such that $U$ can be symplectically
squeezed into the cylinder $C^{2n}(R)$. These quantities are small
modifications of the standard symplectic capacities. In particular,
they are symplectic capacities in the symplectic category $\cal{OH}$
whose objects are open subsets of $\R^{2n}$ and morphisms are
symplectic embeddings induced by compactly supported Hamiltonian
diffeomorphisms (see Cieliebak, Hofer, Latschev and Schlenk
\cite[Section 2.1]{CHLS} for the definition of a symplectic capacity
in a symplectic category). Define the {\it contact squeezing number}
${sq}(U)$ as the infimum of $b \in \R^+$ such that the domain
$\wh{b^{-1/2}\cdot U} \subset V$ is negligible.  Let us emphasize that
since $V$ is the prequantization space of $\R^{2n}$, every compactly
supported Hamiltonian isotopy of $\R^{2n}$ lifts to a compactly
supported contact isotopy of $V$.  Therefore the contact squeezing
number is invariant under compactly supported Hamiltonian
diffeomorphisms. The next result is an immediate consequence of
Theorems \ref{thm-nonsq} and \ref{thm-sq-strong}.

\begin{thm}\label{thm-capacities}
$$\uc (U) \leq {sq}(U) \leq \oc(U)$$
for every bounded domain $U \subset \R^{2n}$.
\end{thm}

As an immediate consequence of the theorem we get that the contact
squeezing number is a symplectic capacity in category $\cal{OH}$.

\subsection{Preliminaries in contact geometry}
\label{subsec-geom-ze}

Let $(P,\eta)$ be a contact manifold with a co-oriented contact
structure. Its symplectization $SP$ is defined as the set of all
non-zero covectors from $T^*P$ whose kernel equals the contact
hyperplane and which agree with its co-orientation. One checks
that $SP$ is a symplectic submanifold of $T^*P$ if and only if
$\eta$ is a contact structure. Thus it inherits from $T^*P$ the
canonical Liouville 1--form $\alpha$ whose differential $\omega =
d\alpha$ is the symplectic form. Note also that $SP$ is a
(trivial) principal $\R_+$--bundle over $P$. The vector field $L$
generating the $\R_+$--action is called the Liouville field on
$SP$. It satisfies $i_{L}\omega = \alpha$, and hence $\R_+$ acts
by conformally symplectic transformations. Any contact form for
$\eta$ is a section of the bundle $SP \to P$. Its graph forms a
hypersurface in $SP$ which is called a {\it starshaped}
hypersurface.

Let $\beta$ be a contact form on $P$. Then using the $\R_+$--action
one can identify $SP$ with $P \times \R_+$: the point $(x,\theta)
\in SP$, where $x \in P$ and $\theta$ is a contact covector at
$x$, corresponds to $(x, \theta/\beta) \in P \times \R_+$. After
such an identification we write a point of $SP$ as $(x,u) \in P
\times \R_+$ and call $(x,u)$ {\it the canonical coordinates on
$SP$ associated to the contact form $\beta$}. In canonical
coordinates $\alpha = u\beta$, $\omega = du \wedge \beta + u
d\beta$, $L=
\partial/\partial{u}$ and $\mathrm{graph}(\beta) = \{u=1\}$.

 The following example is crucial for understanding what
is going on below:
\begin{exam} \label{SS}
{\rm Consider the standard contact sphere $(S^{2n-1}, \eta)$ where the
sphere $S^{2n-1}$ is identified with $\partial B^{2n}(1) \subset
\R^{2n}$ and $\eta = \mathrm{Ker}(\alpha\big{|}_{S^{2n-1}})$. Its
symplectization can be identified by an $\R_+$--equivariant
symplectomorphism with $(\R^{2n}\setminus\{0\}, dp \wedge dq)$, where
the $\R_+$--action on $\R^{2n}\setminus\{0\}$ is given by $z \to
\sqrt{c}z$ for all $c \in \R_+$. Here every contact covector
$(z,s\alpha) \in SP$, where $z \in P$ and $s > 0$, corresponds to the
point $\sqrt{s}z \in \R^{2n}\setminus\{0\}$.}
\end{exam}

$\R_+$--equivariant Hamiltonian functions on $SP$ are called {\it
contact Hamiltonians}.  Every contactomorphism of $(P,\eta)$ uniquely
lifts to an $\R_+$--equivariant symplectomorphism of $SP$. Moreover,
there is a one-to-one correspondence between paths $\{f_t\}_{t \in
[0,1]}$ of compactly supported contactomorphisms with $f_0 = \id$ and
contact Hamiltonian functions $F\co SP \times [0;1] \to \R$ which
vanish outside $\tau^{-1}(K) \times [0;1]$ where $K \subset P$ is a
compact subset and $\tau\co SP \to P$ the natural projection.  We say
the isotopy $\{f_t\}$ is generated by the contact Hamiltonian $F$.

Note that to every contact form $\beta$ on $P$ corresponds a
unique contact Hamiltonian on $SP$ which equals $1$ on
$\mathrm{graph}(\beta)$. The projection of its Hamiltonian field to
$P$ is a well defined vector field $R$ called {\it the Reeb vector
field} of $\beta$. It is determined by the conditions $\beta(R)=1$
and $i_R d\beta = 0$.

Let $Q \subset P$ be a hypersurface which is transversal to the
contact structure. Every contact plane $\eta(x),\; x \in P,$
carries a conformally canonical symplectic structure given by the
differential of any contact form. The symplectic complement in
$\eta(y), \; y \in Q$ to the hyperplane $T_xQ \cap \eta(x) \subset
\eta(x)$ defines a field of lines $l(x) \subset T_xQ$. This field
of lines integrates to a one-dimensional foliation which is called
{\it the contact characteristic foliation} of $Q$.

\subsection{A partial order on contact transformations }
\label{subsec-geom}

The results on contact (non)-squeezing presented in
\fullref{subsec-sq} above are closely related to the geometry of the
group of contactomorphisms of the standard contact sphere
$S^{2n-1}$. In order to present these applications we need to recall
some preliminaries from Eliashberg and Polterovich \cite{EP}.

 Let $(P,\eta)$ be a contact manifold. Denote by $\Cont_0 (P,\eta)$
the identity component of the group of all compactly supported
contactomorphisms of $(P,\eta)$. Let $\tCont _0(P,\eta)$ be its
universal cover. Given $\tilde{f},\tilde{g} \in \tCont _0(P,\eta) $ we
say that $\tilde{f}\succeq \tilde{g}$ if the element
$\tilde{f}\tilde{g}^{-1}$ is represented by a path generated by a
non-negative contact Hamiltonian. In other words, $\tilde{f}\succeq
\tilde{g}$ if $\tilde{f}$ can be reached from $\tilde{g}$ by moving
every point in the non-negative direction with respect to the contact
structure $\eta$.  Clearly, the relation $\succeq$ on $\tCont
_0(P,\eta)$ is reflexive and transitive. As shown in \cite{EP}, for
certain closed contact manifolds (eg, the unit cotangent bundle of the
$n$--torus) it defines a genuine partial order. For the purpose of
this discussion we shall call such manifolds {\it orderable}. The next
proposition from \cite{EP} gives a useful necessary and sufficient
condition for orderability of a closed contact manifold.

\begin{prop} \label{prop-nonderability}
Let $(P,\eta)$ be a closed contact manifold. The following
conditions are equivalent:
\begin{itemize}
\item[(i)] $(P,\eta)$ is non-orderable;
\item[(ii)] there exists a contractible loop
$\phi\co  S^1 \to \Cont _0 (P,\eta)$ with $\phi(0) = \id$ which is
generated by a strictly positive contact Hamiltonian.
\end{itemize}
\end{prop}

It follows from Givental's theory of the non-linear Maslov
index that the standard contact projective space $\R P^{2n-1}$ is
orderable (see \cite{Gi,EP}). In view of this, the authors
of \cite{EP} tended to believe that its double cover $S^{2n-1}$ is
also orderable. Interestingly enough, this is not the case:

\begin{thm} \label{thm-nogo} Assume that $2n \geq 4$.
There exists a positive contractible loop of contactomorphisms of
the standard contact sphere $S^{2n-1}$. In particular, the sphere
is not orderable.
\end{thm}

 After the first version of this paper appeared in the
arXiv, E Giroux informed us that \fullref{thm-nogo} can
be extracted from the existing literature as follows. The standard
contact structure on $S^{2n-1}$ coincides with the field of
maximal complex subspaces tangent to the sphere, where we identify
$$S^{2n-1} = \partial B^{2n} \subset \C^n\;.$$
Thus the group of complex automorphisms of the unit ball $B^{2n}$
acts by contact transformations on the boundary sphere. This group
is isomorphic to $PU(n,1)$. The cone of non-negative contact
Hamiltonians on the sphere restricts to a tangent cone $C$ in the
Lie algebra $pu(n,1)$, which is invariant under the adjoint
representation.  It turns out that 25 years ago G Olshanskii
\cite{Ol}  completely characterized  those cones which give
rise to a genuine partial order on the universal cover of the
group. Applying Olshanskii's criterion, one gets that $C$ does not
generates a genuine partial order. Thus, {\it a fortiori}
$S^{2n-1}$ is non-orderable. For reader's convenience, we present
more details on the Oshanskii criterion and illustrate its
application to non-orderability of $S^3$ in Appendix B. In the
paper, we chose another route and use a method which enables us to
extend  \fullref{thm-nogo} to more general contact manifolds
(cf \fullref{rem-aptw-tw} in Appendix B). This extension is
presented in \fullref{sec-more-orderability} below (see
\fullref{thmtw-stabilization} whose proof occupies \fullref{sectw-stab}).

Existence of positive contractible loops is a manifestation
of ``symplectic flexibility". However, such loops themselves
exhibit a rigid behavior. We illustrate this in the case of the
standard contact sphere $S^{2n-1}$ with $n \geq 2$. Let
$$\Delta = \{f_{t,s}\},\; t \in S^1, s \in [0,1],$$ be a homotopy
of a positive contractible loop $\{f_{t,1}\}$ to the constant loop
$\{f_{t,0}\} \equiv \id$. Assume that $f_{0,s} = \id$ for all $s$.
Write $F_s ,\; s \in[0,1],$ for the contact Hamiltonian on
$(\R^{2n}\setminus\{0\}) \times S^1$ generating the loop
$\{f_{t,s}\},\; t \in S^1$. Put
\begin{equation}\label{eq-mu}
\mu(\Delta) := -\min_{z,s,t} \frac{F_s(z,t)}{\pi|z|^2}\,.
\end{equation}
\begin{thm} \label{mu} Assume that $2n \geq 4$. Then
\begin{itemize}
\item[(i)]$\mu(\Delta)\geq 1$ for every homotopy $\Delta$ of a positive
contractible loop of contactomorphisms of the sphere $S^{2n-1}$ to
the constant loop;
\item[(ii)] moreover, this estimate is sharp:
$\inf_{\Delta} \mu(\Delta) = 1$.
\end{itemize}
\end{thm}

We shall see in \fullref{subsec-cont-dom} that the
inequality $\mu(\Delta) \geq 1$ follows from Non-Squeezing \fullref{thm-nonsq}.

\subsection{Liouville manifolds}\label{secLiouville}

In this section we introduce the class of Liouville-fillable
contact manifolds which are in the focus of our study. Let $(M,
\omega)$ be a connected symplectic manifold which satisfies the
following conditions:
\begin{itemize}
\item[(i)] There exists a complete vector field $L$ on $M$ such
that $\omega=d\alpha$ with $\alpha= i_L\omega$. This implies that
the flow $L^t, \; t \in \R$, of $L$ acts on $M$ by conformally
symplectic diffeomorphisms.
\item[(ii)] There exists a closed connected hypersurface $P \subset M$ which
is transversal to $L$ and bounds an open domain $U \subset M$ with
compact closure such that $M = U \sqcup \bigcup_{t\geq 0} L^tP$.
\end{itemize}
The vector field $L$ is called a {\it Liouville field}, and a
symplectic manifold $(M,\omega)$ with a fixed Liouville field $L$
(or, which is the same, with a fixed primitive $\alpha=i_L\omega$
which is called a {\it Liouville form}) is called a {\it Liouville
manifold}. We will call {\it starshaped} any hypersurface $P$ and
any domain $U$ in the Liouville manifold $(M,\omega,L)$ which
satisfy the condition (ii). Given any starshaped domain $U\subset
M$ its repeller $\Core_P(M) =\bigcap_{t\in\R^+}L^{-t}(U)$ is
called the {\it core} of the Liouville manifold $(M,\omega, L)$.

Put $\beta:= \alpha|_P$. The transversality condition in (ii) is
equivalent to the requirement that $\eta:= \mathrm{Ker}{\beta}$ is a
contact structure on $P$. Moreover, the symplectization $SP$ can
be naturally symplectically identified with the set $ M_{*,P} =
\bigcup_{t \in \R} L^{t}(P)$ as follows: Consider the splitting
$SP = P \times (0,\infty)$ associated to the form $\beta$ (see
\fullref{subsec-geom-ze} above) and identify a point $(x,u) \in
SP$ with $L^{\log u}x \in M $. Under this identification the form
$\alpha$ and the vector field $L$ on $M$ correspond to the
Liouville form and Liouville vector field on $SP$ respectively.
Note that in coordinates $(x,u)$ we have $\alpha = u\beta$.

We claim that the decomposition $M = M_{*,P} \sqcup \Core_P(M)$
does not depend on the choice of a starshaped hypersurface $P$.
Indeed, let $Q$ be another starshaped hypersurface. Consider the
subset $X = M_{*,P} \cap M_{*,Q}$, and denote by $Y$ its
projection to $P$ along the trajectories of $L$. Note that $X$ is
non-empty: otherwise we have $M_{*,Q} \subset U$ which is absurd
since $M_{*,Q}$ has infinite volume and $U$ has compact closure by
condition (ii) above. Hence $Y$ is an open non-empty subset of
$P$. Furthermore, since $M \setminus U =\bigcup_{t\geq 0} L^tP$,
every trajectory of $L$ starting on $P$ leaves any compact subset
$K \supset U$ in finite time. The same holds true for $Q$.
Therefore, all segments of trajectories of $L$ with endpoints on
$P \cup Q$ have uniformly bounded length (understood as the length
of the corresponding time-interval). This readily yields that $Y$
is closed, and hence $Y=P$, since $P$ is connected. Hence, $Q$ is
a starshaped hypersurface in $M_{*,P} = SP$, and the claim
follows. In particular, the core of $M$ is independent of $P$. We
write $\Core(M)$ for $\Core_P(M)$, and $M_*$ for $M \setminus
\Core (M)$.

Furthermore, the projection $\pi_{P,Q}$ of $Q$ to $P$ along the
trajectories of $L$ establishes a contactomorphism between the
contact structures $\mathrm{Ker}(\alpha|_{TP})$ and
$\mathrm{Ker}(\alpha|_{TQ})$. Thus we associated to a
Liouville manifold $(M,\omega,L)$ a {\it canonical contact manifold},
  defined as a family of contact manifolds
$\{(P,\Ker{\alpha|_{TP}})\}$ and contactomorphisms $\pi_{P,Q}\co  Q
\to P$ satisfying $\pi_{P,Q}\circ \pi_{Q,R} = \pi_{P,R}$ where
$P,Q,R$ run over the set of all starshaped hypersurfaces in $M$.
This contact manifold admits a more geometric description. Note
that $\R_+$ acts freely on $M\setminus\Core(M)$ by the formula
$c*x=L^{\log c}x$ and we have $\left(L^{\log
c}\right)^*\alpha=c\alpha$. Hence the plane field $\{\alpha=0\}$
on $M_*$ is invariant under the $\R_+$--action. It descends to a
contact plane field $\eta_{\infty}$ on $P_{\infty} := M_*/\R_+$.
We will call $(P_{\infty},\eta_{\infty})$ the {\it ideal contact
boundary} of the Liouville manifold, and in the sequel will use
both of its descriptions. Let us emphasize that the
symplectization of $P_{\infty}$ is canonically identified with
$M_*$, and hence we have a canonical decomposition
\begin{equation} \label{eq-decomposition}
M = SP_{\infty} \sqcup \Core(M)\;.
\end{equation}
Contact manifolds $(P_{\infty},\eta_{\infty})$ arising in this way
are called {\it Liouville-fillable}, and we refer to
$(M,\omega,L)$ as a {\it Liouville filling} of
$(P_{\infty},\eta_{\infty})$.

\begin{example}\label{ex-Liouv-mfd-on} {\rm The standard symplectic
linear space $(\R^{2n}, dp \wedge dq)$ equipped with the vector
field $L = \frac{1}{2}(p\frac{\p}{\p p} + q \frac{\p}{\p q})$ is a
Liouville manifold. The Liouville form $\alpha$ equals
$\frac{1}{2}(pdq-qdp)$. It follows from \fullref{SS} above
that the ideal contact boundary of $\R^{2n}$ is the standard
contact sphere $S^{2n-1}$, and the core equals $\{0\}$. Our
convention is that the zero-dimensional space (that is, the point)
is a Liouville symplectic manifold: it coincides with its core,
and its ideal contact boundary is empty.}
\end{example}

\begin{example}\label{ex-Liouv-mfd-tw} {\rm The cotangent bundle
$T^*X$ of a closed manifold $X$ equipped with the standard
symplectic form $dp \wedge dq$ and the Liouville vector field
$p\frac{\p}{\p p}$ is a Liouville manifold. Its ideal contact
boundary is called {\it the space of co-oriented contact elements}
of $X$ and is denoted $\PP_+T^*X$. The core coincides with the
zero section.}
\end{example}

\begin{example}[Weinstein manifolds]\label{ex-Liouv-mfd-th}
{\rm A Liouville manifold $(M,\omega,L)$ is called {\it Weinstein} if
the vector field $L$ is gradient-like for an exhausting (that is,
proper and bounded from below) Morse function $h$ on $M$.  Contact
manifolds $(P,\eta)$ arising as ideal contact boundaries of Weinstein
manifolds are called Weinstein-fillable \cite{Eliashberg??}. This
class of contact manifolds has an alternative description in terms of
complex geometry, see \cite{Eliash-Stein}: namely, they appear as
strictly pseudo-convex boundaries of Stein domains equipped with the
field of complex tangent subspaces of maximal dimension. (Recall that
a Stein domain is a sublevel set of an exhausting plurisubharmonic
function on a Stein manifold.) All critical points of the function $h$
as above have Morse index $\leq n=\frac{1}{2} \dim W$. If all the
indices do not exceed $n-k$ for $1 \leq k \leq n$, the Weinstein
manifold $M$ is called {\it $k$--subcritical.} Otherwise, it is called
{\it critical.} For instance, $\R^{2n}$ is $n$--subcritical (and hence,
of course, $k$--subcritical for every $k \in [1,n]$), while $\PP_+T^*X$
is critical. A Weinstein-fillable contact manifold is called
$k$--subcritical if it admits a $k$--subcritical Weinstein filling, and
critical otherwise.}\end{example}

\begin{example}[Stabilization]\label{ex-Liouv-mfd-fo} {\rm
Let $(M,\omega,L)$, $(M',\omega',L')$ be two Liouville manifolds.
Their product $(M \times M', \omega \oplus \omega', L \oplus L')$
is again a Liouville manifold. In the case when
$$(M',\omega',L')= \left( \R^{2n}, dp \wedge dq, \frac{1}{2} (p\ppp
+q\ppq)\right)\;$$ the obtained Liouville manifold is called the
{\it $n$--stabilization} of $(M,\omega,L)$. The significance of
this notion is due to the following result by K Cieliebak
\cite{Cieliebak}: Every $k$--subcritical Weinstein manifold is a
$k$--stabilization of another Weinstein manifold.}
\end{example}

\subsection{On orderability of Liouville-fillable manifolds} \label{sec-more-orderability}

 The next theorem is a generalization of \fullref{thm-nogo}:

\begin{thm}\label{thmtw-stabilization}
For any Liouville manifold $(M,\omega,L)$ the ideal contact
boundary of its $n$--stabilization is not orderable provided that
$n\geq 2$.
\end{thm}

 The proof is given in \fullref{sectw-stab}
below. \fullref{thm-nogo} corresponds to the case when $M$ is
a point. Our discussion in \fullref{ex-Liouv-mfd-fo} above
yields the following corollary.

\begin{cor}\label{thmsubcrit}
Weinstein-fillable $2$--subcritical contact manifolds are not
orderable.
\end{cor}

It is interesting to confront this result with the following

\begin{thm}\label{thmcont-elements}
Let $X$ be a closed manifold. Assume that either $\pi_1(X)$ is
finite, or $\pi_1(X)$ has infinitely many conjugacy classes. Then
the space $\PP_+T^*X$ of co-oriented contact elements of $X$ with
its canonical contact structure is orderable.
\end{thm}

 This result generalizes Theorem 1.3.B. from \cite{EP}.  We
refer to \fullref{seccotangent} for the proof.  Interestingly enough,
it is unknown whether \fullref{thmcont-elements} covers all closed
manifolds. This depends on the answer to the following long-standing
open question in group theory due to J Makowsky \cite{Makowsky}: does
there exist an infinite finitely-presented group with finitely many
conjugacy classes? (See Baumslag, Myasnikov and Shpilrain \cite{zebra}
for further discussion.) Note that spaces of co-oriented contact
elements, or unit cotangent bundles, are examples of contact manifolds
with critical Weinstein filling.  We refer the reader to
\fullref{sec-discussion-open-problems} for further discussion on
fillability and orderability.

\subsection{Orderability versus squeezing} \label{subsec-cont-dom}

Let $(M,\omega, L)$ be a Liouville manifold with the ideal contact
boundary $(P,\eta)$. Write $\alpha$ for the Liouville 1--form on
$M$. Consider the prequantization space of $M$ which is defined as
the contact manifold $(V = M \times S^1, \xi =
\mathrm{Ker}(\alpha-dt))$. In this section we present a basic link
between contact (non)-squeezing in $V$ and the orderability of
$P$.

\medskip
{\bf Convention}\qua In what follows we identify $SP$ with
$M_* = M \setminus \Core(M) $ and extend all contact Hamiltonians
from $M_*$ to $M$ by setting them equal to $0$ on $\Core(M)$. Of
course, this extension is not, necessarily, smooth along
$\Core(M)$.

\medskip
An open domain of $V$ with smooth compact boundary is
called {\it fiberwise starshaped} if it intersects the fibers $M
\times \{t\},\; t \in S^1$ along starshaped domains and its
boundary is transversal to the fibers.

Suppose that $P$ is non-orderable, that is it admits a positive
contractible loop of contactomorphisms. Denote by $\{f_t\}$ the
corresponding loop of $\R_+$--equivariant symplectomorphisms of
$M_*$. Let $$\{f_{t,s}\},\; t \in S^1, s \in [0,1],$$ be a
homotopy of $\{f_{t,1}\} = \{f_t\}$ to the constant loop
$\{f_{t,0}\} \equiv \id$. Assume that $f_{0,s} = \id$ for all $s$.
Write $F_s\, ,\; s \in[0,1],$ for the contact Hamiltonian on $V$
generating the loop $\{f_{t,s}\},\; t \in S^1$. Let $E\co  M \to \R$
be a positive contact Hamiltonian. For each $R>0$ let us consider
a domain
$$A(R) = \{E < R\} \times S^1 \subset V \;.$$
Roughly speaking, the next theorem shows that every homotopy of a
positive contractible loop of contactomorphisms of $P$ to the
constant loop serves as a ``squeezing tool" for the domain $A(R)$
provided that $R$ is small enough.

\begin{thm} \label{thm-corresp}
Suppose there exists $\mu > 0$ such that
$$F_s(z,t) > -\mu E(z) $$
for all $z \in M_*, t \in S^1, s \in [0;1]$. Then
\begin{itemize}
\item[(i)] there exists
$\gamma >0 $ such that $A(R)$ can be contactly squeezed into
$A\Big{(}\frac{R}{1+\gamma R}\Big{)}$ for all $R < \mu^{-1}$;
\item[(ii)] $A(R)$ can be contactly squeezed into itself
inside $A(\rho)$ for all
$$R < \mu^{-1} \;\;\;\text{and}\;\;\; \rho >
\frac{1}{R^{-1}-\mu}\;.$$
\end{itemize}
Moreover the squeezing in (i) and (ii) can be performed in the
class of fiberwise starshaped domains.
\end{thm}

 The proof is based on an elementary geometric
construction from \cite{EP} which links together positive paths of
contactomorphisms of $P$ and fiberwise starshaped domains in $V$.
Let $\phi=\{f_t\},\; t \in [0;1], f_0 = \id,$ be a path of
contactomorphisms of $P$ generated by a positive time-periodic
contact Hamiltonian $F\co V \to \R$. Define a domain $$U(\phi) :=
\{F< 1\} \subset V\;.$$ Vice versa, every fiberwise starshaped
domain of $V$ is of the form $\{F<1\}$ for a unique positive
time-periodic contact Hamiltonian $F$, and thus corresponds to the
positive path of contactomorphisms of $P$ generated by $F$.

\begin{exam}[Standard rotation]\label{ex-rotat}
{\rm The rotation $e_t(z) = e^{2\pi i t}z, \; t \in \R$, defines a
positive path of contact transformations of $P=S^{2n-1}$ generated
by the contact Hamiltonian $\pi|z|^2$. Fix $R>0$ and consider the
path $\epsilon := \{e_{t/R}\}$. Clearly, the corresponding
starshaped domain $U(\epsilon)$ is simply $$\wh B^{2n}(R)=
{B}^{2n}(R) \times S^1 \subset \C^n \;.$$}
\end{exam}

 A crucial feature of the correspondence $\phi \to
U(\phi)$ is as follows:

\begin{lemma}\label{lmpath-domain}
Let $\{\phi_s\},\; s \in [0;1]$ be a homotopy of paths of
contactomorphisms of $P$ through positive paths with fixed end
points. Then there exists an ambient contact isotopy $\Psi_s\co V\to
V$ with $\Psi_0 = \id$ and $\Psi_s(U(\phi_0))=U(\phi_s)$ for all
$s \in [0;1]$.
\end{lemma}

 The proof virtually repeats the one of Lemma 3.1.B in
\cite{EP}.

\begin{rem}\label{remdomain-path}
{\rm Note, however, that existence of a global contact isotopy
which deforms $U(\phi_0)$ into $U(\phi_1)$ through starshaped
domains implies only that there exists a homotopy of positive
paths of contactomorphisms
$$\phi_s=\{f_s^t\}_{t\in[0,1]},\,s\in[0,1],$$ such
that $f_s^0=\id$ and $f_s^1=g_s^{-1}\circ f_0^1\circ g_s$ for a
certain path of contactomorphisms $g_s\co P\to P,\, s\in[0,1]$. }
\end{rem}

\begin{proof}[Sketch of the proof of \fullref{thm-corresp}
(ii)] Denote by $e_t$ the Hamiltonian flow of $E$ on $M_*$.
Denote by $\epsilon$ the path of contactomorphisms of $P$
corresponding to $e_{t/R},\; t \in [0;1]$. With this notation
$A(R) = U(\epsilon)$.

Consider the following homotopy of the path $\epsilon :=
\{e_{t/R}\}$ with fixed endpoints:
$$\epsilon_s = \{e_{t/R}f_{t,s}\},\; t \in S^1, s \in [0;1].$$
The Hamiltonian $H_s$ generating $\epsilon_s$ is given by
\begin{equation} \label{eq-red}
H_s(z,t) = R^{-1}E(z) + F_s(e_{-t/R}z,t) \;.
\end{equation}
The assumptions of the theorem guarantee that
$$H_s (z,t)> R^{-1}E(z)-\mu E(e_{-t/R}z)\;.$$
The energy conservation law for $\{e_t\}$ yields $E(e_{-t/R}z) =
E(z)$, and therefore $$H_s(z,t) > (R^{-1} -\mu)E(z) >
\rho^{-1}E(z)\;.$$ Hence, by \fullref{lmpath-domain} the family
of closed domains $\Cl(U(\epsilon_s))$ provides a contact isotopy
of $\Cl(U(\epsilon_0)) = \Cl(A(R))$ inside $A(\rho)$. Recall that
we started with a positive loop, hence $F_1
> 0$. Therefore $H_1 > R^{-1}E(z)$ which yields $\Cl (U(\epsilon_1)) \subset
U(\epsilon_0) = A(R)$, and thus we get the desired squeezing of
$A(R)$ into itself inside $A(\rho)$. \end{proof}

 The argument above contains a gap: Hamiltonians
$H_s(z,t)$ are in general not 1--periodic in time, and hence
formally speaking domains $U(\epsilon_s)$ are not well defined.
This can be corrected by an appropriate ``smoothing of corners"
argument: see \fullref{sec-sq-ord} below for the complete
proof of \fullref{thm-corresp}.

\fullref{thm-corresp} applied to the case when $M=\R^{2n}
= \C^n$ and $E(z) =\pi |z|^2$ enables us to reduce \fullref{mu}(i) to \fullref{thm-nonsq} and \fullref{thm-sq} to
\fullref{mu}(ii). Note that in this case $A(R) = B^{2n}(R)$.

\begin{proof}[\fullref{thm-nonsq} $\Rightarrow$ \fullref{mu}(i)] Assume on the contrary that $\mu(\Delta) < 1$. Then
\fullref{thm-corresp}(i) gives us a squeezing of $\Cl
(\wh{B}(1))$ into its interior, in contradiction to \fullref{thm-nonsq}.
\end{proof}

\begin{proof}[\fullref{mu}(ii) $\Rightarrow$ \fullref{thm-sq}] Assume that $0< R_2 < R_1 < 1$. Choose a homotopy
$\Delta = \{f_{t,s}\}$ of a positive contractible loop
$\{f_{t,1}\}$ of contactomorphisms of the sphere to the constant
loop $\{f_{t,0}\}\equiv \id$ with $\mu (\Delta) < 1/R_1$. \fullref{thm-corresp}(i) guarantees that for every $R \leq R_1$, the
domain $\wh{B}^{2n}(R)$ can be squeezed into $\wh{B}^{2n}(v(R))$
where the function $v(R)$ is given by
$$v(R) = \frac{R}{1+\gamma R}$$
for some $\gamma > 0$. Note that $0 < v(R) < R$ for $R>0$, and the
$N$-th iteration $v^{(N)}(R)= (v \circ...\circ v)(R)$ satisfies
$v^{(N)}(R) \to 0$ as $N \to +\infty$. Choose $N$ large enough so
that $v^{(N)}(R_1) < R_2$. Hence, iterating our construction $N$
times we get a squeezing of $\wh{B}^{2n}(R_1)$ into
$\wh{B}^{2n}(R_2)$. This completes the proof. \end{proof}

\begin{rem} \label{rem-sharp}{\rm
The second part of \fullref{mu} states that for every
$\varepsilon>0$ there exists a homotopy $\Delta$ of a positive
contractible loop to the constant loop with $\mu(\Delta) \leq
1+\varepsilon$. Hence, \fullref{thm-nosmallsq} provides a
sharp restriction on any squeezing of $\wh{B}^{2n}(\frac{1}{k})$
into itself, where $k>1$ is an integer. Indeed, apply \fullref{thm-corresp}(ii) with $R = \frac{1}{k}$ and $\mu(\Delta) \leq
1+\varepsilon$. We get a squeezing of $\wh{B}^{2n}(\frac{1}{k})$ into
itself inside $\wh{B}^{2n}(\rho)$ provided $\rho > 1/(k-\mu)$.
Taking $\varepsilon$ small enough we see that $\rho $ can be
chosen arbitrarily close to $\frac{1}{k-1}$.}
\end{rem}

Let us mention also that there is an alternative way to use a
positive (not necessarily contractible) loop of contactomorphisms
of the sphere $S^{2n-1}$ for producing contact embeddings of
domains in $\R^{2n} \times S^1$, see \fullref{subseq-loopc}
below. A special (trivial) case of this construction enables us to
prove the following result which we already have mentioned in
\fullref{subsec-sq} above.

\begin{prop} \label{thm-all-sq}
For all $R_1,R_2 > 0$ there exists a contact embedding of
$\wh{B}^{2n}(R_1)$ into $\wh{B}^{2n}(R_2)$. For $n>1$ this
embedding can be chosen isotopic through smooth embeddings to the
natural inclusion $$\wh{B}^{2n}(R_1)\hookrightarrow \R^{2n}\times
S^1\;.$$
\end{prop}

 Of course, in general this isotopy cannot be made
contact. Furthermore, assume that $R_2 \leq m < R_1$ for some $m
\in \N$. Then the embedding guaranteed by the proposition does not
extend to a compactly supported contactomorphism of $V$ in view of
the Non-Squeezing Theorem above.

\begin{proof}
For any integer $N>0$ define a map $F_N\co \R^{2n}\times
S^1\to\R^{2n}\times S^1$ by the formula
$$(z,t)\mapsto\left(v(z)e^{2\pi Ni t}z,t\right),$$
where
\begin{equation}\label{eqsimple-twist}
v(z)=\frac1{\sqrt{1+N\pi\sum\limits_1^n|z_j|^2}}\;.
\end{equation}
It is straightforward to check that $F_N$ is a contactomorphism of
$\R^{2n}\times S^1$. It maps $\wh{B}^{2n}(R_1)$ onto
$\wh{B}^{2n}(R)$ with
$R=\frac{R_1}{1+NR_1}\mathop{\to}\limits_{N\to\infty} 0$. For
$n>1$ and an even $N$ the map $F_N$ is isotopic to the inclusion.
\end{proof}

The following fact is well known to specialists. However, we provide it with a proof for a reader's convenience.
\begin{cor}\label{cor-small-nbhd}
For every $R>0$ there exists a contact embedding of $\hat{B}(R)$
into an arbitrarily small neighbourhood of a point in
any contact manifold.
\end{cor}

\begin{proof} By the contact Darboux theorem it suffices
to show that for every $R>0, \epsilon>0$ there exists a contact
embedding of $\hat{B}(R)$ into Euclidean ball of radius $\epsilon$
in the standard contact space
$$(\R^{2n+1}, \Ker(dz+\frac{1}{2}(pdq-qdp))\;.$$
Consider the Euclidean circle of radius $\epsilon/2$ with the
center at $0$ lying in the $(p_1,q_1)$--plane. This circle is
transversal to the contact structure, and hence by the relative
Darboux theorem (eg see \cite{Ge}, Example 2.33) its sufficiently
small neighbourhood is contactomorphic to $\hat{B}(r)$ with some
$r>0$. But $\hat{B}(R)$ can be contactly embedded into
$\hat{B}(r)$ in view of \fullref{thm-all-sq}. This
completes the proof.
\end{proof}

\subsection{Contact homology} \label{subsec-ch}

The proof of Non-Squeezing \fullref{thm-nonsq} is based on
contact homology theory. Here we present a brief outline of
contact homology adjusted to our purposes and refer the reader to
\fullref{secflavors} below for more details and
generalizations.\footnote{The more general contact homology groups
constructed in \fullref{secflavors} have a more sophisticated
grading which involves the set of free homotopy classes of loops
in $M$. The simpler version presented in this section corresponds
to the component associated to the class of contractible loops.}
Let $(M,\omega,L)$ be a Liouville manifold whose first Chern class
$c_1(TM)$ vanishes on $\pi_2(M)$. Put $V = M \times S^1$ and set
$\lambda = dt-\alpha$, where $\alpha$ is a Liouville form on $M$.
Equip $V$ with the contact structure $\xi = \mathrm{Ker}(\lambda)$.

Consider the set $\cU$ consisting of all fiberwise starshaped open
domains in $V$ and their images under the group $\cG = \mathrm{Cont}
(V,\xi)$ of compactly supported contactomorphisms of $(V,\xi)$.
Contact homology $\CH(U)$ of a domain $U \in \cU$ is a $\Z$--graded
vector space over $\Z_2$. Every inclusion $U_1 \subset U_2$ gives
rise to a morphism
$$\incl \co  \CH(U_1) \to \CH(U_2).$$ Every contactomorphism
$\Phi \in \cG$ induces an isomorphism
$$\Phi_{\sharp}\co  \CH(U) \to \CH(\Phi(U)).$$ These morphisms preserve the grading and
have the following properties (the diagrams below are
commutative):
\begin{equation} \label{eq-functon}
\xymatrix{ \CH(U_1) \ar[r]^{\incl} \ar[dr]^{\incl} &
\CH(U_2) \ar[d]^{\incl} \\
&\CH(U_3) }
\end{equation}
provided $U_1 \subset U_2 \subset U_3$ are domains from $\cU$;
\begin{equation} \label{eq-functtw}
\xymatrix{ \CH(U_1) \ar[r]^{\incl} \ar[d]^{\Phi_{\sharp}} &
\CH(U_2) \ar[d]^{\Phi_{\sharp}} \\ \CH(\Phi(U_1)) \ar[r]^{\incl} &
\CH(\Phi(U_2)) }
\end{equation}
provided $U_1 \subset U_2$ are domains from $\cU$;
\begin{equation} \label{eq-functth}
(\Phi\circ \Psi)_{\sharp}= \Phi_{\sharp}\circ\Psi_{\sharp};
\end{equation}
for all $\Phi,\Psi \in \cG$;
\begin{equation} \label{eq-functfo}
\incl_U = \id \;\;\; \text{and}\;\;\; \id_{\sharp} = \id\;,
\end{equation}
where $\incl_U$ is the natural inclusion $U\to U$.

 In a more formal language, $\CH$ is a $\cG$--functor from the
category, whose objects are domains from $\cU$ and morphisms
correspond to inclusions, to the category of $\Z$--graded vector spaces
over $\Z_2$.\footnote{Let $\cU$ be a category, and a group $\cG$ acts
on $\cU$ by functors. A $\cG$--functor is a functor $F$ from $\cU$ to
another category, and a family of natural transformations $g_{\sharp}\co 
F \to F \circ g,\; g \in \cG$ such that $(gh)_{\sharp} = g_{\sharp}
\circ h_{\sharp}$ for all $g,h \in \cG$. The terminology is borrowed
from Jackowski and S{\l}omi{\'n}ska \cite{JS}.}

A typical application of contact homology to non-squeezing
is given in the following statement:

\begin{prop} \label{prop-applch} 
Let $U_1,U_2$ and $W$ be domains from $\cU$ such that $\Cl
(U_1) \subset W$ and $\Cl (U_2) \subset W$. Assume that that the
inclusion morphism $\CH_k(U_1) \to \CH_k(W)$ does not vanish,
while the inclusion morphism $\CH_k(U_2) \to \CH_k(W)$ vanishes
for some $k \in \Z$. Then $U_1$ cannot be mapped into $U_2$ by a
contactomorphism $\Phi \in \cG$ with $\Phi(W) = W$.
\end{prop}

\begin{proof} Assume on the contrary that there exists
a contactomorphism $\Phi \in \cG$ such that $\Phi(W) = W$ and
$\Phi(U_1) \subset U_2$. Then the diagrams \eqref{eq-functon} and
\eqref{eq-functtw} yield the following commutative diagram:
\[
\xymatrix{ \CH_k(U_1) \ar[r]^{\incl\neq 0}
\ar[d]^{\Phi_{\sharp}} & \CH_k(W) \ar[dr]^{\Phi_{\sharp}} & \\
\CH_k(\Phi(U_1)) \ar[r]^{\incl} & \CH_k(U_2)\ar[r]^{\incl = 0} &
\CH_k(W) }
\]
But $\Phi_{\sharp}$ is an isomorphism, and we get a contradiction
which proves the desired non-squeezing.
\end{proof}

Our next result relates contact homology of a split domain
$\wh D = D \times S^1$, where $D$ is a starshaped domain of $M$,
with filtered symplectic homology $\SH_*(D)$. We use a version of
symplectic homology associated to negative compactly supported
Hamiltonians as in \cite{BPS,CGK,GG}, see
\fullref{secflavors} below for the definition and conventions
about the grading.

\begin{thm} \label{thm-chsh}
Assume that the characteristic foliation of $\partial D$ has no
closed orbits $\gamma$ with $\int_{\gamma}\alpha = 1$. Then {\rm
$$ \CH_*(\wh{D}) = \SH^{(-\infty;-1)}_* (D).$$} Moreover, the
above correspondence between symplectic and contact homology is
functorial in the sense that it commutes with the morphisms
induced by inclusions and symplectic/contact diffeomorphisms of
domains.
\end{thm}

 We restate and prove this result in \fullref{subseq-calcul-split} below (see \fullref{thm-calcul-split} and \fullref{propfunctoriality}).
\fullref{thm-chsh} enables us to calculate contact homology
for domains of the form $\wh{E}(N,R)$, where $E(N,R)$ is an
ellipsoid given by
\begin{equation} \label{eq-ellip}
E(N,R):= \{\pi |z_1|^2 + \frac{\pi}{N} \sum_{i=2}^n |z_i|^2 < R\},
\; N \in \N, \frac{1}{R} \notin \N
\end{equation}
Put
\begin{equation} \label{eq-k}
k(N,R) = -2\Big{
[}\frac{1}{R}\Big{]}-2(n-1)\Big{[}\frac{1}{NR}\Big{]},
\end{equation}
where $[s]$ stands for the integer part of a positive real number
$s$.

\begin{thm} \label{thm-ch-ellip}
{\rm $\CH_i(\wh{E}(N,R)) = \Z_2$} for $i = k(N,R)$ and {\rm
$\CH_i(\wh{E}(N,R)) = 0$} for all $i \neq k(N,R)$.
\end{thm}

Using the correspondence between contact and symplectic
homologies, it is easy to calculate the morphisms induced by
natural inclusions of balls (that is $N=1$):

\begin{thm} \label{thm-ball-morph}
Assume that $$\frac{1}{k} < R_1 <R_2 < \frac{1}{k-1}$$ for some $k
\in \N$. Then the inclusion $B^{2n}(R_1) \hookrightarrow
B^{2n}(R_2)$ induces an isomorphism of contact homologies.
\end{thm}

 Theorems \ref{thm-ch-ellip} and \ref{thm-ball-morph} are
proved in \fullref{subsec-ellip} below. Now we are ready to
prove Non-Squeezing \fullref{thm-nonsq}.

\begin{proof}[ Proof of \fullref{thm-nonsq}] Without loss of
generality assume that $R_2 < m < R_1$ for some $m \in \N$.

First of all we reduce the problem to the case $m=1$ as follows.
Assume on the contrary that there exists a contactomorphism $\Phi
\in \cG$ so that
$$\Phi(\mathrm{Closure}(\wh{B}^{2n}(R_1))) \subset
\wh{C}^{2n}(R_2)\;.$$ Define the $m$--fold covering
$$\tau\co  \R^{2n} \times S^1 \to \R^{2n} \times S^1, \;\; (z,t)
\mapsto (\sqrt{m}z,mt)\;.$$ Observe that $\tau^*(dt-\alpha) =
m(dt-\alpha)$, and therefore
\begin{equation} \label{eq-lift}
\tau_*\xi=\xi\;.
\end{equation}
Furthermore,
$$\tau^{-1}(\wh{B}^{2n}(R_1))=
\wh{B}^{2n}(R_1/m)\;\;\;\text{and}\;\;\;
\tau^{-1}(\wh{C}^{2n}(R_2)) = \wh{C}^{2n}(R_2/m)\;.$$
Recall that the diffeomorphism $\Phi$ is compactly supported. In
particular, $\Phi$ acts trivially on the fundamental group of
$\R^{2n} \times S^1$, and the complement to the support is
connected. Thus $\Phi$ admits a unique lift (with respect to $\tau$)
to a compactly supported diffeomorphism of the covering space.
  In view of \eqref{eq-lift} this lift, denoted $\Phi'$, is
a contactomorphism.

Summing up, passing to the covering space, we get that
$$\Phi'(\mathrm{Closure}(\wh{B}^{2n}(R_1/m))) \subset
\wh{C}^{2n}(R_2/m)\;$$ where $\Phi' \in \cG$ and $R_2/m < 1 <
R_1/m$. This completes the reduction to the case $m=1$.

In view of the reduction above we assume that $R_2 < 1 < R_1$.
Furthermore, let us suppose that $1/R_2 \notin \Z$. It suffices to
show that $\wh{B}^{2n}(R_1)$ cannot be mapped by a
contactomorphism $\Phi \in \cG$ into domain $\wh{E}(N,R_2)$ with
$N$ large enough, which corresponds to a ``long" ellipsoid, so
that $\Phi$ keeps invariant the domain $\wh{B}^{2n}(R_3)$ with
large $R_3$. Using Theorems \ref{thm-ch-ellip} and
\ref{thm-ball-morph} we see that $\CH_0 (\wh{B}^{2n}(R_1)) =
\CH_0(\wh{B}^{2n}(R_3)) = \Z_2$ and, moreover, the natural
inclusion induces an isomorphism in contact homology. Here we used
that $R_3
> R_1 > 1$. On the other hand $\CH_0(\wh{E}(N,R_2)) = 0$ since
$R_2 < 1$. The desired result follows now from \fullref{prop-applch}.
\end{proof}

\begin{proof}[ Proof of \fullref{thm-nosmallsq}]
Exactly as in the previous proof, we reduce the problem to the
case $m=1$. Assume that
$$\frac{1}{k+1}< R_2 < \frac{1}{k} < R_1 < R_3 < \frac{1}{k-1}.$$
We have to show that $\wh{B}^{2n}(R_1)$ cannot be mapped into
$\wh{B}^{2n}(R_2)$ by a contactomorphism from $\cG$ which keeps
invariant the domain $\wh{B}^{2n}(R_3)$. Put $N= -2n(k-1)$. Using
Theorems \ref{thm-ch-ellip} and \ref{thm-ball-morph} we see that
$$\CH_N (\wh{B}^{2n}(R_1)) = \CH_N(\wh{B}^{2n}(R_3)) = \Z_2,$$ and
moreover, the natural inclusion induces an isomorphism in contact
homology. On the other hand $\CH_N(\wh{B}^{2n}(R_2)) = 0$. The
desired result follows now from \fullref{prop-applch}.
\end{proof}

\subsection{A guide for the reader}

In the next diagram we present logical interrelations between the
main results, methods and phenomena described above. We abbreviate
PC for {\it positive contractible} and CH for {\it contact
homology}. The arrows are labeled by the numbers of the
corresponding theorems and/or propositions.

\[
\xymatrix{ \fbox{SQUEEZING} & \fbox{$\exists$
PC-LOOP}\ar[l]_{\ref{thm-corresp}}\ar[dr]\ar[r]^{\ref{prop-nonderability}}
& \fbox{NO ORDER} \ar[l] \\
\fbox{CH}\;\;\;\;\ar[r]^-{\ref{prop-applch}}&\;\;\;
\fbox{NON-SQUEEZING}\;\;\; \ar[r]^-{\ref{mu}} \ar[d] &
\fbox{PC-LOOPS RIGIDITY}\\ & \fbox{ORDER} }
\]

 The rest of the paper is organized as follows:

In \fullref{sec-sq-ord} we continue our discussion started in
\fullref{subsec-cont-dom} where we suggest to use loops of
contactomorphisms as a squeezing tool for contact domains. In
particular, we prove \fullref{thm-corresp} and generalize (and
demystify) the formulas which appeared in the proof of \fullref{thm-all-sq}.

In \fullref{sectw-stab} we prove \fullref{thmtw-stabilization} which states that the ideal contact
boundary of the $n$--stabilization of any Liouville manifold is
non-orderable, provided that $n\geq 2$. To do this, we present a
rather explicit construction of a positive contractible loop of
contactomorphisms. Analyzing an explicit homotopy of our positive
loop to the constant loop, we prove \fullref{mu}(ii). Finally,
we prove \fullref{thm-sq-strong} on negligibility of the
cylinder.

To prove our results on contact non-squeezing and orderability we
shall need to understand the relation between cylindrical contact
homology, which is a special case of the symplectic field theory,
and Hamiltonian Floer homology. In \fullref{secflavors} we
introduce a version of Floer homology which contains both theories
as its special cases.

In \fullref{sec-calcul} we complete the proof of several
``hard" results stated in \fullref{sec-intro}. First of all,
we express contact homology of the prequantization of a symplectic
domain in terms of its filtrated symplectic homology, see \fullref{thm-calcul-split} and \fullref{propfunctoriality}
below which together form a slightly more explicit version of
\fullref{thm-chsh}. Our approach is based on generalized Floer
homology theory developed in \fullref{secflavors}. We apply
this result to various calculations with contact homology. In
particular, we prove Theorems \ref{thm-ch-ellip} and
\ref{thm-ball-morph} on contact homology of prequantizations of
ellipsoids and balls in $\R^{2n}$. Furthermore, we study contact
homology of prequantizations of unit ball bundles of closed
manifolds in terms of cohomology of free loop spaces. As a result,
we prove \fullref{thmcont-elements} on orderability of spaces
of contact elements.

In \fullref{sec-discussion-open-problems} we touch
miscellaneous topics related to the geometry of contact domains
and transformations: We discuss (non)-squeezing of contact domains
from the viewpoint of quantum mechanics. Then we speculate on
links between orderability and fillability of contact manifolds.
Finally, we introduce a canonical semigroup associated to a
contact manifold. This semigroup carries footprints of a partial
order even when the manifold is non-orderable.

\section{Loops of contactomorphisms as a squeezing tool} \label{sec-sq-ord}

Here we prove \fullref{thm-corresp} by filling the gaps in the
sketch presented in \fullref{subsec-cont-dom} above, and
generalize the formulas which appeared in the proof of \fullref{thm-all-sq}.

\subsection[Proof of \ref{thm-corresp}]{Proof of \fullref{thm-corresp}}

Denote by $e_t$ the Hamiltonian flow of $E$ on $M_*$. Denote by
$\epsilon$ the path of contactomorphisms of $P$ corresponding to
$e_{t/R},\; t \in [0;1]$. With this notation we have $A(R) =
U(\epsilon)$. Fix $\delta>0$ small enough so that
\begin{equation} \label{eq-vsp-corresp-on}
F_1(z,t) > \delta R^{-1} E(z) \;\;\;\text{and}\;\;\; F_s(z,t) >
-(1-\delta) \mu E(z)
\end{equation}
for all $z \in M_*, t \in S^1, s \in [0,1]$. Choose a
non-decreasing function $\tau \co [0,1] \to [0,1]$ such that $\tau
\equiv 0$ near $0$, $\tau \equiv 1$ near $1$ and $\tau'(t) <
1+\delta$ for all $t \in [0,1]$. Put
$$\theta(t)=(1+\delta)t-\delta\tau(t)\;.$$
Note that
\begin{equation}\label{eq-vsp-corresp-tw}
\theta'(t) \equiv
1+\delta\;\;\;\text{for}\;\;\;t\;\;\;\text{near}\;\;\;0\;\;\;\text{and}\;\;\;1\;,\;\;\text{and}\;\;\;
\theta'(t)\geq 1-\delta^2
\end{equation}
for all $t$.

Let us first prove part (ii) of the theorem. Our strategy is as
follows. We describe a homotopy $\epsilon_s ,\; s \in [-1,1]$, of
the path $\epsilon = \epsilon_{-1}$ with fixed end-points such
that
\begin{equation}\label{eq-vsp-corresp-th}
U(\epsilon_s) \subset A(\rho)
\end{equation}
for all $s$ and $\Cl (U(\epsilon_1)) \subset A(R)$. The special
attention is paid to the fact that the paths $\epsilon_s$ have to
be generated by 1--periodic in time Hamiltonians: the functions
$\tau$ and $\theta$ are designed just to guarantee
time-periodicity.

Each $\epsilon_s$ corresponds to a path $\{h_{s,t}\}$ of
$\R_+$--equivariant symplectomorphisms of $M_*$ generated by
contact Hamiltonians $H_s(z,t)$. These symplectomorphisms are
defined as follows.

\eject
{\bf Step 1\qua $s \in [-1,0]$}

This part of our homotopy consists of time-reparameterizations of
$\{e_{t/R}\}$. Put
$$\kappa(s,t) = (-st + (s+1)\theta(t))R^{-1}\;,$$
and set $h_{s,t} = e_{\kappa(s,t)}$. Then
$$H_s(z,t) = (-s + (s+1)\theta'(t))R^{-1}E(z)\;.$$
It follows from formula \eqref{eq-vsp-corresp-tw} that $H_s(z,t)$
descends to $M_*\times S^1$ and inclusion \eqref{eq-vsp-corresp-th}
holds provided
\begin{equation} \label{eq-vsp-corresp-delta-on}
\delta^2 < 1-\frac{R}{\rho}\;.
\end{equation}

{\bf Step 2\qua $s \in [0,1]$}

Put $h_{s,t}= e_{\theta(t)/R}\circ f_{\tau(t),s}\;.$ The
corresponding Hamiltonians are given by
$$H_s(z,t) = \theta'(t)R^{-1}E(z) + \tau'(t)
F_s(e_{-\theta(t)/R}z,\tau(t))\;.$$ Taking into account that
$\tau'(t)\equiv 0$ and $\theta'(t) \equiv 1+\delta$ for $t$ near
$0$ and $1$ we get that $H_s$ descend to Hamiltonians on $M_*
\times S^1$. Combining the first inequality in
\eqref{eq-vsp-corresp-on} with the energy conservation law for
$\{e_t\}$ we get that
$$H_1(z,t) > R^{-1}\theta'(t)E(z)+ \tau'(t)R^{-1}\delta\cdot
E(e_{-\theta(t)/R}z)$$ $$ = R^{-1}E(z)(\theta'(t) +
\delta\tau'(t))= R^{-1}(1+\delta)E(z)\;.$$ Thus
\begin{equation}\label{eq-vsp-corresp-incl}
U(\epsilon_1) \subset A\Big{(}\frac{R}{1+\delta}\Big{)}\;.
\end{equation}
It remains to verify inclusion \eqref{eq-vsp-corresp-th}. Combining
the second inequality in \eqref{eq-vsp-corresp-on} with the energy
conservation law for $\{e_t\}$ we get that
$$H_s(z,t) > (R^{-1}\theta'(t) -\mu(1-\delta)\tau'(t)) E(z)\;.$$
But $\tau'(t) < 1+\delta$ by assumption and $\theta'(t) \geq
1-\delta^2$ due to formula \eqref{eq-vsp-corresp-tw}, and hence
$$H_s(z,t) \geq (R^{-1} -\mu)(1-\delta^2) E(z)\;.$$
This yields inclusion \eqref{eq-vsp-corresp-th} provided that
\begin{equation}\label{eq-vsp-corresp-delta-tw}
\delta^2 < 1- \frac{1}{\rho(R^{-1}-\mu)}\;.
\end{equation}
This completes the proof of part (ii) of the theorem.

In order to prove part (i) we put $\rho=+\infty$ and note that the
restrictions \eqref{eq-vsp-corresp-delta-on} and
\eqref{eq-vsp-corresp-delta-tw} on $\delta$ read $\delta < 1$.
Therefore, the choice of $\delta$ is governed by inequalities
\eqref{eq-vsp-corresp-on}. Assume that $$F_1(z,t) > c E(z)
\;\;\;\text{and}\;\;\; F_s(z,t) > -(1-c) \mu E(z)\;$$ for some $c
\in (0;1)$. Then the arguments above will work with $\delta =
\min(cR,c)$. Put $\gamma = \min(c,c\mu)$. Since $R\mu < 1$ we have
$\delta \geq \gamma R$. Thus inclusion \eqref{eq-vsp-corresp-incl}
yields $$U(\epsilon_1) \subset A\Big{(}\frac{R}{1+\gamma R
}\Big{)}\;,$$ which proves part (i) of the theorem. \qed

\subsection[Loops in Cont0(P) as contactomorphisms of SPxS1$]{Loops in $\Cont_0(P)$ as contactomorphisms of $SP\times S^1$}
\label{subseq-loopc} Let $(P,\eta)$ be a contact manifold. Write
$\alpha$ and $L$ for the Liouville form and the Liouville vector
field on $SP$, respectively. Consider a contact manifold
$(SP\times S^1, \xi)$ where $\xi = \mathrm{Ker}(dt-\alpha)$. In this
section we study a general construction which relates loops in
$\Cont_0(P)$ and contact embeddings of domains in $SP\times S^1$.
This construction is behind the explicit formula in the proof
\fullref{thm-all-sq} of the Introduction. Let $\{h_t\},\;
t \in S^1,\; h_0 = \id$ be any loop of equivariant
symplectomorphisms of $SP$ generated by a contact Hamiltonian
$H\co SP \times S^1 \to \R$. In what follows the expression of the
form $z/c$ with $z \in SP$ and $c \in \R_+$ is understood in the
sense of the canonical $\R_+$--action on $SP$.

\begin{prop}\label{prop-loop-cont}
The map
$$\Psi(z,t) = \Big{(} \frac{h_t z}{1+ H(h_t z,t)},t
\Big{)}$$ is a contact embedding of the domain $\{1+ H(h_t z)>0\}
\subset SP\times S^1$ into $SP \times S^1$.
\end{prop}

\begin{proof} Recall that for any contact Hamiltonian $H$
on $SP$ we have $$\alpha(\sgrad H) = H.$$ This is simply the Euler
formula for homogeneous functions. Furthermore, $\alpha(L) = 0.$
We denote by $A_c \co SP \to SP$ the action of $c \in \R_+$ on $SP$.
Note that $A_c^*\alpha = c\alpha$. Put $$c =(1+H(h_tz,t))^{-1}.$$
For every tangent vector $v \in T_z(SP)$
$$\Psi_* v = \mathrm{const}\cdot L + A_{c*}
h_{t*} v,
$$
hence
\begin{equation} \label{eq-vsp-loop-on}
(dt-\alpha)(\Psi_*v) = c(dt- \alpha)(v)\;.
\end{equation}
Furthermore,
$$\Psi_* \frac{\partial}{\partial t} = \mathrm{const} \cdot L +c\cdot\sgrad
H_t+\frac{\partial}{\partial t},$$ where $H_t(z):= H(z,t)$. Hence
\begin{equation} \label{eq-vsp-loop-tw}
\begin{split}
(dt-\alpha) (\Psi_* \frac{\partial}{\partial t})& =1 -
c\alpha(\sgrad H_t)\\&= 1-cH(h_tz,t)\\
& = c(dt-\alpha) (\frac{\partial}{\partial t})\;.\\
\end{split}
\end{equation}
Formulas \eqref{eq-vsp-loop-on} and \eqref{eq-vsp-loop-tw} yield
$\Psi^*(dt -\alpha) = c(dt-\alpha)$ which means that $\Psi$ is a
contactomorphism on its domain of definition. \end{proof}

\begin{example}\label{ex-extension}
{\rm When $(P,\eta)$ is the standard contact sphere $S^{2n-1}$ and
$h_t(z)=e^{2\pi iNt}$ then in view of our identification of $SP$
with $\C^n\setminus 0$ the map $\Psi$ from \fullref{prop-loop-cont} takes the form
$$\Psi_N(z,t) = \Big{(}\frac{e^{2\pi i N t}z}{\sqrt{1+\pi N|z|^2}},t\Big{)}\;,$$
which coincides with the formula \eqref{eqsimple-twist} in the
proof of \fullref{thm-all-sq}. }
\end{example}

\section{Non-orderability of stabilizations}\label{sectw-stab}

In this section we prove \fullref{thmtw-stabilization} which
states that the ideal contact boundary $P_{stab}$ of the
$n$--stabilization $M_{stab} = M \times \C^n$ of any Liouville
manifold $(M,\omega,L)$ is non-orderable if $n\geq 2$. To do this,
we present a rather explicit construction of a positive
contractible loop of contactomorphisms of $P_{stab}$. The desired
result then follows from \fullref{prop-nonderability}.
When $M$ is a point, $P_{stab}$ is simply the standard
$(2n-1)$--dimensional contact sphere. Analyzing an explicit
homotopy of our positive loop to the constant loop, we prove
\fullref{mu}(ii). Finally, we combine our squeezing techniques
with the construction from \fullref{subseq-loopc} and get
\fullref{thm-sq-strong} on negligibility of the cylinder.

\subsection{An ambient isotopy}

We start with a version of the ambient isotopy theorem in contact
geometry which will be used below. Let $(P,\eta)$ be a closed
contact manifold with a co-oriented contact structure, and let
$\Sigma \subset SP$ be a starshaped hypersurface.

Consider a Hamiltonian isotopy $h_t\co  W \to SP, t \in [0,1],\; h_0 =
\id$, of a neighborhood $W$ of $\Sigma$ in $SP$. We denote by
$H_t$ the Hamiltonian function generating $h_t$, by $\sigma_t$ the
characteristic foliation on $\Sigma_t=h_t(\Sigma)$, and write
$\sgrad H$ for the Hamiltonian vector field of a function $H$ on
$SP$.

\begin{lemma} [Ambient isotopy] \label{lem-ambisot}
Assume that hypersurfaces $\Sigma_t$ are starshaped for all $t$.
Then there exists a path $\{f_t\}, \; t \in [0,1]$, of
$\R_+$--equivariant symplectomorphisms of $SP$ with $f_0 = \id$
such that $f_t(\Sigma) = \Sigma_t$ and $\sgrad F_t-\sgrad H_t$ is
tangent to $\sigma_t$, where we denote by $F_t$ the Hamiltonian
which generates the path $f_t$ .
\end{lemma}

Let $\beta$ be the contact form on $SP$ whose graph equals
$\Sigma$. In what follows we need a formula for the Hamiltonian
vector field of a function $H$ in canonical local coordinates
$(x,u) \in P \times \R_+$ associated to $\beta$, see \fullref{subsec-geom-ze} above. We assume that $H = H(x,u)$ is defined
in a neighborhood of $\Sigma$, that is in the set $\{x \in P,
|u-1| < \delta\}$. For every $x \in P$ and $u \in
(1-\delta;1+\delta)$ there exists a unique vector $\sgrad _{\beta}
H \in \eta(x)$ such that $d\beta(v,\sgrad _{\beta} H)
=-dH|_{(x,u)}(v)$ for all $v \in \eta(x)$. Let $R$ be the Reeb
vector field of $\beta$. A straightforward calculation shows that
\begin{equation}\label{sgrad}
\sgrad H = -dH (R) \frac{\partial}{\partial u} + u^{-1}\sgrad
_{\beta} H + \frac{\partial H}{\partial u} R\;\;.
\end{equation}

\begin{proof} [ Proof of \fullref{lem-ambisot}] Let $H_t\co  h_t(W)
\to \R$ be the Hamiltonian function generating $h_t$. Denote by
$F_t$ the $\R_+$--equivariant Hamiltonian which equals $H_t$ on
$\Sigma_t$. It generates a path $f_t$ of equivariant
symplectomorphisms. Note that, for fixed $t$, in canonical
coordinates near $\Sigma_t$ we have $F_t(x,u) = uH_t(x,1)$. Using
the fact that $\Sigma_t = \{u=1\}$ and formula \eqref{sgrad} we
get that at each point of $\Sigma_t$ the difference $\sgrad F_t -
\sgrad H_t$ is tangent to $\Sigma_t$. Since this is true for all
$t \in [0,1]$, we conclude that $f_t(\Sigma) = \Sigma_t$. Note
that the two first terms in formula \eqref{sgrad} depend only on
the restriction of the Hamiltonian to the hypersurface and,
therefore, these terms coincide for $F_t$ and $H_t$. On the other
hand, the third term is tangent to $\sigma_t$, and hence so is
$\sgrad F_t-\sgrad H_t$.
\end{proof}

\subsection{Fundamental Lemma on unitary transformations}

In this section we consider unitary transformations of $\C^n$ as
symplectomorphisms. In particular, every path in $U(n)$ is
generated by unique quadratic Hamiltonian. We represent a point $z
\in \C^n$ as $z=(z_1,...,z_n) = (z_1,w)$, where $w =
(z_2,...,z_n)$, and introduce radial coordinates $\rho_j = \pi
|z_j|^2$ for $j=1,...,n$ and $\varrho= \rho_2+...+\rho_n$.

The following loop, which lies in $SU(n)$ (and hence is necessarily contractible)
is an important character of our story:
\begin{equation} \label{eq-loop-Phi}
f_t\co  (z_1,w) \to (e^{2\pi i (n-1)t }z_1, e^{-2\pi i t}w)\;.
\end{equation}

\begin{lemma}[Fundamental Lemma] \label{lem-fund}
The loop $\{f_t\}$ given by formula \eqref{eq-loop-Phi} admits a
homotopy $\{f^{(s)}_t\}$, $f^{(1)}_t = f_t, f^{(0)}_t \equiv \id$
through loops on $SU(n)$ based at $\id$ to the constant loop so
that the corresponding Hamiltonians $F^{(s)}(z,t)$ satisfy
\begin{equation} \label{eq-lem-Phi-on}
F^{(s)}(f^{(s)}_tz,t) \geq -\varrho \;
\end{equation}
for all $s \in [0;1], t\in S^1, z \in \C^n \setminus\{0\}$.
\end{lemma}

\begin{proof} Denote by $b_{j,t}$ the unitary
transformation which preserves all the coordinates besides the
$j$-th one, and multiplies the $j$-th coordinate by $e^{2\pi i
t}$. Let $I^{(s)}_j, \; j \geq 2, s\in [0,\pi/2]$, be the unitary
transformation which preserves all the coordinates besides the
first and the $j$-th ones, and in the $(z_1,z_j)$--plane its action
is given by
$$(z_1,z_j) \mapsto (\cos s\cdot z_1 -\sin s\cdot z_j\;,\; \sin s\cdot z_1 + \cos
s\cdot z_j)\;.$$ $$h^{(s)}_{j,t} = I^{(s)}_j
b_{j,t}(I^{(s)}_j)^{-1}b_{j,-t}\;\leqno{\hbox{Put}} $$ 
$$f^{(s)}_{m,t}z =h^{(s)}_{2,t} \circ ... \circ h^{(s)}_{m,t}\;,
\;\; m = 2,...,n\;.\leqno{\hbox{and}}$$ For $m = n$ this family of loops gives a
homotopy of $\{f_t\}$ to the constant loop. To verify that this
homotopy satisfies inequality \eqref{eq-lem-Phi-on} we prove the
following more general inequality for every $m=2,...,n$:

{\bf Claim}\qua The Hamiltonian $F^{(s)}_m$ generating the
loop $f^{(s)}_{m,t}$ satisfies
\begin{equation} \label{eq-lem-Phi-tw}
F^{(s)}_m(f^{(s)}_{m,t}z,t) \geq -(\rho_2+...+\rho_m)\;.
\end{equation}

 In the proof of the claim we use the following fact:

{\sl Suppose that Hamiltonian $G(z,t)$ generates a flow $g_t$.
Then the Hamiltonian\break $-G(g_tz,t)$ generates the flow $g_t^{-1}$. }

We use notation $(u)_j$ for the $j$-th coordinate of a vector $u
\in \C^n$. Write $H^{(s)}_j(z,t)$ for the Hamiltonian generating
$h^{(s)}_{j,t}$. Calculating, we get that
\begin{equation} \label{eq-lem-Phi-th}
H^{(s)}_j(h^{(s)}_{j,t}z,t) = - \pi |z_j|^2 + \pi
\Big{|}\Big{(}(I_j^{(s)})^{-1}b_{j,-t}z\Big{)}_j\Big{|}^2 \geq
-\pi |z_j|^2 \;.
\end{equation}
Another observation used below is that
\begin{equation} \label{eq-lem-Phi-fo}
(h^{(s)}_{j,t}z)_l = z_l \;\;\text{for}\;\; l \neq1,j\;.
\end{equation}

\begin{proof}[Proof of the Claim] We use induction in $m$. The case $m=2$
immediately follows from inequality \eqref{eq-lem-Phi-th}. Let us
verify the induction step $m \mapsto m+1$. Note that
$f^{(s)}_{m+1,t} = f^{(s)}_{m,t}\circ h^{(s)}_{m+1,t}$. Hence,
putting $y = h^{(s)}_{m+1,t}z$, we get
\begin{equation} \label{eq-lem-Phi-fi}
F^{(s)}_{m+1,t}(f^{(s)}_{m+1,t}z,t) =
H^{(s)}_{m+1}(h^{(s)}_{m+1,t}z,t) + F^{(s)}_m(f^{(s)}_{m,t}y,t)\;.
\end{equation}
Applying the induction hypothesis and inequality
\eqref{eq-lem-Phi-th} we get that the left hand side of
\eqref{eq-lem-Phi-fi} is greater than or equal to
$$C:= -\pi |z_{m+1}|^2 - \pi \sum_{j=2}^m |(y)_j|^2.$$
Applying \eqref{eq-lem-Phi-fo} we see that $(y)_j = z_j$ for
$j=2,...,m$. Hence $$C = -\pi (\rho_2+...+\rho_m),$$ and
inequality \eqref{eq-lem-Phi-tw} follows.
\end{proof}
This completes the proof of the Claim, and hence of \fullref{lem-fund}.
\end{proof}

\subsection{Some preparations}

Let $(M,\omega,L)$ be a Liouville manifold with the ideal contact
boundary $(P,\eta)$. We fix a contact form on $P$ and equip the
symplectization $SP$ with coordinates $(x,u) \in P \times \R_+$ as
in \fullref{subsec-geom-ze}. Thus we have decomposition
\eqref{eq-decomposition}
$$M = (P \times \R_+) \sqcup \Core(M)\;.$$
We extend the function $u$ by $0$ to $\Core(M)$. As before, we
introduce the functions $\rho_1 = \pi|z_1|^2$ and $\varrho =
\sum_{j=2}^n \pi|z_j|^2$ on $\C^n$, and denote $\rho := \rho_1
+\varrho$.

In what follows we have to work with starshaped domains and
hypersurfaces in the stabilization $M_{stab}= M \times \C^n$. A
user-friendly class of these objects arises, roughly speaking, as
follows: Take any closed starshaped domain $D \subset \C^n$ with
boundary $\Gamma = \p D$. Then $\{u \leq 1\} \times D$ is a closed
starshaped domain in $M_{stab}$ whose boundary is a starshaped
hypersurface. The problem is that these objects (the domain and
the hypersurface) have corners unless $M$ is a point, and hence do
not fit into our setup. However, this can be easily corrected by
an appropriate smoothing of corners. Below we describe this
construction in a slightly more general context when the domain
$D$ has the form $D = \{c_1\rho_1 + c\varrho \leq 1\}$ and
$\Gamma=\{c_1\rho_1 + c\varrho = 1\}$ with $c_1,c \in \R$. (Note
that $D$ and $\Gamma$ are not necessarily compact.)

We fix a real concave function $\theta\co [0,1] \to [0,1]$ (a
{\it mollifier}) which equals $1$ near $0$ and equals $0$ at $1$.
We assume, in addition, that its inverse $\theta^{-1}$, which is
defined near $0$, is flat at $0$, ie all its derivatives at $0$
vanish. The graph of $\theta$ is a smoothing of the union of
northern and eastern edges of the unit square on the coordinate
plane. In addition, we need another technical condition. Fix a
number $\lambda \in (0,1)$ which is sufficiently close to $1$ and
assume that
\begin{equation} \label{eq-mollifier}
\theta'(u) < -\lambda^2/4 \;\;\;\text{provided}\;\;\; \theta (u)
\leq \lambda^2\;.
\end{equation}
The desired smoothing of corners gives rise to a domain $$D_{ext}
= \{c_1\rho_1 + c\varrho \leq \theta(u)\}$$ and to a hypersurface
$$\Gamma_{ext}=\{c_1\rho_1 + c\varrho =\theta(u)\}$$ in $M_{stab}$.
The subscript $ext$ indicates that we extended an object defined
on $\C^n$ to $M \times \C^n$. Note that with this notation the
hypersurface $S^{2n-1}_{ext}$ represents the ideal contact
boundary of $M_{stab}$.

\subsection{Distinguished $\R_+$--equivariant symplectomorphism}
\label{subsec-dist-cont}

Consider the following subsets of $\C^{n}$: the domain
$$W := \{(2n+1)\varrho -(n-1)\rho_1 \leq 1 \}\;,$$
the ellipsoid
$$E:= \{\nu \rho_1 +\varrho
=1\;\}\;,$$ and the ball $$B= \{\rho \leq 1\}.$$
For $c>0$ let us denote by $Y_c\co  \C^n \to \C^n$ the shift by
$c$ along the $\mathrm{Re}(z_1)$--axis, and write $\wt{Y}_c$ for its
trivial extension $\id \times Y_c$ to a diffeomorphism of $M
\times \C^n$.

 \begin{lemma}\label{lem-shift} There exist $c>0$ large
enough and $\nu>0$ small enough such that all ellipsoids
$Y_s(E),\; s \in [0,c]$, enclose $0$ and hence are starshaped, and
in addition
\begin{equation} \label{eq-ellip-vsp-on}
Y_b(E) \subset \mathrm{Interior}(W)\;\;\;\hbox{for all}\; b \geq
c\;.
\end{equation}
\end{lemma}
\begin{proof} Let us make analysis of the problem: suppose that
$c$ and $\nu$ are as required. Take any $z = (z_1,w) \in E$. Put
$z' = Y_b(z) = (z_1+b,w)$ with $b \geq c$. We have to check that
$$(n-1)\pi|z_1+b|^2 > (2n+1)\pi|w|^2 -1\;.$$
Substituting $\pi |w|^2 = 1- \pi \nu |z_1|^2$ we rewrite this
inequality as
$$(n-1)\pi|z_1+b|^2 + (2n+1)\nu\pi|z_1|^2 > 2n\;.$$
Write $z_1 = p+iq$ and observe that it suffices to achieve
\begin{equation}
\label{eq-lem-shift-vsp} (n-1)\pi(p+b)^2 + (2n+1)\nu\pi b^2 > 2n
\end{equation}
for all $b \geq c$. The absolute minimum of the function in the
left hand side of this inequality equals
$$\frac{(n-1)(2n+1) \pi\nu b^2}{(n-1)+(2n+1)\nu}\;.$$
Substitute this into the left hand side of
\eqref{eq-lem-shift-vsp}. Replacing $b$ by $c$ and performing an
elementary algebraic manipulation we get
\begin{equation}\label{eq-lem-shift-vsp-tw}
\pi\nu c^2 > \frac{2n}{2n+1} + \nu \frac{2n}{n-1}\;.
\end{equation}
On the other hand, $ Y_s(E)$ encloses $0$ for all $s \leq c$
provided that
\begin{equation}
\label{eq-nu-c-on} \pi \nu c^2 < 1\;.
\end{equation}
To complete the proof it remans to notice that inequalities
\eqref{eq-lem-shift-vsp-tw} and \eqref{eq-nu-c-on} are compatible if
$\nu$ is small enough. \end{proof}

 Take $c$ and $\nu$ from the lemma and choose $\lambda
\in (0,1)$ so that
\begin{equation}
\label{eq-nu-c} \pi \nu c^2 \leq \lambda^2 < 1\;
\end{equation}
(see inequality \eqref{eq-nu-c-on} above). In what follows we
assume that this number $\lambda$ appears in the condition
\eqref{eq-mollifier} on the mollifier $\theta$.

\begin{lemma}\label{lmCLAIM} All domains $\wt{Y}_s(E_{ext}),\; s
\in [0,c], $ are starshaped in $M_{stab}$, and in addition
\begin{equation} \label{eq-ellip-vsp}
\wt{Y}_c(E_{ext}) \subset \mathrm{Interior}(W_{ext})\;.
\end{equation}
\end{lemma}
\begin{proof}
Write $z_j = p_j +iq_j$, and note that the Liouville field on
$M_{stab}$ can be written as
$$L=\partial_u \oplus \frac{1}{2}\sum_{j=1}^n (p_j \frac{\p}{
\partial{p_j}}
+ q_j \frac{\p}{ \partial{q_j}})\;.$$ Choose $s \in [0,c]$. The
equation of $\wt{Y}_s (E_{ext})$ is
$$ F(p,q,u) := \pi \nu (p_1-s)^2 + \pi \nu q_1^2 + \pi\cdot \sum_{j=2}^n (p_j^2 +q_j^2) - \theta(u)\ =0\;.$$
We claim that $dF(L)>0$ at every point of the hypersurface
$\wt{Y}_s (E_{ext})$. Indeed, calculating $dF(L)$ at such a point
we get
$$dF(L) = \pi \nu s (p_1-s) +\theta(u) - \theta'(u).$$
On $\wt{Y}_s (E_{ext})$ we have $|p_1-s| \leq
\sqrt{\theta(u)/(\pi\nu)}$. Thus inequality \eqref{eq-nu-c} yields
$$dF(L) \geq \theta(u) - \lambda \sqrt{\theta(u)} -\theta'(u)\;.$$
The positivity of the right hand side, say $A$, of this
expression follows from assumption \eqref{eq-mollifier} on the
mollifier $\theta$. Indeed, denote by $I$ the segment
$$\{u\;:\; \theta (u)\leq \lambda^2\}\;.$$ For $u \notin I$ we have $\theta(u) - \lambda
\sqrt{\theta(u)}
>0$ and $\theta'(u) \leq 0$ so $A>0$. On $I$ we have
$\theta(u) - \lambda \sqrt{\theta(u)} \geq -\lambda^2/4$ and
$-\theta'(u) > \lambda^2/4$. Thus again $A>0$.

The ``singular" case when $u=1$ (and so $\theta'(u)$ is not
defined) can be easily checked separately. The claim follows.

Thus $\wt{Y}_s(E_{ext}), s \in [0,c]$, is an isotopy of an
(obviously) starshaped hypersurface $E_{ext}$ in the class of
hypersurfaces transversal to the Liouville field $L$. Therefore
all these hypersurfaces are starshaped.

 In order to verify inclusion \eqref{eq-ellip-vsp}, we use
a scaling type argument: fix $u$ and put $\varrho' =
\varrho/\theta(u), \rho_1' = \rho_1/\theta(u), c' =
c/\sqrt{\theta(u)}$. Then the desired inclusion is equivalent to
inclusion \eqref{eq-ellip-vsp-on} with $b = c' \geq c$ since
$\theta(u) \leq 1$. This completes the proof.
\end{proof}

Applying Ambient Isotopy \fullref{lem-ambisot} we get a
family $\{a^{(s)}\},\; s \in [0,c]$, of equivariant
symplectomorphisms of $M_{stab}$ such that $a^{(s)} (E_{ext}) =
\wt{Y}_s(E_{ext})$ for all $s \in [0,c]$. Symplectomorphism
$a:=a^{(c)}$ plays an important role below.\eject
\begin{lemma}\label{lma-properties}
The equivariant symplectomorphism $a$ satisfies the following
conditions:
\begin{gather} \label{eq-sq-pr-vspon}
a(\partial B_{ext}) \subset \mathrm{Interior}(W_{ext})\;;\\
\label{eq-firstintegral} a^*u = u,\;\;\; a^*\varrho = \varrho\;.
\end{gather}
\end{lemma}
\begin{proof}
The property \eqref{eq-sq-pr-vspon} is straightforward. Let us
verify formulas \eqref{eq-firstintegral}. Note that the
characteristic foliation $\sigma_0$ of $E_{stab}$ is generated by
$$\nu \cdot \mathrm{sgrad}_{\C^n} \rho_1 + \mathrm{sgrad}_{\C^n} \varrho- \theta'(u)\mathrm{sgrad}_{M}
u\;.$$ The functions $u$ and $\varrho$ are constant along
$\sigma_0$, and, moreover, the shifts $\wt{Y}_s$ preserve them.
Thus, using Ambient Isotopy \fullref{lem-ambisot}, we get
$u(a^{(s)}x) = u(x)$ and $\varrho(a^{(s)}x) = \varrho(x)$ for all
$x \in E_{stab}$. In view of the homogeneity of $a$ we conclude
\eqref{eq-firstintegral}.
\end{proof}

\subsection{The main construction} \label{subsec-main-consttw}

 {\bf Convention}\qua{\em All transformations $g$ from $SU(n)$
smoothly extend to $M \times \C^n$ by $\id \times g$. We denote
the extensions by the same letter as the original
transformations.}

\medskip
 Put $ e_tz:= e^{2\pi i t }z$. Let $f_t$ be the loop of
unitary transformations constructed in Fundamental \fullref{lem-fund}. Let $a$ be the distinguished $\R_+$--equivariant
symplectomorphism of $M_{stab}$ constructed in \fullref{subsec-dist-cont}.

\begin{thm} \label{thm-main-consttw}
The loop
$$\varphi_t: = e_{-t}f_{3t}ae_{t}a^{-1}$$ is a positive contractible
loop of $\R_+$--equivariant symplectomorphisms of $$SP_{stab} =
M_{stab}\setminus \Core(M_{stab}).$$
\end{thm}

  In view of \fullref{prop-nonderability} this
yields \fullref{thmtw-stabilization}. In the case of the
standard contact 3--sphere somewhat similar positive contractible
loops can be defined by an explicit analytic formula. Furthermore,
"a cousin" of the distinguished map $a$ admits a transparent
geometric meaning. We refer to \fullref{rem-aptw-on} in Appendix
B for the details.

\begin{proof} Let $f_t$ and $f_t^{(s)}$ be the loop of
unitary transformations and its homotopy from Fundamental \fullref{lem-fund}, and let $F(z)$ and $F^{(s)}(z,t)$ be the
corresponding Hamiltonians considered as functions on $M \times
\C^n$.

\medskip
{\bf (1)}\qua Put $g_t:= f_t ae_{t}a^{-1}$ and $g^{(s)}_t:=
f^{(s)}_t ae_{t}a^{-1}$, and denote by $G(x,z,t)$ and
$G^{(s)}(x,z,t)$ their Hamiltonians. Here $x \in M, z \in \C^n$.
We claim that
\begin{equation} \label{eq-pr-main-on}
G^{(s)} \geq 0
\end{equation}
and
\begin{equation} \label{eq-pr-main-tw}
G(x,z,t) > 2n\varrho(z)
\end{equation}
for all $(x,z) \in M_{stab}, t \in S^1, s \in [0,1]$. Indeed,
$$G^{(s)}(x,z,t) = F^{(s)}(z,t) + \rho(a^{-1}(x,(f^{(s)}_t)^{-1}z)).$$
Recall that $\rho = \rho_1 + \varrho$. Put $y
=(f^{(s)}_t)^{-1}z$. In view of $\R_+$--equivariance, it suffices
to verify inequalities \eqref{eq-pr-main-on} and
\eqref{eq-pr-main-tw} for fixed $s$ and $t$ assuming that $(x,z)$
runs over any given starshaped hypersurface. Thus, we can assume
that
\begin{equation}\label{eq-Bext}
(x,y) \in a(\partial B_{ext})\;.
\end{equation}
By Fundamental \fullref{lem-fund} and equations
\eqref{eq-firstintegral} we have
$$G^{(s)}(x,z,t) = F^{(s)}(f^{(s)}_ty,t) + \rho(a^{-1}(x,y)) \geq -\varrho(y)+
\rho(a^{-1}(x,y))$$ $$\geq -\varrho(y)+ \varrho(a^{-1}(x,y))
=-\varrho(y)+ \varrho(y)= 0,$$ which proves inequality
\eqref{eq-pr-main-on}.

Let us now turn to the case $s=1$. Inclusion
\eqref{eq-sq-pr-vspon} guarantees that
\begin{equation}\label{eq-pr-main-th}
\theta(u(x,y))+ (n-1)\rho_1(y) -(2n+1)\varrho(y) > 0\;.
\end{equation}
Note that when $s=1$ we have $y = f_{-t}z$, and hence
\begin{equation} \label{eq-pr-main-vspon}
\rho_1(y)=\rho_1(z) \;\;\text{and}\;\;\varrho(y)=\varrho(z).
\end{equation}
Using equations \eqref{eq-firstintegral}, inequality
\eqref{eq-pr-main-th} and inclusion \eqref{eq-Bext} we get
$$G(x,z,t) = F(y)+ \rho(a^{-1}(x,y)) = (n-1)\rho_1(y) -\varrho(y)+\theta(u(a^{-1}(x,y)))$$
$$=(n-1)\rho_1(y) -\varrho(y)+ \theta(u(x,y)) > 2n\varrho(y) = 2n\varrho(z).$$
This proves inequality \eqref{eq-pr-main-tw}, and the claim
follows.

\medskip
{\bf (2)}\qua Put $h_t := f_{2t}g_t$, $h_t^{(s)} =
f_{2t}^{(s)}g_t$ and let $H(x,z,t)$ and $H^{(s)}(x,z,t)$ be
corresponding Hamiltonians. We claim that
\begin{equation} \label{eq-pr-main-onon}
H^{(s)} \geq 0\;
\end{equation}
and
\begin{equation} \label{eq-pr-main-ontw}
H(x,z,t) > (2n-2)\rho(z)\;
\end{equation}
for all $(x,z) \in M_{stab}, t \in S^1, s\in [0,1]$. Indeed,
$$H^{(s)}(x,z,t) = 2F^{(s)}(z,2t) +
G(x,(f^{(s)}_{2t})^{-1}z,t).$$ Put $y = (f^{(s)}_{2t})^{-1}z$.
Applying Fundamental Lemma and \eqref{eq-pr-main-tw} we have
$$H^{(s)}(x,z,t)=2F^{(s)}(f^{(s)}_{2t}y,2t) + G(x,y,t) \geq
-2\varrho(y) + 2n\varrho(y) \geq 0\;,$$ which proves inequality
\eqref{eq-pr-main-onon}.

When $s=1$ we apply formulas \eqref{eq-pr-main-vspon} and
\eqref{eq-pr-main-tw} and get
$$H(x,z,t) = 2(n-1)\rho_1(z) -2\varrho(z) + G(x,y,t) >
(2n-2)(\rho_1(z)+\varrho(z))\;.$$ This proves inequality
\eqref{eq-pr-main-ontw}. The claim follows.

\medskip
{\bf(3)}\qua Note that the loop $\{\varphi_t\}$ presented in
the formulation of the theorem is given by $\varphi_t= e_{-t}h_t$.
Let $\Phi(x,z,t)$ be the Hamiltonian of $\{\varphi_t\}$. Then
using inequality \eqref{eq-pr-main-ontw} we get
$$\Phi(x,e_{-t}z,t) = -\rho(z) + H(x,z,t) >
(2n-3)\rho(z) \geq 0$$ and hence
\begin{equation}\label{eq-positiveloop-bound}
\Phi(x,z,t) > (2n-3)\rho(z) \geq 0
\end{equation}
for all $(x,z) \in M_{stab}$ and $t \in S^1$. This proves
positivity of the loop $\{\varphi_t\}$.

\medskip
{\bf (4)}\qua Let us present a homotopy $\Delta$ of the loop
$\{\varphi_t\}$ to the constant loop. For that purpose write
$$\varphi_t = e_{-t} f_{2t} f_t ae_{t}a^{-1}\;.$$
The homotopy $\Delta$ is given in three steps.

{\bf Step 1}\qua Contract the term $f_{2t}$ using homotopy
$f_{2t}^{(s)}$ given by the Fundamental Lemma and arrive to
$$\varphi'_t = e_{-t} f_t ae_{t}a^{-1}\;.$$

{\bf Step 2}\qua Contract the term $f_{t}$ using homotopy
$f_{t}^{(s)}$ given by the Fundamental Lemma and arrive to
$$\varphi''_t = e_{-t}ae_ta^{-1}\;.$$

{\bf Step 3}\qua Use the path $a^{(s)}$ joining the
identity with $a$ described in \fullref{subsec-dist-cont}. We
get the homotopy
$$e_{-t}a^{(s)}e_t(a^{(s)})^{-1},$$
which equals the constant loop when $s= 0$.

 This completes the proof. \end{proof}

 Let us emphasize that the proof above works when $M =
\{\mathrm{point}\}$ with the following modifications:
\begin{itemize}
\item We do not need any mollifier $\theta$, and work directly with $W,E,B
\subset \C^n$, omitting sub-index $ext$;
\item The shift $\wt{Y}_s$ coincides with $Y_s$.
\end{itemize}
Taking into account these remarks we are ready to prove \fullref{mu}(ii).

\begin{proof}[Proof of \fullref{mu}(ii)] Look at the
homotopy $\Delta$ of the loop $\varphi_t$ described in Part {\bf
4} of the proof of \fullref{thm-main-consttw} above. It
suffices to show that the Hamiltonians generating loops obtained
in the process of this homotopy are $\geq-\pi|z|^2$. In Steps 1
and 2 of the homotopy we get loops $e_{-t} h^{(s)}_t$ and $e_{-t}
g^{(s)}_t$, respectively. Inequalities \eqref{eq-pr-main-onon} and
\eqref{eq-pr-main-on} guarantee that in both cases the
corresponding Hamiltonians are $\geq -\pi |z|^2 .$ In Step 3 the
Hamiltonian again is $\geq-\pi|z|^2$. This shows that $\mu(\Delta)
\leq 1$ (and in fact, one can see that $\mu(\Delta) = 1$, as it
should be in accordance with \fullref{mu}(i).) This completes
the proof. \end{proof}

\subsection[Proof of \ref{thm-sq-strong}]{Proof of \fullref{thm-sq-strong}} \label{sec-sq-strong}

We work in $\C^n$ with coordinates $(z_1,...,z_n)$. Put $\rho_j =
\pi|z_j|^2$, $\varrho = \sum_{j=2}^n \pi|z_j|^2$ and $\rho =
\sum_{j=1}^n \pi |z_j|^2$. We will prove that the cylinder $C =
\{\rho_2 < 1\} \subset \C^n$ is negligible. 
\begin{gather*}
D =\big\{\varrho < \frac{1}{n}\big\}\tag*{\hbox{Put}} \\ 
W = \{(n+1)\rho_2 + n\rho_3 +...+n\rho_n -\rho_1 < 1\}\;.
\tag*{\hbox{and}}
\end{gather*}
Introduce the loops of transformations from $SU(n)$
\begin{gather*}f_t(z) = (e^{2\pi (n-1)i t}z_1, e^{-2\pi i t}z_2,...,e^{-2\pi i
t}z_n)\\
g_t(z) = (e^{2\pi i t}z_1, e^{-2\pi i t}z_2, z_3,...,z_n)\;
\tag*{\hbox{and}}
\end{gather*}
which are generated by the Hamiltonians 
\begin{gather*}F(z) = (n-1)\rho_1(z) -
\varrho(z)\\G(z) = \rho_1(z) -\rho_2(z),
\tag*{\hbox{and}}
\end{gather*} respectively.
Our main squeezing tools will be contact embeddings
\begin{gather*}\Phi(z,t) = \frac{f_t(z)}{\sqrt{1+F(z)}}\\
\Psi(z,t) = \frac{g_t(z)}{\sqrt{1+G(z)}}\tag*{\hbox{and}}
\end{gather*}
associated to loops $\{f_t\}$ and $\{g_t\}$, respectively, in
accordance to \fullref{prop-loop-cont}. (The square root
in these formulas reflects the fact that the natural $\R_+$--action
on $\C^n \setminus \{0\}$ considered as the symplectization of
$S^{2n-1}$ is given by $z \to \sqrt{c}z$ for $c \in \R_+$.)

In view of Fundamental \fullref{lem-fund} the loop $\{f_t\}$ is
homotopic to a point through loops whose Hamiltonians are $\geq
-\varrho$. In particular, the embedding $\Phi\co  D \to V$ is well
defined and contactly isotopic to the inclusion.

In view of Claim \eqref{eq-lem-Phi-tw} inside the proof of
Fundamental Lemma, the loop $\{g_t\}$ is homotopic to a point
through loops with Hamiltonians $\geq -\rho_2$. In particular, the
embedding $\Psi\co  C \to V$ is well defined and contactly isotopic
to the inclusion.

Let $K \subset C$ be any compact subset. Arguing exactly as in
\fullref{subsec-dist-cont} we find a contact isotopy of $K$ to
its shift $Y(K)$ along the $\mathrm{Re} z_1$--axis such that $Y(K)
\subset C \cap W$. Take $z \in C \cap W$ and put $z' = \Psi(z)$.
Then
$$\varrho(z') = \frac{\varrho(z)}{1+\rho_1(z)-\rho_2(z)} <
\frac{1}{n}\;,$$ where the equality follows from the explicit
formula for $\Psi$ and the inequality follows from the definition
of $W$. We conclude that $\Psi(Y(K)) \subset D$. Take now any $u
\in D$ and put $u' = \Phi(u)$. Then
$$\rho(u') = \frac{\rho_1(u) + \varrho(u)}{1+(n-1)\rho_1(u) -
\varrho (u)} \leq \frac{\rho_1(u) + \frac{1}{n}}{1 +
(n-1)\rho_1(u) -\frac{1}{n}} = \frac{1}{n-1}\;,$$ where the
equality on the left follows from the explicit formula for $\Phi$
and the inequality follows from the definition of $D$. We conclude
that
$$\Phi(\Psi(Y(K))) \subset \wh{B}^{2n}(\frac{1}{n-1})\;.$$
Thus, we contactly squeezed $K$ into the set
$\wh{B}^{2n}(\frac{1}{n-1})$, which is negligible in view of
\fullref{thm-sq}. This shows that any compact subset of $C$ is
negligible, as required.\qed

\section{Different flavors of Floer homology}\label{secflavors}

To prove our results on contact non-squeezing and orderability we
need to understand the relation between cylindrical contact
homology which is a special case of the symplectic field theory
(see \cite{SFT,Ustilovsky,Yau,B}), and
(periodic) Floer homology for the symplectic action functional
(see \cite{Floerr,MS}). In order to do that we introduce a more
general version of Floer homology for {\it Hamiltonian structures}
which contains both theories as its special cases.

\subsection{Stable Hamiltonian structures}\label{secstable}

Two differential 1--forms on a manifold are called {\it equivalent}
if they differ by an exact 1--form. We denote by $(\Theta)$ the
equivalence class of a 1--form $\Theta$. The 2--form $d\Theta$ does
not depend on the choice of a representative $\Theta \in (\Theta)$
and will be denoted by $d(\Theta)$. A {\it Hamiltonian structure}
$\cH$ on an odd-dimensional oriented, possibly non-compact
manifold $V$ is an equivalence class $(\Theta)$ of 1--forms such
that $\Omega=d(\Theta)$ has the maximal rank. The tangent line
field $\ell=\Ker\Omega$ is called the {\it characteristic} line
field. The field $\ell$ integrates to a 1--dimensional {\it
characteristic foliation} of $\Omega$. Note that $\Omega$ defines
a fiberwise symplectic structure (and hence an orientation) on the
bundle $TV/\ell$. Thus the line bundle $\ell$ is equipped with an
orientation. We will call {\it characteristic} any vector field
$R$ which generates $\ell$ and respects its orientation.

 Any co-orientable hypersurface $V$ in an exact
symplectic manifold $(W,\wt\Omega=d\wt\Theta)$ inherits a
Hamiltonian structure $(\Theta)=(\wt\Theta|_{V})$. Conversely, any
Hamiltonian structure $(V,(\Theta))$ embeds as a hypersurface in a
symplectic manifold $(V\times(-\e,\e),d\wt\Theta)$ where the form
$\wt\Theta$ can be constructed as follows. Let $\lambda$ be any
$1$--form which is not vanishing on $\ell$, and $s$ be the
coordinate along the second factor. Then we set
$\wt\Theta=\Theta+s\lambda$. Note that by Darboux' theorem the
Hamiltonian structure $(V,(\Theta))$ determines the symplectic
structure in its neighborhood uniquely up to a diffeomorphism
fixed on $V$. We call (a germ along $V$ of ) the symplectic
structure $\wt\Omega$ on $V\times{(-\e,\e)}$ the {\it symplectic
extension} of $(V,(\Theta))$.

 A Hamiltonian structure $(V,(\Theta))$ is called {\it stable} (see \cite{HZ})
if its symplectic extension can be realized by a form $\wt\Omega$
on $V\times{(-\e,\e)}$, such that all the Hamiltonian structures
$\cH_s,\,s\in(-\e,\e)$ induced on $V$ by the inclusions $V=V\times
s\hookrightarrow V\times(-\e,\e)$ have the same characteristic
line field $\ell$.
\begin{proposition}\label{propstability}
A Hamiltonian structure $(V,(\Theta))$ is stable if and only if
there exists a $1$--form $\lambda$ and a characteristic vector
field $R$ such that
\begin{equation}\label{eqframing}
\lambda(R)=1 \;\;\;\text{and}\;\;\;i_R d\lambda=0 \,.
\end{equation}
\end{proposition}

 Note that in view of the Cartan formula we have
$L_R\lambda=d(\lambda(R))+i_R d\lambda,$ and hence the second
condition can be restated as invariance of $\lambda$ under the
flow of $R$.

\begin{proof}
Suppose that there exists $\lambda$ which satisfies the above
conditions. Set $\wt\Theta=\Theta+s\lambda$ and
$\wt\Omega=d\wt\Theta$. Then for $\wt\Omega_s:=\wt\Omega|_{V\times
s}=\Omega+sd\lambda$ we have
$$i_R\wt\Omega_s=i_R\Omega+si_R d\lambda=0\,,$$ and hence
$(V, (\Theta))$ is stable.

 Conversely, suppose that a Hamiltonian structure $(V,
(\Theta))$ is stable. Then according to the definition it has a
symplectic extension $(V\times(-\e,\e),\wt\Omega=d\wt\Theta)$ such
that $i_R\wt\Omega_s=0$, where $\wt\Omega_s=\wt\Omega|_{V\times
s}$ and $R$ is a characteristic vector field for $\Omega$. Let us
write $\wt\Omega$ as $\Omega+\eta_s+\lambda_s\wedge ds$, where
$\lambda_s$, $\eta_s$ are families of $1$-- and $2$--forms on $V$
and $\eta_0=0$. We have
$$0=d\wt\Omega=d\Omega+d\eta_s-\dot{\eta}_s\wedge ds
+d\lambda_s\wedge ds= d\eta_s-\dot{\eta}_s\wedge
ds+d\lambda_s\wedge ds\,,$$ where
$\dot{\eta}_s=\frac{d\eta_s}{ds}$. Hence,
$d\lambda_s=\dot{\eta}_s$. On the other hand,
$$0=i_R\wt\Omega_s= i_R\eta_s=0\,.$$
Differentiating with respect to $s$ we get $i_R\dot{\eta}_s=0$,
and therefore, $i_R d\lambda_s=0$. Note that $\lambda_0(R)$ does
nor vanish. Indeed along $V\times 0$ we have
$\wt\Omega^{n+1}=\Omega^{n}\wedge \lambda_0\wedge ds,$ and thus
$0\neq i_R{\wt\Omega}^{n+1}=\lambda_0(R)\Omega^{n}\wedge ds$.
(Here $\dim V = 2n+1$.) Normalizing $R$ in such a way that
$\lambda_0(R)\equiv 1$, we get that $ \lambda_0(R)=1$ and $i_R
d\lambda_0=0$, as required.
\end{proof}

 The structures characterized by existence of a $1$--form
$\lambda$ satisfying the conditions of \fullref{propstability} appeared in \cite{SFT-compact}, but were not
identified there with stable hypersurfaces in symplectic manifolds
which were first studied in \cite{HZ}.

 The $1$--form $\lambda$ as in \fullref{propstability} is called a {\it framing}
of the stable Hamiltonian structure $\cH= (V,(\Theta))$. The
hyperplane field $\xi=\{\lambda=0\}$, called a {\it cut} of $\cH$,
and the characteristic vector field $R$, called a {\it Reeb field}
of the Hamiltonian structure, are uniquely determined from
\eqref{eqframing} by the framing $\lambda$. The triple
$(V,(\Theta),\lambda)$ is called {\it a framed Hamiltonian
structure} and is denoted by $\ocH$. Here are examples of stable
Hamiltonian structures taken from \cite{SFT-compact}.

\begin{exam}\label{exstable-Hamiltonian}\ 

{\bf 1\qua Contact forms}\qua Let $\xi$ be a contact
structure on $V$ and $\lambda$ a corresponding contact form. Then
$(V,(\lambda),\lambda)$ is a framed stable Hamiltonian structure
with cut $\xi$. In this case $R$ the usual Reeb field of the
contact form $\lambda$ and $\Omega=d\lambda$.

{\bf 2\qua Hamiltonian functions}\qua Let $(M,\omega=d\alpha)$
be an exact symplectic manifold and $H_t\co M\to\R$, $t\in
S^1=\R/\Z$, a $1$--periodic time-dependent Hamiltonian function.
Write $\sgrad{H_t}$ for its Hamiltonian vector field. Put
$V=M\times S^1$. Set $\Theta=-\alpha+H_tdt$ and $\lambda=dt$. Then
$(V,(\Theta),\lambda)$ is a framed stable Hamiltonian structure.
Its Reeb vector field is given by $R=\frac{\p}{\p t}+\sgrad{H_t}$
and its cut $\xi=\{\lambda=0\}$ is formed by tangent spaces to
$M\times t, t\in S^1$.

{\bf 3\qua Prequantization spaces}\qua Let $(M,\omega)$ be a
symplectic manifold with the integral cohomology class $[\omega]$
of the symplectic form $\omega$. Consider the corresponding
prequantization space, that is a principal $S^1$--bundle $p\co V\to M$
with the first Chern class $[\omega]$. In this case the lift
$\Omega:= p^*\omega$ of the symplectic form is exact: $\Omega =
d\Theta$. Then $(V,(\Theta))$ is a stable Hamiltonian structure.
Indeed, one can choose any $S^1$--connection form $\lambda$ as its
framing. The corresponding Reeb vector field $R$ is the
infinitesimal generator of the $S^1$--action, and the cut $\xi$ is
the horizontal distribution of the connection.
\end{exam}

\subsection{Special Hamiltonian structures}\label{subsec-ham-struc}

As it is stated above, the goal of \fullref{secflavors} is to
develop the Floer homology theory for Hamiltonian structures. To
avoid unnecessary technicalities we will restrict the theory to a
special class of stable Hamiltonian structures which will be
sufficient for all the applications considered in this paper. In
\fullref{subsec-disc-nonspecial} we briefly discuss possible
generalizations of the theory.

 Let $(M^{2n},\omega,L)$ be a Liouville manifold with the Liouville
form $\alpha$. We will work on the manifold $V = M \times S^1$,
write $t \;(\mathrm{mod} 1)$ for the coordinate on $S^1=\R/\Z$, and
orient $V$ by the volume form $(-\omega)^n \wedge dt$.

 Let us denote by $\digamma_{\sC,\sK}$ the set of all
1--forms on $V$ which coincide with $\sC dt-\sK\alpha$ ($\sC,\sK
\in \R$) outside a compact subset of $V$. We identify the set
$\digamma:= \bigcup_{\sC, \sK \in \R} \digamma_{\sC,\sK}$ with $\R
\times \R \times \digamma_{0,0}$ and equip it with the product
topology, where $\digamma_{0,0}$ is considered with the strong
Whitney topology. It induces a topology on the set of equivalence
classes of 1--forms.

 A Hamiltonian structure $\cH=(V,(\Theta))$ on $V$ is
called {\it special} if $\Theta \in \digamma_{\sC,\sK}$ , where
$\sC$ and $\sK$ are positive constants which are called {\it the
structure constants} of $\cH$, and the following Axioms 1--3 are
satisfied.

\medskip
 {\bf Axiom 1}\qua {\sl The first Chern class of the
symplectic vector bundle $(TV/\ell,\Omega=d(\Theta))$, where
$\ell$ is the characteristic line field of $\cH$, vanishes on any
2--cycle represented by a mapping $\mathbb{T}^2 \to V$. }

\medskip
 Recall that the characteristic foliation $\cT$ of $\cH$
is equipped with an orientation. The {\it action} of a periodic
orbit $\g$ of $\cT$ is defined by
\begin{equation}\label{eqaction}
\cA(\gamma)=\int \limits_{\g}\Theta\;.
\end{equation}
By Stokes' formula, the action does not depend on the specific
choice of a 1--form $\Theta$ representing $(\Theta)$. In view of
the conditions ``at infinity" imposed on $\Theta$, the
trajectories of $\cT$ outside a sufficiently large compact subset
of $V$ are circles $\{\mathrm{point}\} \times S^1$, and their action
equals to the structure constant $\sC$.

\medskip
 {\bf Axiom 2}\qua {\sl The characteristic foliation $\cT$
has no contractible closed orbits of action $\leq \sC$.}

\medskip
{\bf Axiom 3}\qua {\sl There exist numbers $\mu \geq 0, \sP
> 0,\sQ \geq 0 $ and a 1--form $\lambda \in \digamma_{\sP,\sQ}$ such that
$\lambda - \mu \Theta$ is a closed 1--form representing cohomology
class $(\sP-\mu \sC)[dt]$ and $\lambda$ is positive on the
characteristic foliation $\cT$ with respect to its natural
orientation.}

\medskip  Note that in view of Axiom 3 we have $d\lambda=\mu\Omega$, which
implies stability of $\cH$. In particular, $\lambda$ is a framing
of the Hamiltonian structure $\cH$. Such a framing is called {\it
special}. Special framings form a convex cone. The constant $\mu$
is called the {\it parameter}, and the numbers $\sP$ and $\sQ$
{\it the structure constants} of the framing. Note that $\sQ=\mu
\sK$. Clearly, the cut $\xi$ is a contact structure on $V$ if $\mu
>0$, and $\xi$ integrates to a foliation defined by a closed $1$--form if $\mu=0$.
Sometimes, in order to emphasize the dependence of the introduced
objects on $\cH$ we will write
$(\Theta)_{\cH},\Omega_{\cH},\cT_{\cH},\cA_{\cH}$ etc.

 Consider a special Hamiltonian structure
$(V,(\Theta))$ with a special framing $\lambda$. We start with the
following useful {\it period-action equation}, which is an
immediate consequence of the definitions. Let $\gamma$ be a
$T$--periodic closed orbit of the Reeb vector field $R$. Then
\begin{equation}\label{eq-per-act-on}
T = \mu \cA(\gamma) + (\sP-\mu \sC)\langle[dt],[\gamma]\rangle \;.
\end{equation}

\begin{prop}\label{prop-special-stab}
Let $(\Theta)$ be a special Hamiltonian structure on $V$. Then
every 1--form $\Theta'$ in a sufficiently small neighborhood of
$\Theta$ defines a special Hamiltonian structure. Moreover, if
$\lambda$ is a special framing of $\Theta$ then there exists a
special framing $\lambda'$ of $\Theta'$ which is sufficiently
close to $\lambda$.
\end{prop}

\begin{proof} If $\Theta'$
is sufficiently close to $\Theta$ then Axiom 1 for $\Theta'$ is
fulfilled automatically. To check Axiom 3, note that the 1--form
$$\lambda' = \lambda + \mu(\Theta'-\Theta)$$ is a
framing of $(\Theta')$ with the same parameter $\mu$. It remains
to verify Axiom 2. Note that if $\mu=0$ then in view of equation
\eqref{eq-per-act-on} the characteristic foliation $\cT$ of
$(\Theta')$ cannot have contractible closed orbits. Let us assume
that $\mu > 0$ and suppose on the contrary that there exists a
sequence of 1--forms $\Theta_i \to \Theta$ whose characteristic
foliations admit contractible closed orbits of action $\leq
\sC_i$, where $\sC_i$ is the structure constant of $\Theta_i$.
Applying again \eqref{eq-per-act-on} we see that the periods of
$\gamma_i$ are bounded, and thus these orbits converge to a
contractible closed orbit of $(\Theta)$ of action $\leq \sC$. This
contradicts to Axiom 2 for $\Theta$.
\end{proof}

 We will be mostly dealing in this paper with special
Hamiltonian structures, and hence, if not otherwise noted,

\medskip
\cl{\fbox{\parbox{0.8 \linewidth}{ \centerline{{\it all considered
Hamiltonian structures and }} \centerline{{\it their framings will
be assumed special.}} }}}

\subsection{Basic examples}\label{subsec-exa}
The following two basic examples of special Hamiltonian structures
which will appear below (see \cite{SFT-compact}) are
specifications of Examples \ref{exstable-Hamiltonian}.1 and
\ref{exstable-Hamiltonian}.2.

\begin{exam}[Contact forms]\label{ex-cont}
{\rm Let $F\co V \to \R$ be a positive function
which is equal to a constant $\sC$ outside a compact set. Consider
a contact form $\lambda:= F(dt-\alpha)$ on $V$ whose Reeb vector
field has no contractible closed orbits of period $\leq \sC$. We
call such a form {\it admissible}. Then $(\lambda)$ is a
Hamiltonian structure on $V$. The same 1--form $\lambda$ serves as
its framing, and so the framing parameter $\mu$ equals 1. All the
structure constants are equal to $\sC$. The cut $\xi =
\mathrm{Ker}(\lambda)$ is a contact structure on $V$, and the Reeb
vector field is the usual Reeb field of $\lambda$. One can easily
check that the symplectic subbundles $(\xi, \Omega=d\lambda)$ and
$(TM,-\omega)$ of $TV$ are homotopic, and thus have the same first
Chern classes (see \fullref{seccontact} for more details).
Hence, Axiom 1 is equivalent to the requirement that $c_1(TM)$
vanishes on any cycle represented by a mapping $\mathbb{T}^2 \to
M$. Verification of other axioms is straightforward. We denote the
obtained framed Hamiltonian structure by $\ocH_{\lambda}$, as it
is fully determined by the contact form $\lambda$ on $V$. }
\end{exam}

\begin{exam}[Hamiltonian functions]\label{ex-sympl}{\rm
Let $F\co V=M\times S^1 \to \R$ be a time-dependent $1$--periodic
Hamiltonian function on $M$ which is equal to a positive constant
$\sC$ outside a compact set. Let us introduce a Hamiltonian
structure $(F dt-\alpha)$ on $V$. The framing is chosen as
$\lambda:= dt$, and hence the framing parameter $\mu$ vanishes.
Thus the cut $\xi$ is the tangent bundle to the fibers $M
\times\{\mathrm{point}\}$, and the Reeb vector field is
$R=\frac{\p}{\p t}+ \mathrm{sgrad}{F}$. As in the previous example
it is straightforward that all axioms are satisfied provided that
$c_1(TM)$ vanishes on any cycle represented by a mapping
$\mathbb{T}^2 \to M$. This framed Hamiltonian structure is fully
determined by the function $F$ and will be denoted by $\ocH_{F}$.
}
\end{exam}

\subsection{Periodic orbits of the characteristic foliation
}\label{subsec-per-char}

Let $\cT$ be the oriented characteristic foliation of a
Hamiltonian structure $\cH = (\Theta)$. Write $\sC$ for the
structure constant. Denote by $\cS$ the set of free homotopy
classes of loops $S^1 \to V$ which project to the class
$[\mathrm{point} \times S^1]$ under the natural projection $V=M
\times S^1 \to S^1$. We will focus on those orbits of $\cT$ with
action $< \sC$ which represent classes from $\cS$, and denote the
set of such orbits by $\Pc$.

\begin{defin} \label{def-Ham-reg}{\rm
A Hamiltonian structure is called {\it regular} if all periodic
orbits from $\Pc$ are non-degenerate, that is their linearized
Poincar\'e return maps do not contain $1$ in their spectra.}
\end{defin}

\begin{prop} \label{prop-orb-regular} Regular Hamiltonian
structures form a Baire set in the space of all Hamiltonian
structures.
\end{prop}

\begin{proof}[Outline of the proof] A little nuance here is that
the ``flat" orbits at infinity, whose action equals $\sC$, are
certainly degenerate. One goes round this difficulty as follows.
Choose a starshaped hypersurface $P \subset M$ such that $M
=\mathrm{Core}(M) \bigsqcup SP$, and introduce coordinates $(x,u)$
on $SP= P \times \R_+$. In these coordinates $\alpha = u\beta$
where $\beta$ is the contact form on $P$ defined by the
restriction of $\alpha$. Denote by $R_{\beta}$ the Reeb field of
$\beta$ and by $T$ the minimal period of a closed orbit of
$R_{\beta}$.

 Let $\cH = (\Theta)$ be any Hamiltonian structure with
the structure constants $\sC$ and $\sK$. We can assume without
loss of generality that $\sK=1$ and $\Theta = \sC dt-u\beta$ for
$u \geq 1/2$. We claim that there exists an arbitrarily small
perturbation $\Theta_1$ of $\Theta$ such that $\Theta_1 = \Theta$
for $u \geq 2$ and the region $A:= \{1\leq u < 2\}$ does not
contain closed orbits of the characteristic foliation of
$(\Theta_1)$ representing classes from $\cS$. Indeed, write
$\Theta_1 = G(u)\Theta$, where $G(u) \equiv 1$ for $u \geq 2$. Put
$H(u) = \sC G'/(Gu)'$. Assume that for $u \in [1,2)$ the function $G$
is strictly increasing and $C^1$--close to the constant function
$u\equiv 1$. Then the characteristic foliation $\cT_1$ of
$(\Theta_1)$ on $A$ is generated by the field $\ppt +
H(u)R_{\beta}$ with $0< H(u) < 1/T$. Thus $\cT_1$ has no closed
orbits in $A$ representing classes from $\cS$. Now one can extend
$\Theta_1$ to the rest of $V$ to get the desired perturbation.

 It remains to perturb $\Theta_1$, keeping it fixed on the set
$\{u>3/2\}$, in such a way that all closed orbits of the
characteristic foliation outside this set will be non-degenerate.
This can be done by a standard adaptation of the Kupka-Smale type
argument, see eg \cite{Palis-DiMelo}. \end{proof}

 \begin{rem} \label{rem-finite} {\rm Period-action
equation \eqref{eq-per-act-on} guarantees that the periods of all
orbits from $\Pc$ are bounded. Indeed, given an orbit from $\Pc$
of period $T$ we have
$$T \leq \mu \sC + (\sP-\mu \sC)\cdot 1 = \sP\;.$$
Thus for a regular Hamiltonian structure and for any $\e>0$ the
set $\{\g\in\cP|\cA(\g)\leq \sC-\e\}$ is finite.}
\end{rem}

 Consider the functional $\cA_{\cH} (\gamma)= \int_\gamma
\Theta$ defined on the space of smooth loops $S^1 \to V$
representing classes from $\cS$. Define the {\it action spectrum}
$\Spec \cH $ of a Hamiltonian structure $\cH$ as the set of values
of $\cA_{\cH}$ on orbits from $\Pc$. Note that this set can be
interpreted as the set of critical values of $\cA_{\cH}$ lying in
the interval $(-\infty;\sC)$. Indeed, the critical points of
$\cA_{\cH}$ correspond to closed orbits of the characteristic
foliation parameterized in an arbitrary way by the circle $S^1$.
Using period-action equation \eqref{eq-per-act-on} one can show (by
a lengthy, but quite straightforward modification of the standard
argument, cf \cite{HZ}) that even for a non-regular Hamiltonian
structure $\cH$ the set $\Spec \cH \cup \{\sC\}$ is a closed
nowhere dense subset of $\R$. The points from $(0,\sC)
\setminus\Spec\cH$ are called {\it non-critical} values of
$\cA_{\cH}$. For two non-critical values $a,b$ we denote by
$\Pc^{(a;b)}$ the subset of $\Pc$ which consists of periodic
orbits $\g$ with $a<\cA_{\cH}(\gamma)<b$.

\subsubsection{The Conley-Zehnder index and the grading}\label{subsubsec-CZ}

Recall that the Conley-Zehnder index is an integer number
associated to a path of symplectic matrices $A=A(t), t \in [0,1]$,
where $A(0) = \id$ and $A(t)$ does not contain $1$ in its
spectrum, see eg \cite{Robbin-Salamon}. It is denoted by $\CZ
(A)$. Different authors use different conventions on the sign and
normalization (up to an additive constant) of the Conley-Zehnder
index. Our convention is as follows:
\begin{itemize}
\item If a path $A$ of symplectic matrices is generated by a sufficiently small
quadratic Hamiltonian $F$ on $\R^{2n}= \C^n$ then the
Conley-Zehnder index of $A$ equals the Morse index of $F$, that is
the number of negative squares;
\item Given a path $\gamma_1$ and a loop $\gamma_2$ of symplectic
matrices, the Conley-Zehnder index of their concatenation
satisfies
$$\CZ(\gamma_1 \sharp \gamma_2)= \CZ(\gamma_1) -
\mathrm{Maslov}(\gamma_2)$$ (note the minus sign!), where $\mathrm{Maslov}$
stands for the Maslov index of a loop;
\item The Maslov index is normalized by
$$\mathrm{Maslov}(\{e^{2\pi i t}\}_{t\in[0,1]}) = +2\;.$$
\end{itemize}

 The Conley-Zehnder index plays a crucial role in the
definition of a grading on the space $C = \mathrm{Span}_{\Z_2}(\Pc)$
which in turn is a basic ingredient in the definition of
generalized Floer homology. Let us discuss this grading in more
details.\footnote{This material is essentially known to the
experts.} First of all, given a free homotopy class $e \in \cS$
let us denote by $\Pc_e$ the set of all orbits from $\Pc$ which
represent $e$. Set $C_e = \mathrm{Span}_{\Z_2}(\Pc_e)$. Then we have
a decomposition $C = \oplus_{e \in \cS} C_e$.

Fix a class $e \in \cS$, and take any loop $\phi_0\co  S^1 \to V$
representing $e$. Pick a symplectic trivialization $g_0\co 
\phi_0^*(TV/\ell) \to \R^{2n} \times S^1$ of the symplectic vector
bundle $(\phi_0^*(TV/\ell), d(\Theta))$ over $S^1$. Given any
other loop $\phi_1\co  S^1 \to V$ representing $e$, choose a homotopy
$\Phi\co  S^1 \times [0;1] $ between $\phi_0$ and $\phi_1$ with
$\Phi|_{S^1 \times \{i\}} = \phi_i,\; i \in \{0,1\}$. Extend $g_0$
to a symplectic trivialization of $\Phi^*(TV/\ell)$ and denote by
$g_1$ its restriction to $S^1 \times \{1\}$. We get a symplectic
trivialization of the bundle $(\phi_1^*(TV/\ell), d(\Theta))$.
Axiom 1 in the definition of a special Hamiltonian structure
guarantees that $g_1$ does not depend on the choice of the
homotopy $\Phi$ up to multiplication by a contractible loop of
symplectic matrices. We say that trivializations $g_0$ and $g_1$
are equivalent, and call the equivalence class {\it a coherent
trivialization} of the bundle $TV/\ell$ in the class $e$. The set
of coherent trivializations in the class $e$ is denoted by
$\JJ_e$. Importantly, this set carries a natural action of $\Z$.
Indeed, consider a coherent trivialization $j \in \JJ_e$ as above
induced by a pair $(\phi,g)$, where $\phi \co  S^1 \to V$ is a loop
and $g\co \phi^*(TV/\ell) \to \R^{2n}\times S^1$ is a symplectic
trivialization.  \marginpar{\small\it Erratum: replace $2k$
by $-2k$}Take any loop $A(t)$ of symplectic matrices with
the \hypertarget{Corr1}{Maslov index} $2k,\; k \in \Z$. Define a map
$$\bar{A}\co  \R^{2n} \times S^1 \to \R^{2n} \times S^1$$
by $\bar{A}(x,t) = (A(t)x,t)$. Consider a new trivialization $g'=
\bar{A} \circ g$ of $\phi^*(TV/\ell)$, and denote by $j'$ the
induced coherent trivialization. By definition, $j'$ is the result
of the action of $k \in \Z$ on $j$, which is denoted by $j' =
j+k$. One can easily check that this $\Z$--action on $\JJ_e$ is
well defined, transitive and free.

Choose any framing of a regular Hamiltonian structure
$(V,(\Theta))$. Denote by $R^t\co V \to V$ the flow of the
corresponding Reeb vector field. The linearization $R^t_*$ of this
flow acts by symplectic automorphisms of the bundle
$(TV/\ell,d(\Theta))$. Let $p(t)$ be a closed orbit of $R^t$ from
$\Pc_e$ with period $T$. Put $\tilde{p}(t) = p(Tt)$. Every
coherent trivialization $j \in \JJ_e$ gives rise to a symplectic
trivialization $$g \co  {\tilde p}^*(TV/\ell) \to \R^{2n} \times
S^1\;.$$ \marginpar{\small\it Erratum:
replace $g_tR_*^{tT}g_0^{-1}$ by
$(g_tR_*^{tT}g_0^{-1})^{-1}$.}
Denote by $g_t$ its restriction to the fiber over $t \in
S^1$. Define the Conley-Zehnder index $\CZ(p,j)$ of the orbit $p
\in \Pc_e$ with respect to the coherent trivialization $j$ as the
Conley-Zehnder index of the \hypertarget{Corr2}{loop} $g_t R^{tT}_*g_0^{-1}$ of
symplectic matrices. One readily checks that this definition does
not depend on the choice of a special framing. Taking into account
our convention on the Conley-Zehnder index and the definition of
the $\Z$--action on $\JJ_{e}$ we have
\begin{equation}\label{eq-CZshift}
\CZ(p,j+k) = \CZ(p,j)-2k \;\;\hbox{for all}\;\; k \in \Z\;.
\end{equation}
As an immediate consequence we get that for each pair of
orbits $p,p' \in \Pc_e$ {\it the difference} $\CZ(p,j) -\CZ(p',j)$
does not depend on the choice of a coherent trivialization $j \in
\JJ_e$. Thus we often denote it by $\CZ(p)-\CZ(p')$ even though
each term in this expression depends on $j$. We conclude that the
space $C_{e}$ is equipped with {\it a $\Z$--grading up to an
additive shift}. Formula \eqref{eq-CZshift} yields, however, that
the even and the odd parts of $C_e$ are well defined.

 Instead of speaking about {\it a $\Z$--grading up to an additive
shift} we prefer to adopt the following more concrete viewpoint.
Let us call a $\Z$--set any set equipped with a free transitive
$\Z$--action. Given such a set, say $Y$, one introduces the notions
of a $Y$--graded vector space and a $Y$--graded chain complex. The
homology groups of a $Y$--graded chain complex are naturally
$Y$--graded.

Introduce a structure of a $\Z$--set on
$${\II_e}:= \JJ_e \times \{0;1\}$$
as follows: put $(j,0)+1 = (j,1)$ and $(j,1)+1 = (j+1,0)$. For an
element $p \in \Pc_e$ define its degree $\deg(p)\in {\II_e}$ as
$(j,0)$ if $\CZ(p,j) = 0$ and as $(j,1)$ if $\CZ(p,j)=1$. One can
easily check that the degree is well defined and
$$\deg(p') = deg(p) + (\CZ(p')-\CZ(p))$$
for all $p,p' \in \Pc$. Thus we endowed the space $C_e$ with
${\II_e}$--grading. Put $$\II = \{(e,i)\;:\; e \in \cS, i \in
\II_e\}\;.$$ Then $$C = \oplus_{(e,i) \in \II} C_{e,i}\;.$$

\begin{defin} \label{def-gradingset}{\rm
The set $\II$ (respectively, $\II_{e}$) is called {\it the grading
set} of the Hamiltonian structure $\cH$ (respectively, in the
class $e$).}
\end{defin}

In \fullref{subsec-floer-framed} below we shall define
generalized Floer homology of a Hamiltonian structure $\cH$ as the
homology of $(C,d)$ where $d\co  C \to C$ is a differential which
maps $C_{e,i} \to C_{e,i-1}$. Hence the homology will inherit the
same $\II$--grading.

 We write $\JJ_e (\cH)$, $\II_e(\cH)$ and $\II(\cH)$
whenever we wish to emphasize dependence of the set of coherent
trivializations and the grading sets on the Hamiltonian structure
$\cH$.

 \begin{rem}\label{rem-grading-constant} {\rm In the case
when $e = [\mathrm{point} \times S^1]$ the set $\II_e$ has a
distinguished element. Namely, take any point $x \in M$ outside a
sufficiently large compact subset of $M$ where $\Theta = \sC dt -
\sK \alpha$. Consider the loop $\phi(t) = (x,t)$ and notice that
$(\phi^* TV/\ell, d\Theta)$ can be naturally identified with
$(T_xM \times S^1, -\sK d\alpha)$. Denote by $j_0 \in \JJ_e$ the
corresponding coherent trivialization, and put $i_0 = (j_0,0)$.
Identify $\Z$ with $\II_e$ via the map $k \to i_0+k$. With this
identification, $\deg (p) = \CZ(p,j_0)$ for every $p \in \Pc_e$.
Thus we have a canonical $\Z$--grading of $C_e$. This remark is
especially useful when $M$ is simply connected, for instance $M =
\R^{2n}$. Let us mention also that if we wish to restrict
ourselves to periodic orbits from the class $[\mathrm{point} \times
S^1]$ only, we can relax Axiom 1 of special Hamiltonian structures
and require that the first Chern class of $(TV/\ell, d(\Theta))$
vanishes on all spherical $2$--cycles. }
\end{rem}

\begin{rem} \label{rem-grading-Lagrangian}{\rm
Another possibility to introduce a natural $\Z$--grading on the
spaces $C_e$ for all $e$ appears in the case when the bundle
$(TV/\ell, d(\Theta))$ is equipped with a Lagrangian subbundle,
say, $L$. In this case we choose any almost complex structure $J$
on $TV/\ell$ compatible with $d(\Theta)$. Given any loop $\phi\co S^1
\to V$ representing the class $e$, we trivialize $\phi^*(TV/\ell)$
by choosing the homotopically unique $J$--unitary frame in
$\phi^*L\subset\phi^*(TV/\ell)$. This procedure is applied, for
instance, when $(M,\omega)$ is the cotangent bundle of a closed
manifold $X$ equipped with the standard symplectic form, and the
Hamiltonian structure $(V,(\Theta))$ is associated to a
time-dependent Hamiltonian function on $M$ as in \fullref{ex-sympl} above. Here the bundle $(TV/\ell, d(\Theta))$ can
be identified with $(TT^*X,-\omega)$ and the Lagrangian subbundle
$L \subset TT^*X$ is simply formed by the tangent spaces to the
fibers of $T^*X \to X$. }\end{rem}

\subsection{Adjusted almost complex structures}\label{seccomplex}
Next, we introduce a class of almost complex structures compatible
with a framed Hamiltonian structure $ \ocH = ((\Theta),\lambda)$.
These almost complex structures are defined on a manifold $\bV = V
\times \R$.

 Take any starshaped hypersurface $P\subset M$, and
denote by $\eta = \mathrm{Ker}(\alpha|_{TP})$ the contact structure
on $P$. In view of the canonical splittings
$M\setminus\Core(M)=P\times \R_+$ and $\bV = M \times S^1 \times
\R$ we will consider the cut of $\ocH$, denoted by $\xi$, and the
contact plane field $\eta$ on $P$ as subbundles of $T\bV$.

 The group $\R$ acts on $\bV$ by translations $(z,s)
\mapsto (z,s+c)$ where $s,c \in \R$ and $z \in V$.

\begin{defin}\label{def-adj-ac}
{\rm An $\R$--invariant almost complex structure $J$ on $\bV$ is
called {\it adjusted} to the framed Hamiltonian structure $\ocH$
if
\begin{itemize}
\item the cut $\xi$ of $\ocH$ is $J$--invariant, and $d(\Theta)|_\xi$ is compatible with $J$;
\item $ J\frac{\p}{\p s}=R$, where $R$ is the characteristic (Reeb) vector field;
\item the contact plane field $\eta$ is $J$--invariant outside $Z
\times \R$ where $Z \subset M$ is a sufficiently large compact
subset.
\end{itemize}
}
\end{defin}

 An important feature of the almost complex manifold
$(\bV,J)$ with the adjusted $J$ is that it carries a natural
foliation by $J$--holomorphic curves. It has the form $\cT_J=
\cT\times\R$ where $\cT$ is the characteristic foliation of $V$
generated by $R$. We call $\cT_J$ {\it the characteristic
$J$--holomorphic foliation} of $(\bV,J)$.

\begin{rem}\label{remscaling-framing}
{\rm Suppose that $J$ is adjusted to a framed Hamiltonian
structure $(\cH,\lambda)$. Given $c>0$ let us define a new almost
complex structure $J_c$ by setting $J_c|_{\xi}=J|_{\xi}$ and
$J\pps=\frac1c R$. Then $J_c$ is adjusted to the framed
Hamiltonian structure $(\cH,c\lambda)$, and the map
$$(x,s)\mapsto (x,cs)$$
is a biholomorphism
$$(\overline V=V\times\R,J)\to(\overline V, J_c).$$
}
\end{rem}

 The next lemma will allow us to feel comfortable with
pseudo-holomorphic curves in $(\bV,J)$. Recall that a
(co-)oriented hypersurface, say $\Sigma$, in an almost complex
manifold $(W,J)$ is called {\it weakly $J$--convex} if there exists
a 1--form $\sigma$ on $\Sigma$ such that
$$T_{\C}\Sigma = \mathrm{Ker}\sigma,$$
$d\sigma(v, Jv) \geq 0$ for $v \in T_{\C}\Sigma$, and the
orientation of $T_{\C}\Sigma$ determined by the orientation of
$\Sigma$ and its co-orientation by $\sigma$ coincides with its
complex orientation.

 Write $z \mapsto c*z$, $c\in \R_+, z \in M$ for the natural
$\R_+$--action on $M$ associated to the Liouville vector field.
Consider a hypersurface
$$\Sigma_u = (u * P) \times S^1 \times \R \subset \bV\;.$$

\begin{lemma}\label{lmcontact-convex}
For all sufficiently large $u$ the hypersurface $\Sigma_u$ is
weakly $J$--convex in $(\bV,J)$.
\end{lemma}
\begin{proof}
Suppose that $u$ is large enough. Then
$$\Theta = \sC dt -\sK\alpha,\;\xi = \mathrm{Ker}(\sP dt - \sQ\alpha), \;R = \frac{1}{\sP}\ppt\;,$$
where $\sC,\sK,\sP,\sQ$ are the structure constants of $\ocH$. By
definition, $$J\pps = R = \frac{1}{\sP}\ppt\;,$$ and hence the
plane $\zeta:= \mathrm{Span}(\pps,\ppt)$ is $J$--invariant.
Furthermore, $\eta$ is also $J$--invariant. Thus the hyperplane
$$H_u :=\eta \oplus \zeta \subset T\Sigma_u$$ is the complex
tangent space to $\Sigma_u$: $$H_u = T_{\C}\Sigma_u\;.$$
 Note that $H_u = \mathrm{Ker}(\alpha|_{T\Sigma_u})$.
Observe that $\eta$ is contained in $\xi$, since both $dt$ and
$\alpha$ vanish on $\eta$. The 2--form $d\alpha$ coincides with the
restriction of $-\sK^{-1}d\Theta$ to $\Sigma_u$. Therefore, since
$d\Theta|_\xi$ is compatible with $J$, we obtain that $-d\alpha$
restricted to $\eta$ is also compatible with $J$. Finally,
$d\alpha$ vanishes on $\zeta$. We conclude that the hyperplane
$T_{\C}\Sigma_u \subset T\Sigma_u$, which coincides with $H_u$,
and the 1--form $-\alpha$ on $\Sigma_u$ satisfy the definition of
weak $J$--convexity given above. We leave it to the reader to
verify the orientation condition. \end{proof}

\subsection{$J$--holomorphic cylinders in $\bV$
}\label{subsec-cylind}

Let $\ocH$ be a framed Hamiltonian structure. Take an adjusted
almost complex structure $J$ on $\bV:= V\times \R$. Write
$\Upsilon = \R \times S^1$ for the standard cylinder with
coordinates $\sigma \in \R$ and $\tau \in S^1$. We endow
$\Upsilon$ with the standard complex structure by introducing the
complex coordinate $z=\sigma+i\tau$. Given two orbits
$\g_\pm\in\Pc$, we consider $J$--holomorphic cylinders
\begin{equation}
\label{eqdbar} F=(f,\varphi): \Upsilon \to V \times \R, \;\;\;\;
\frac{\p F}{\p \sigma} + J\frac{\p F}{\p \tau} = 0
\end{equation}
with asymptotic boundary conditions
\begin{equation}\label{eqboundary-conditions}
\begin{split}
\varphi(\sigma,\tau)\mathop{\to}\limits_{\sigma\to\pm\infty}&=\pm\infty,\cr
f(\sigma,\tau)\mathop{\to}\limits_{\sigma\to\pm\infty}&=\gamma_\pm(T_\pm\tau),\cr
\end{split}
\end{equation}
where the orbits $\gamma_\pm$ are parameterized by the Reeb field
$R$, and $T_\pm$ are their periods.

 Denote by $\hat{\Mc}(\g_+,\g_-)$ the set of solutions of
this system. The group $\Gamma:= \R \times \R \times S^1$ acts on
$\hat{\Mc}(\g_+,\g_-)$ as follows:
$$(a,b,\vartheta)* F= G, \;\; \text{where} \;\; G(\sigma,\tau) =
(F(\sigma+a,\tau+\vartheta), \varphi(\sigma+a,\tau+\vartheta)+ b
).\;$$ The action is free provided $\g_-\neq \g_+$. This is a
straightforward consequence of the fact that the orbits $\g_-$ and
$\g_+$ represent a homotopy class from $\cS$, and therefore are
{\it simple} (ie not multiply covered). Let us mention also that
if $\g=\g_-=\g_+$ then the solution set $\hat{\Mc}(\g,\g)$
consists of maps of the form
$$(\sigma, \tau) \mapsto (\gamma(T(\tau+ b)),\sigma+a),\;$$
which represent the leaf of the $J$--holomorphic characteristic
foliation of $\bV$ corresponding to $\gamma \times \R$.

\begin{defin}\label{def-J-reg}
{\rm An adjusted almost complex structure $ J$ is called {\it
regular} if at every $J$--holomorphic cylinder which satisfies
\eqref{eqboundary-conditions} the linearized
$\overline\p$--operator is surjective. }
\end{defin}

\begin{prop}\label{propregularity}
Regular almost complex structures are generic, that is they form a
Baire set with respect to a natural topology.
\end{prop}
 The proof of this statement is scattered in the
literature for various classes of framed Hamiltonian structures
(see eg \cite{HWZ,Dr,B2}). The cornerstone of the proof is the
following lemma, which we prove in detail in Appendix A. The
reader is referred to the above-cited papers for the rest of the
argument. Consider a $J$--holomorphic map $F = (f,\varphi)\co 
\Upsilon \to \bV$ which converges to two distinct simple periodic
orbits $\g_{\pm}$ as $\sigma \to \pm\infty$. A point $z \in
\Upsilon$ is called {\it an injectivity point} of $F$ if $d_zf$ is
injective and $f^{-1}(f(z))=z$.

\begin{lemma}\label{lem-somewhereinjective} There is an open dense set in $\Upsilon$
consisting of the points of injectivity.
\end{lemma}
 The proof is given in Appendix A.

Assume that $J$ is a regular almost complex structure.
Then for a pair of distinct periodic orbits $\g_+,\g_- \in \Pc$
the moduli space
$${\Mc}(\g_+,\g_-)=\hat{\Mc}(\g_+,\g_-)/\Gamma\;$$
is a smooth manifold of dimension $(\CZ(\g_+)-\CZ(\g_-))-1$. In
particular, when $\CZ(\g_+)-\CZ(\g_-)=1$ the corresponding moduli
space is zero-dimensional. Moreover, this set is compact, and
hence consists of a finite number of points. The compactness can
be proved along the following lines. First, by using \fullref{lmcontact-convex} above one shows that $J$--holomorphic
cylinders in question project to a compact subset of $V$. Then, as
it was proved in \cite{SFT-compact}, the compactness follows from
the absence of contractible closed orbits of the characteristic
vector field of action $\leq \sC$.

\subsection{Floer homology of framed Hamiltonian structures}
\label{subsec-floer-framed}

Let $\ocH$ be a regular framed Hamiltonian structure with an
adjusted regular almost complex structure $J$. The Floer complex
can be defined in a usual way. Assuming $\CZ(\g_+)-\CZ(\g_-)=1$ we
denote by $\nu(\g_+,\g_-)$ the $\mathrm{mod}\,2$ number of elements of the
finite set $\Mc(\g_+,\g_-)$ and consider the $\Z_2$--space $C $
generated by the orbits from $\Pc $. The boundary operator $d\co 
C\to C $ is defined by
\begin{equation}\label{eqboundary}
d\gamma_+=\mathop{\sum\limits_{\g_-\in\Pc
}}\limits_{\CZ(\g_+)-\CZ(\g_-)=1}\nu(\g_+,\g_-)\g_-\,.
\end{equation}
The space $C$ is graded by elements of the grading set
$\II(\cH)$, see \fullref{def-gradingset}. The differential
$d$ maps $C_{e,i}$ to $C_{e,i-1}$ for every $(e,i) \in \II$.

The proof of the following theorem repeats the standard arguments
in the Floer homology theory; one have to take into account \fullref{lmcontact-convex} and compactness results from
\cite{SFT-compact}.

\begin{thm}\label{thmdtw} $d^2=0$.
\end{thm}

Note that the boundary operator $d$ preserves the
filtration defined by the action functional $\cA_{\cH}$. Hence,
for any $a<b<\mathtt{C}$ we get a chain complex
$\left(C^{(a,b)}(\ocH,J),d\right)$, where
$$C^{(a,b)}=(C\cap\{\cA(\g)<b\})/(C\cap\{\cA(\g)\leq a\})\;.$$ Its
$\II$--graded homology groups
$$\GFH^{(a,b)}(\ocH,J)=\Ker\,
d/{\mathrm{Im}}\, d$$ are called {\it generalized Floer homology}
of the framed Hamiltonian structure $\ocH$. As we will see below
(see \fullref{cor-indep}) the homology is independent of
$J$ and of the choice of a framing.

Suppose that a diffeomorphism $\Phi\co V\to V$ is an equivalence
between two framed Hamiltonian structures $\ocH_1$ and $\ocH_2$.
Let $J_1$ be an almost complex structure adjusted to $\ocH_1$, and
$J_2 = \Phi_*(J_1)$. Then $\Phi$ induces a canonical isomorphism
$$\Phi_\sharp\co \GFH^{(a,b)}(\ocH_1, J_1)\to
\GFH^{(a,b)}(\ocH_2,J_2).$$ Indeed, $\Phi$ sends periodic orbits
to periodic orbits, preserves their action and establishes a
diffeomorphism between the moduli spaces of $J_1$-- and
$J_2$--holomorphic curves. Furthermore, $\Phi$ induces a natural
bijection between the grading sets $\II(\cH_1)$ and $\II(\cH_2)$.
Isomorphism $\Phi_\sharp$ alters the grading of generalized Floer
homology in agreement with this bijection.

The generalized Floer homology are equipped with canonical
{\it varying window homomorpisms}
\begin{equation}\label{eqwindow}
E(a,a')\co \GFH^{(a',b)}\to\GFH^{(a,b)}\;\;\hbox{and}\;\;E(b,b')\co \GFH^{(a,b)}\to
\GFH^{(a,b')}
\end{equation}
for $a'<a$ and $b<b'$. These homomorphisms preserve the grading.

\subsection{Directed concordance of Hamiltonian structures}
\label{concordances}

\subsubsection{Concordance}\label{subsubseq-precob}
Let $\cH_+ = (\Theta_+) $ and $\cH_- = (\Theta_-) $ be two
Hamiltonian structures on $V$. We write $\sC_{\pm},\sK_{\pm}$ for
the corresponding structure constants. The standing assumption of
this subsection is $\sC_+ > \sC_-$.

\begin{defin} \label{def-concord} {\rm {\it A directed
concordance} between a pair of Hamiltonian structures $\cH_-$ and
$\cH_+$ is a pair $(W , (\Pi))$, where $(\Pi)$ is an equivalence
class of $1$--forms on the manifold $W= V\times[a_-,a_+]$ and
$a_- < a_+$ are any real numbers, such that
\begin{itemize}
\item[(i)] $\Pi|_{\p_\pm W}=\Theta_\pm$, where $\p_\pm W = V\times a_\pm $;
\item[(ii)] $d\Pi$ is a symplectic form on $W$;
\item[(iii)] $\Pi= f(s)dt -g(s)\alpha$ outside $(V\setminus Z) \times [a_-;a_+]$,
where $f$ and $g$ are real valued functions with $f'>0$ and $g>0$,
and $Z\subset V$ is a compact set.
\end{itemize}
}
\end{defin}

The conditions $f'>0,g>0$ guarantee that the 2--form
$d\Pi$ is symplectic outside $Z\times[a_-,a_+]$. The assumption
$\sC_+ >\sC_-$ is a consequence of $f'>0$. We call any 1--form
$\Pi$ representing $(\Pi)$ and satisfying condition (iii) {\it the
concordance structure}.

We will sometimes denote a directed concordance as
$\mathrm{Con}(\cH_-,\cH_+)$. Let us emphasize that in this notation
$\cH_-$ is the Hamiltonian structure on the bottom of the
concordance and $\cH_+$ is on its top, where the vertical
direction is oriented by the variable $s$.

\begin{exam} \label{ex-cont-cob}{\rm We work in the
situation described in \fullref{ex-cont} above. Let
$\lambda_{\pm}= F_{\pm}(dt-\alpha)$ be two admissible contact
forms on $V$ with $F_+>F_-$. Let $SV = V \times \R_+$ be the
symplectization of $V$ endowed with the Liouville form $\Pi =
s(dt-\alpha)$. Consider the domain
$$U = \{F_-(x,t) \leq s \leq F_+(x,t)\}\;$$ and define a diffeomorphism
$\varphi\co  V\times[0,1]\to U$ by the formula
\begin{equation}\label{eqscaling-diffeo}
\varphi(x,t,s)=(x,t,F_-(x,t)+ s(F_+(x,t)-F_-(x,t)))\;.
\end{equation}
One readily checks that $(V \times [0,1], (\varphi^*\Pi))$ is a
concordance between the Hamiltonian structures
$\cH_{\lambda_+}=(V,(\lambda_+))$ and
$\cH_{\lambda_-}=(V,(\lambda_-))$. }
\end{exam}

\begin{exam} \label{ex-sympl-cob}{\rm We work in the
situation described in \fullref{ex-sympl} above. Let $F_+$ and
$F_-$ be two time-dependent Hamiltonian functions on $M$ with $F_+
> F_-$. Consider the domain
$$U = \{F_-(x,t) \leq s \leq F_+(x,t)\} \subset V \times \R$$
endowed with the 1--form $\Pi= sdt-\alpha$. Let $\varphi$ be the
diffeomorphism $ V\times[0,1]\to U$ defined by formula
\eqref{eqscaling-diffeo}. Then $(V \times [0,1], (\varphi^*\Pi))$
is a concordance between the Hamiltonian structures $\cH_{F_+}$
and $\cH_{F_-}$. }
\end{exam}

It will be convenient in the sequel to deal with
equivalence classes of concordances. To define it let us stick to
the splitting $V = M \times S^1$.

\begin{defin}\label{def-equiv-conc}{\rm Two concordances
$$(W_1 = V \times [a_-;a_+] , (\Pi_1))\quad\hbox{ and}\quad (W_2 = V \times
[b_-,b_+] , (\Pi_2))$$ are called equivalent if there is a
diffeomorphism $\phi \co W_1\to W_2$ such that
\begin{itemize}
\item[(i)] $(\phi^*\Pi_2)=(\Pi_1)$;
\item[(ii)] $\phi (z,t,a_{\pm}) = (z,t,b_{\pm})$ for all $(z,t) \in M \times
S^1$;
\item[(iii)] there exists a diffeomorphism $u\co  [a_-,a_+] \to [b_-,b_+]$ with
$u(a_{\pm})=b_{\pm}$ and a positive function $v\co  [a_-,a_+] \to
\R_+$ with $v(a_-)=v(a_+) = 1$ such that $\phi(z,t,s) =
(v(s)*z,t,u(s))$ for every point $(z,t)$ lying outside a compact
subset of $V$. Here $c*z$ stands for the natural $\R_+$--action on
$V$ generated by the Liouville field.
\end{itemize}}
\end{defin}

As an immediate consequence of the Darboux theorem, we have the
following useful lemma:

\begin{lemma}[Normal form]\label{lem-normform}
Let $\ocH_\pm = ((\Theta_{\pm}),\lambda_\pm)$ be framed
Hamiltonian structures and let $\Con$ be any concordance between
$\cH_+$ and $\cH_-$. Then $\Con$ is equivalent to a concordance
$(W = V \times [a_+,a_-],(\Pi))$ such that the concordance
structure $\Pi$ satisfies the conditions
\begin{itemize}
\item $ \Pi = \Theta_+ + (s-a_+)\lambda_+ $ on $V \times (a_+-\e,a_+]$
\item $\Pi = \Theta_- + (s-a_-)\lambda_- $ on $V \times [a_-,
a_-+\e)$
\end{itemize}
for some $\e>0$.
\end{lemma}

Let us emphasize that this normal form depends on a
choice of the framings $\lambda_{\pm}$ of $\cH_{\pm}$.

\subsubsection{Concordances and coherent trivializations}
\label{subsubsec-conc-triv}

Let $e \in \cS$ be a free homotopy class of loops $S^1\to V$.
Every directed concordance $\mathrm{Con}(\cH_-,\cH_+)$ induces a
canonical bijection $\kappa_e\co  \JJ_{e}(\cH_+) \to \JJ_{e}(\cH_-)$
between the sets of coherent trivializations in the class $e$
associated to Hamiltonian structures on its top and on its bottom.
Moreover, $\kappa_e$ is equivariant with respect to the canonical
$\Z$--actions on $\JJ_{e}(\cH_+)$ and $\JJ_{e}(\cH_-)$. In fact,
this bijection depends only on the equivalence classes of
concordances. The bijection $\kappa_e$ is defined as follows. Let
$(W = V\times [a_-,a_+] ,(\Pi))$ be our concordance. Denote by
$\ell_{\pm}$ the characteristic line fields of $\cH_{\pm}$. Choose
framings of $\cH_-$ and $\cH_+$ and denote by $R_{\pm}\in
\ell_{\pm}$ the corresponding Reeb vector fields. Set $F_{\pm} =
\mathrm{Span}(\pps,R_{\pm})$, and let $E_{\pm}$ be the
skew-orthogonal complements (in the sense of the symplectic form
$d(\Pi)$) of $F_{\pm}$ at the points of $\p_{\pm}W$. Clearly,
$E_{\pm}$ is naturally identified with $TV/\ell_{\pm}$. Let
$\phi_{\pm}\co S^1 \to V \times \{a_{\pm}\}$ be two loops
representing $e$. Every trivialization $g_{\pm}$ of
$\phi_{\pm}^*E_{\pm}$ gives rise to a uniquely defined
trivialization $h_{\pm}$ of $\phi_{\pm}^*(TW)$. Indeed, at the
points of $\p_{\pm}W$ we have $TW = E_{\pm}\oplus F_{\pm}$, where
the bundle $F_{\pm}$ is equipped with a basis and hence is a
priori trivialized. In the case when $h_+$ and $h_-$
simultaneously extend to a trivialization of $TW$ over a cylinder
in $W$ spanning $\phi_+$ and $\phi_-$, we say that
$\kappa_e([g_+]) = [g_-]$. One readily checks that the map
$\kappa$ is a well defined equivariant bijection. It gives rise to
an equivariant bijection $\chi_e\co  \II_{e}(\cH_+) \to
\II_{e}(\cH_-)$ between the grading sets in the class $e$. Setting
$\chi(e,i) = (e,\chi_e(i))$ we get a bijection $\chi\co \II(\cH_+)
\to \II(\cH_-)$ between the full grading sets.

Assume now that $\gamma_+$ and $\gamma_-$ are closed orbits of
$R_+$ and $R_-$ in the class $e$. Then the difference
$\CZ(\gamma_+,j) - \CZ(\gamma_-, \kappa_e(j))$ does not depend on
the particular choice of a coherent trivialization $j \in
\JJ_e(\cH_+)$. Thus we will denote this difference by
$\CZ(\gamma_+)-\CZ(\gamma_-)$.

In the situation described in \fullref{rem-grading-constant},
that is when $e$ is the class of $\mathrm{point}\times S^1$, the map
$\kappa$ preserves the distinguished element. As far as the case
$M=T^*X$ is concerned (see \fullref{rem-grading-Lagrangian})
the same holds true for all $e \in \cS$ for the concordances which
will arise in the present paper. We will use this without special
mentioning.

\begin{rem} \label{rem-conc-homotopy} {\rm Note that
every homotopy $\Theta_t, \; t \in [0,1]$, between $\cH_0 =
(\Theta_0)$ and $\cH_1 = (\Theta_1)$ through Hamiltonian
structures gives rise to a natural equivariant bijection between
the grading sets of this structures. In order to include this
situation in the context of concordances we say that a concordance
$(V \times [a_-,a_+], (\Pi))$ is {\it homotopy-like } if the
restriction of $\Pi$ to each hypersurface $V \times \{s\}$, where
$s \in [a_-,a_+]$, is a Hamiltonian structure. One can easily
check that for a homotopy-like concordance the bijection
associated to the homotopy $(\Pi|_{V \times \{s\}})$ of
Hamiltonian structures coincides with the bijection $\chi$
introduced in this section. For instance, when $\Theta_1$ is
sufficiently close to $\Theta_0$, the linear segment between these
1--forms passes through Hamiltonian structures (see \fullref{prop-special-stab} and its proof). Therefore, in this case
the full grading sets of $\cH_0$ and $\cH_1$ are naturally
identified. }
\end{rem}

\subsubsection{Adjusted almost complex structures on
concordances}\label{subsubseq-adj-cob}

$\phantom{9}$ Let $\ocH_\pm$ be two\break framed Hamiltonian structures. Choose
adjusted almost-complex structures $J_{\pm}$ for $\ocH_{\pm}$
associated to the splittings $$M \setminus \mathrm{Core}(M) =
P_{\pm}\times \R\;,$$ where $P_{\pm} \subset M$ are starshaped
hypersurfaces. Let $$ \mathrm{Con}= (W = M \times S^1 \times
[a_-,a_+],(\Pi))$$ be a directed concordance between $\cH_+$ and
$\cH_-$. We assume that $W$ lies in $\bV= V \times \R$. Fix real
numbers
$$a''_+ > a'_+ >a_+\;\;\;\text{and}\;\;\;a''_-< a'_-<a_-\;.$$
Choose a smooth family $P_s, s \in (a''_-,a''_+)$, of starshaped
hypersurfaces in $M$ such that $P_s = P_+$ for $s > a'_+$ and $P_s
= P_-$ for $s<a_-$. Denote by $c*z$, where $c \in \R_+$ and $z \in
M$, the natural $\R_+$--action defined by the Liouville field $L$
on $M$. Hypersurfaces $$\Sigma_c= \bigcup_{s\in (a''_-,a''_+)}
(c*P_s) \times S^1 \times \{s\}\,,$$ $c \in \R_+$, foliate the set
$$(M \setminus \mathrm{Core}(M))\times S^1 \times (a''_-,a''_+)\;.$$
Denote by $S$ the vector field on this set obtained as the
projection of $\pps$ to $T\Sigma_c$ along the Liouville field $L$.
Let $f(s),g(s)$ be the functions on $[a_-,a_+]$ such that $f'>0$,
$g>0$ and $\Pi=f(s)dt-g(s)\alpha$ outside a compact subset of $W$
(see \fullref{def-concord}(iii) above). Extend them to the
interval $(a''_-,a''_+)$ keeping $f'>0$ and $g>0$.

Denote by $\eta_s$ the natural contact structure on $P_s$,
considered as a subbundle of $TW$.
\begin{defin} \label{def-adjusted-cob} {\rm An almost
complex structure $I$ on $\bV$ is called {\it adjusted} to the
concordance $\Con$ if it satisfies
\begin{itemize}
\item[(i)] $I=J_{+}$ on $\{s> a'_+\}$ and $I= J_-$ on $\{s<a'_-\}$ ;
\item[(ii)] there exists a 1--form $\wt\Pi$ on $V \times (a''_-,a''_+)$
so that $(\wt\Pi)=(\Pi)$ on $W$ and $d{\wt\Pi}$ is symplectic form
which tames $I$.
\end{itemize}
In addition, we request that there exists a compact subset $Z
\subset V$ such that on $(V \setminus Z) \times (a''_-,a''_+)$ the
following conditions hold:
\begin{itemize}
\item[(iii)]$\wt\Pi = f(s)dt-g(s)\alpha$;
\item[(iv)] the plane field $\eta_s$ is
invariant under $I$;
\item[(v)] $IS = k(s)\ppt$ for some function $k(s)>0$.
\end{itemize}
}
\end{defin}

Conditions (i) and (v) are compatible provided $k(s) =
1/\sP_+$ on $(a'_+,a''_+)$ and $k(s) = 1/\sP_-$ on $(a''_-,a'_-)$,
where $\sP_{\pm}$ are the structure constants of the framings of
$\ocH_{\pm}$.

The space of almost complex structures adjusted to a
concordance is contractible.

Note that the definition of an adjusted almost complex
structure agrees with the notion of equivalence of concordances
introduced above. Let $\Con'=(W',(\Pi'))$ be a concordance which
is equivalent to $\Con=(W,(\Pi))$ via some diffeomorphism $\phi\co W
\to W'$. Assume that $I$ is an almost complex structure on
$\overline V$ adjusted to concordance $\Con$. Consider a piecewise
smooth homeomorphism $\wt\psi\co  \bV \to \bV$ which coincides with
$\phi$ on $W$ and has the form $(z,t,s) \to (z,t,s+e_{\pm})$ on
the connected components of $\bV \setminus W$ for suitable
(unique!) constants $e_{\pm}$. Let $U$ be a sufficiently small
neighborhood of $W$ in $\bV$. Smoothen $\wt\psi$ inside
$U\setminus W$ to a diffeomorphism $\psi\co \bV \to \bV$. The
smoothing can be performed in such a way that
\begin{equation}\label{eq-psiinfty}
\psi(z,t,s) = (v(s)*z,t,u(s))\;\;\; \text{on}\;\;\; (V\setminus K)
\times (a''_-,a''_+)\,,
\end{equation}
where the functions $u$ and $v$ are defined on $(a''_-,a''_+)$ and
$K \subset V$ is a compact subset (see condition (iii) of
\fullref{def-equiv-conc}).

\begin{lemma}\label{lem-equiv-adjust} The push-forward
$I'$ of $I$ under $\psi$ is adjusted to concordance $\Con'$.
\end{lemma}

\begin{proof} The structure $I'$ obviously satisfies
conditions (i) and (ii) of \fullref{def-adjusted-cob},
where one defines $\wt{\Pi}'$ as the push-forward of $\wt{\Pi}$
under $\psi$. Let us verify conditions (iii)-(v) ``at infinity"
using representation \eqref{eq-psiinfty}. Obviously $\wt{\Pi}'$
satisfies (iii). Denote by $\psi_{\infty}$ the restriction of
$\psi$ to $\bV_{\infty}:= (V\setminus K) \times (a''_-,a''_+)$.
Put $b''_{\pm} = u(a''_\pm)$, $P'_s = v(u^{-1}(s))*P_{u^{-1}(s)}$
and
$$\Sigma'_c = \bigcup_{s\in (b''_-,b''_+)}
(c*P_s) \times S^1 \times \{s\}\,$$ where $c \in \R_+$. Then
$\psi_{\infty}(\Sigma_c)=\Sigma'_c$ for $c$ large enough. The
distribution $\eta'_s$, defined as the natural contact structure
on $P'_s$, is the image of $\eta_s$ under $\psi_{\infty}$ which
guarantees (iv).

Finally, writing $L$ for the Liouville vector field on $M$, note
that
\begin{equation} \label{eq-vsp-equiv-adjust}
(\psi_{\infty})_{*}\pps = \gamma(s) \pps + \delta (s)L\;
\end{equation}
with $\gamma(s) > 0$. Denote by $Y$ (respectively, $Y'$) the
projection of $T\bV$ to $T\Sigma_c$ (respectively, to
$T\Sigma'_c$) along $L$ at points of $\bV_{\infty}$ (respectively,
of $\psi_{\infty}(\bV_{\infty})$). Observe that $(\psi_{\infty})_{
*}$ intertwines between $Y$ and $Y'$. Put $S' = Y'\pps$. Then,
applying \eqref{eq-vsp-equiv-adjust} we have
$$I'S' = \frac{1}{\gamma(s)}I'Y'(\psi_{\infty})_{*}\pps
=\frac{1}{\gamma(s)}I'(\psi_{\infty})_{*}Y\pps =
\frac{1}{\gamma(s)}(\psi_{\infty})_{*}IS\;.$$ Using condition (v)
for $I$ we conclude that
$$I'S' = \frac{k(s)}{\gamma(s)}(\psi_{\infty})_{*}\ppt
=\frac{k(s)}{\gamma(s)}\ppt\;$$ which proves that $I'$ satisfies
condition (v). This completes the proof.
\end{proof}

The next result is an analogue of \fullref{lmcontact-convex} in the context of concordances.

\begin{lemma}\label{lmcobord-convex}
Let an almost complex structure $I$ be adjusted to the directed
concordance $(W ,(\Pi))$. Then for a sufficiently large $c>0$ the
hypersurface $\Sigma_c$ is weakly $I$--convex.
\end{lemma}

\begin{proof} Put $\zeta = \mathrm{Span}(S,\ppt)$.
Assume that $c$ is large enough. Then
$$T_{\C}\Sigma_c = \eta_s \oplus \zeta = \mathrm{Ker}(\alpha|_{T\Sigma_c})\;.$$
Clearly, $d\alpha=0$ on $\zeta$. Thus it suffices to show that the
2--form $-d\alpha$ tames $I$ on $\eta_s$. For that purpose, note
that near $\Sigma_c$ one has $\Pi = f(s)dt-g(s)\alpha,$ and hence
$$d\Pi = f'(s)ds\wedge dt - g'(s)ds \wedge \alpha - g(s)
d\alpha.$$ Thus $d\Pi|_{\eta_s} = -g(s)d\alpha|_{\eta_s}$. By
definition, $d\Pi$ tames $I$, and therefore, $-d\alpha$ tames $I$
on $\eta_s$. This yields the desired weak convexity.
\end{proof}

\begin{rem} \label{rem-normform-on}{\rm Let $\ocH_{\pm}$
be two framed Hamiltonian structures, and let $\mathrm{Con} = (V
\times [a_-,a_+],(\Pi))$ be a concordance between them. Assume
that the concordance is in the normal form as in \fullref{lem-normform}. Let $J_{\pm}$ be almost complex structures on
the top and on the bottom adjusted to $\ocH_{\pm}$. Observe that
$d\Pi$ tames $J_+$ for $s \in (a_+-\e,a_+]$ and tames $J_-$ for $s
\in (a_-,a_- + \e]$ if $\e>0$ is small enough. Thus one can choose
an almost complex structure $I$ on $V \times \R$ which is adjusted
to $\mathrm{Con}$, coincides with $J_+$ for $s > a_+-\e$ and
coincides with $J_-$ for $s< a_-+\e$. We say that $I$ {\it
respects the normal form}. }\end{rem}

\subsubsection{Gluing of
concordances}\label{subsubsec-gluingconc}

Suppose that we have three Hamiltonian structures $\cH_1,\cH_2$
and $\cH_3$. Let
$$\mathrm{Con}^{i+1,i}= \Con(\cH_i,\cH_{i+1})= (W^{i+1,i} ,(\Pi^{i+1,i}) )\;\;,i=1,2$$
be directed concordances between $\cH_{i+1}$ and $\cH_{i}$. Up to
equivalence of concordances we can assume that $\Pi^{i+1,i}$ are
in normal forms with respect to some framings $\lambda_i$ on
$\cH_i$ (see \fullref{lem-normform} above) and
$W^{i+1,i}=V\times[a_i,a_{i+1}]$ for $i=1,2$, where $a_1<a_2<a_3$.
Then the formulas
\begin{equation}\label{eqcon-gluing}
\begin{split}
W^{31}=&W^{32}\cup W^{21}=V\times [a_1,a_3];\cr
\Pi^{31}=&\begin{cases} \Pi^{32}&\hbox{on}\;\; W^{32};\cr
\Pi^{21}& \hbox{on}\;\; W^{21}\cr
\end{cases}\cr
\end{split}
\end{equation}
define a smooth directed concordance
$\mathrm{Con}^{31}=(W^{31},(\Pi^{31}))$ between $\cH_3$ and $\cH_1$.
which is called a {\it composition} of $\mathrm{Con}^{21}$ and
$\mathrm{Con}^{32}$. We will denote this operation by
$$\mathrm{Con}^{31} = \mathrm{Con}^{32} \diamond \mathrm{Con}^{21}\;.$$
Note that $\Pi^{31}$ is automatically in a normal form for the
framed Hamiltonian structures $\ocH_3$ and $\ocH_1$. One readily
checks that the equivalence class of a directed concordance
$\mathrm{Con}^{31} $ is determined by the equivalence classes of $
\mathrm{Con}^{32}$ and $ \mathrm{Con}^{21}$.

Assume now that our framed Hamiltonian structures
$\ocH_i$ are equipped with adjusted almost complex structures
$J_i$. Assume further that the concordances $\mathrm{Con}^{21}$ and
$\mathrm{Con}^{32}$ are equipped with adjusted almost complex
structures $I^{21}$ and $I^{32}$ which respect the almost complex
structures on the boundaries and which respect the normal forms of
the concordances (see \fullref{rem-normform-on} above). Define
an almost complex structure $I^{31}$ on $V \times \R$ which
coincides with $I^{21}$ for $s \geq a_2$ and with $I^{32}$ for $s
\leq a_2$. One readily checks that $I^{31}$ is smooth, adjusted to
$\mathrm{Con}^{31}$ and respects $J_1$,$J_3$ and the normal form of
$\mathrm{Con}^{31}$.

\subsection{Directed concordances as morphisms in Floer Homology}\label{sec-cob-morphism}
\subsubsection{Monotonicity homomorphism}\label{secmonotonicity}
Let $\ocH_{\pm}$ be a pair of framed Hamiltonian structures
together with a directed concordance $\mathrm{Con} = (W,(\Pi) )$
between them. Let $\chi \co  \II(\cH_+) \to \II(cH_-)$ be the
bijection between the grading sets on the top and the bottom of
the concordance introduced in \fullref{subsubsec-conc-triv}.
Choose adjusted almost complex structures $J_{\pm}$ and $I$ on
$\bV$. Denote by $\sC_{\pm}$ the corresponding structure constants
and by $\Pc_\pm$ the sets of periodic orbits of the Reeb vector
fields $R_{\pm}$ (see \fullref{subsec-per-char}).

Suppose we are given two orbits $\g_\pm\in\Pc_\pm$. Let
us denote by $\Mc(\Con,I;\g_+,\g_-)$ the moduli space of proper $
I$--holomorphic maps
$$F = (h,\varphi)\co  \Upsilon \to \bV = V \times \R$$ of the cylinder $\Upsilon=\R\times S^1$
which satisfy the following asymptotic conditions:
\begin{equation}\label{eqboundary-conditions-on}
\begin{split}
\varphi(\sigma,\tau)\mathop{\to}\limits_{\sigma\to\pm\infty}&=\pm\infty,\cr
h(\sigma,\tau)\mathop{\to}\limits_{\sigma\to\pm\infty}&=\gamma_\pm(T_\pm\tau),\cr
\end{split}
\end{equation}
where the orbits $\gamma_\pm$ are parameterized by the vector
fields $R_\pm$, and $T_\pm$ are their periods. Similarly to the
cylindrical case we will impose the following generic {\it
regularity} condition: the almost complex structure $I$ is regular
along holomorphic cylinders.

This implies, in particular, that for every pair of
orbits $\g_+,\g_-\in\Pc_\pm $ with $\CZ(\g_+)-\CZ(\g_-)=0$ we have
$$\dim \Mc(\Con,I,\g_+,\g_-)=0.$$ This enables us to define a
{\it monotonicity morphism}
$$\overline{\mon}\co C^{(a,b)}(\ocH_+,J_+)\to C^{(a,b)}(\ocH_-,J_-)$$
by the formula
\begin{equation}\label{eqhomo}
\overline{\mon}(\g_+)=
\mathop{\sum\limits_{\g\in\Pc^{(a,b)}(R_+)}}\limits_{\CZ(\g_+)-\CZ(\g_-)=0}\nu(\g_+,\g_-)\g_-\,,
\end{equation}
where $\nu(\g_+,\g_-)$ is the $\mathrm{mod}\,2$ number of elements of
the finite set $\Mc(\Con,I,\g_+,\g_-)$. This morphism respects the
gradings of the complexes. The next standard fact from the Floer
theory is as follows.

\begin{thm}\label{thmdF}
Under the above assumptions, for any $ a<b < \min(\sC_-,\sC_+)$ and
for all $(e,i) \in \II(\cH_+)$ formula \eqref{eqhomo} defines a
homomorphism of generalized Floer complexes {\rm
$$\overline{\mon}\co \left(C_{e,i}^{(a,b)}(\ocH_+,J_+),d_+\right)\to
\left(C_{\chi(e,i)}^{(a,b)}(\ocH_-,J_-),d_-\right),$$} so that we
have {\rm
$$d_-\circ \overline{\mon}= \overline{\mon}\circ d_+.$$} In particular,
$\overline{\mon}$ defines a homology homomorphism {\rm
$$\mon\co  \GFH_{e,i}^{(a,b)}(\ocH,J_+)\to
\GFH_{\chi(e,i)}^{(a,b)}(\ocH,J_-)\;.$$}
\end{thm}

The homomorphism $\mon$ constructed above a priori
depends on the following data:
\begin{itemize}
\item a concordance $\Con = \Con (\cH_-,\cH_+) = (W,(\Pi))$;
\item an almost complex structure $I$ adjusted to $\Con$ which
respects almost complex structures $J_{\pm}$ adjusted to
$\ocH_{\pm}$.
\end{itemize}
Standard arguments of Floer theory show that the homomorphism
$\mon$ does not change under the following operations:
\begin{itemize}
\item replacing the pair $(\Con,I)$ by a pair $(\Con',I')$
where $\Con'$ is a concordance equivalent to $\Con$ and $I'$ is
the corresponding (see \fullref{lem-equiv-adjust}) almost
complex structure adjusted to $\Con'$;
\item a homotopy of the concordance structure $(\Pi)$ on $W$
through the concordance structures;
\item a homotopy of the adjusted almost complex structure $I$
on $W$ through almost complex structures adjusted to $\Con$ and
respecting $J_{\pm}$.
\end{itemize}

Consider now three framed Hamiltonian structures $\ocH_1,\ocH_2$
and $\ocH_3$ equipped with adjusted almost complex structures
$J_1,J_2$ and $J_3$ respectively. Assume that we are given
directed concordances
$$\mathrm{Con}^{i+1,i}= \Con(\cH_i,\cH_{i+1})\;,\;i=1,2,\;$$
which induce monotonicity morphisms
$$\mon^{i+1,i}\co  \GFH(\ocH_{i+1},J_{i+1}) \to
\GFH(\ocH_{i},J_{i})\;.$$ A standard argument of Floer theory
shows that the gluing $\mathrm{Con}^{32} \diamond \mathrm{Con}^{21}$
of concordances corresponds to the composition $\mon^{32}\circ
\mon^{21}$ of the monotonicity morphisms.

The next result is crucial for our proof of non-squeezing results
in contact geometry. It will enable us to reduce calculations in
contact homology to more traditional calculations in Floer
homology.

\begin{prop} \label{cor-indep} For a regular Hamiltonian structure
the generalized Floer homology does not depend on the choice of a
framing and of an adjusted almost complex structure.
\end{prop}

We start with two auxiliary lemmas whose proof is
elementary and is left to the reader.

\begin{lemma} \label{lem-conc-algebraic} Let $A_i^{\delta}$, where
$\delta \in \{-,0,+\}$ and $i \in \{0,1\}$ be a collection of six
linear spaces. Suppose that we are given eight morphisms between
them such that the following diagram commutes:
\begin{equation}
\label{eq-diagram-conc-algebraic} \xymatrix{ A_0^+ \ar[r] \ar[dr]
& A_1^0 \ar[r] \ar[dr] & A_0^- \\
A_1^+ \ar[r] \ar[ur] & A_0^0 \ar[r] \ar[ur] & A_1^- }
\end{equation}
Assume that the diagonal arrows are isomorphisms. Then all the
arrows are isomorphisms.
\end{lemma}

\begin{lemma} \label{lem-conc-geometric} Let $E$ be a topological
linear space and let $\Lambda \subset E$ be a convex cone. Given
any $\sigma_0,\sigma_1 \in E$ with $\sigma_1-\sigma_0 \in
\mathrm{Interior}(\Lambda)$ and any $\lambda_0,\lambda_1 \in
\Lambda$, there exists a smooth path $\sigma_s$, $s \in [0,1]$,
connecting $\sigma_0$ with $\sigma_1$ such that its derivative
$\dot{\sigma_s}$ satisfies the following conditions:
\begin{itemize}
\item[(i)] $\dot{\sigma_s} = \lambda_0$ for $s$ near $0$ and
$\dot{\sigma_s} = \lambda_1$ for $s$ near $1$;
\item[(ii)] $\dot{\sigma_s} \in \Lambda$ for all $s$.
\end{itemize}
\end{lemma}

Let us mention that it suffices to prove the lemma in
dimension at most $3$ by passing to $E' =
\mathrm{Span}(\sigma_1-\sigma_0,\lambda_0,\lambda_1)$ and $\Lambda'
= \Lambda \cap E'$.

Let $(\Theta)$ be a Hamiltonian structure and let
$\Lambda$ be the convex cone of its framings. Define the set $D$
of 1--forms $\sigma = c \lambda$, where $c \in \R$ and $\lambda \in
\Lambda$ so that $(\Theta + \sigma)$ is again a Hamiltonian
structure. Recall that according to Axiom 3 for any
$\lambda\in\Lambda$ there exists $\mu\geq 0$ such that
$d\lambda=\mu d\Theta$. Hence, we have $-\epsilon \lambda \in D$
if $\mu\epsilon<1$. Then, given $\sigma \in D$ we have $d\sigma =
\mu' d\Theta$ with $1+\mu' > 0$, and hence Hamiltonian structures
$(\Theta)$ and $(\Theta +\sigma)$ share the same characteristic
foliation and the set of framings, and $\sigma +\Lambda \subset D$
for every $\sigma \in D$.

Suppose that we are given $\sigma_0,\sigma_1 \in D$ such that
$\sigma_1-\sigma_1 \in \mathrm{Interior}(\Lambda)$. For any
$\lambda_0,\lambda_1 \in \Lambda$ consider framed Hamiltonian
structures $\ocH_i = ((\Theta+\sigma_i),\lambda_i)$, where $i \in
\{0,1\}$, equipped with adjusted almost complex structures $J_i$.
We wish to associate to these data a monotonicity morphism $$\mon\co 
\GFH (\ocH_1,J_1) \to \GFH(\ocH_0,J_0)$$ associated to a certain
directed concordance $(V \times [0,1], (\Pi))$ between $\cH_0$ and
$\cH_1$ equipped with an adjusted almost complex structure $I$
which respects $J_0$ and $J_1$. Here the form $\Pi$ is defined by
$\Pi = \Theta+\sigma_s$, where the family of 1--forms $\sigma_s$ is
chosen from \fullref{lem-conc-geometric}. Since $\sigma_s \in D$
for all $s$ in view of condition (ii) of the lemma, we have
$d\sigma_s=\mu_s d\Theta$ with $1+\mu_s
> 0$. Then the 2--form
$$d\Pi = (1+\mu_s) d\Theta + ds \wedge \dot{\sigma_s}$$
is symplectic since the 1--form $\dot{\sigma_s}$ lies in $\Lambda$
and hence is positive on the characteristic foliation of
$d\Theta$. Take $\kappa>0$ small enough and define an almost
complex structure $I$ on $V \times \R$ as $J_0$ for $s<\kappa$ and
as $J_1$ for $s > 1-\kappa$. It follows from condition (i) of
\fullref{lem-conc-geometric} that $d\Pi$ tames $I$ near the top
and the bottom of our concordance. Extending $I$ in an arbitrary
way to an adjusted almost complex structure on $V \times \R$ which
is tamed by $d\Pi$, we complete the construction. For the sake of
further reference, we call the constructed morphism {\it the
preparatory morphism}.

\begin{rem}\label{rem-prep-iso}
{\rm Assume that in the above setting $\lambda_0 = \lambda_1$ and
$J_0 = J_1=J$. Then the adjusted almost complex structure $I$ on
the constructed concordance can be chosen as $I = J$. In this case
all $I$--holomorphic cylinders with $0$--dimensional moduli spaces
are of the form $\g \times \R$ where $\g$ is a closed orbit of the
common characteristic foliation of $\cH_0$ and $\cH_1$. Therefore,
our concordance induces the identity morphism in the generalized
Floer chain complexes $CF(\ocH_1,J)$ and $CF(\ocH_0,J)$. Let $a <
b$ be non-critical values of the Hamiltonian structure $(\Theta)$
so that the intervals $(a-\kappa,a+\kappa)$ and
$(b-\kappa;b+\kappa)$ do not intersect the action spectrum
$\Spec((\Theta))$ for some $\kappa>0$ small enough. Suppose now
that the 1--forms $\sigma_0$ and $\sigma_1$ are sufficiently small.
Then the sets $\Pc^{(a,b)}(\cH_1)$ and $\Pc^{(a,b)}(\cH_0)$
consist of the same periodic orbits. Thus our concordance gives an
isomorphism between $\GFH^{(a,b)}(\ocH_1,J)$ and
$\GFH^{(a,b)}(\ocH_0,J)$. }
\end{rem}

\begin{proof}[Proof of \fullref{cor-indep}] Let $a < b$
be non-critical values of the Hamiltonian structure $(\Theta)$.
Let $\lambda_0,\lambda_1$ be two different framings of the
Hamiltonian structure $(\Theta)$, and $J_i$ an almost complex
structure on $V \times \R$ adjusted to the framed Hamiltonian
structure $\ocH^0_i = ((\Theta),\lambda_i)$, $i =0,1$. Pick any
framing $\lambda \in \mathrm{Interior}(\Lambda)$, and take $\epsilon
> 0$ small enough so that $-\e\lambda$ lies in the
set $D$. Consider framed Hamiltonian structures
$$\ocH^{\pm}_i = ((\Theta \pm
\e\lambda),\lambda_i),\;\; i \in \{0,1\}\;,$$ and put
$$A^{\delta}_i := \GFH^{(a,b)} (\ocH^{\delta}_i,J_i),\;\;
\delta \in \{-,0,+\},\;i \in \{0,1\}\;.$$ Consider the diagram
\eqref{eq-diagram-conc-algebraic} whose arrows are defined as the
preparatory morphisms. One can easily check that this diagram is
commutative: in order to prove commutativity of a parallelogram
formed by its arrows, one has to verify that the compositions of
the concordances corresponding to these arrows have homotopic
concordance data. We leave the details to the reader. Decreasing
if necessary $\epsilon$, we get that the diagonal arrows are
isomorphisms, see \fullref{rem-prep-iso}. Hence by \fullref{lem-conc-algebraic} all the arrows are isomorphisms and in
particular the spaces $A^0_0=\GFH^{(a,b)}((\Theta),\lambda_0)$ and
$A^0_1=\GFH^{(a,b)}((\Theta),\lambda_1)$ are isomorphic. Note that
these isomorphisms preserve the grading since all the concordances
involved into their construction are homotopy-like in the sense of
\fullref{rem-conc-homotopy} above. This completes the proof.
\end{proof}

\begin{rem} \label{rem-independ-canonical}{\rm
One can readily check that the isomorphism between generalized
Floer homologies provided by \fullref{cor-indep} is
canonical in the following sense. Let $\lambda_i$, $i \in
\{1,2,3\}$, be three framings of $(\Theta)$ and let $J_i$ be
almost complex structures adjusted to $\ocH_i =
((\Theta),\lambda_i)$. Let
$$F_{i,j} \co  \GFH(\ocH_{j},J_{j}) \to \GFH(\ocH_{i},J_{i})\;, \;\; i,j \in \{1,2,3\},$$
be our isomorphism. Then $F_{13} = F_{12}\circ F_{23}$.}
\end{rem}

It follows from \fullref{cor-indep} and \fullref{rem-independ-canonical} that we can remove the framing and
the almost complex structure from the notation and write \break
$\GFH^{(a,b)}(\cH)$ for the generalized Floer homology of a
Hamiltonian structure $\cH$.

Consider now two regular Hamiltonian structures $\cH_0=
(\Theta_0)$ and $\cH_1 = (\Theta_1)$. Recall that the full grading
sets of sufficiently close Hamiltonian structures can be naturally
identified (see \fullref{rem-conc-homotopy}).

\begin{proposition}\label{propisomorphismon}
Take any two non-critical values $ a<b$ of $\cA_{\cH_0}$. If
$\Theta_1$ is sufficiently close to $\Theta_0$, there is a
grading-preserving isomorphism
$$\GFH^{(a,b)}(\cH_0)\to \GFH^{(a,b)}(\cH_1).$$
\end{proposition}

\begin{proof}
Fix a framing $\lambda_0$ of $\cH_0$ with the framing parameter
$\mu$. Then $\lambda_1 = \lambda_0 + \mu (\Theta_1-\Theta_0)$ is a
framing of $\cH_1$. Take $\epsilon > 0$ small enough so that
$(\Theta_0 -\e \lambda_0)$ is a Hamiltonian structure. If
$\Theta_1$ is sufficiently close to $\Theta_0$, the 1--form
$\Theta_1-\e \lambda_1$ also defines a Hamiltonian structure.
Consider Hamiltonian structures $\cH_i^{\pm} = (\Theta_i \pm
\e\lambda_i)$, where $i \in \{0,1\}$. Re-denote $\cH^{0}_i :=
\cH_i$ and put
$$A^{\delta}_i := \GFH^{(a,b)} (\cH^{\delta}_i),\;\;
\delta \in \{-,0,+\},\;i \in \{0,1\}\;.$$ Our next goal to
describe the arrows in the diagram
\eqref{eq-diagram-conc-algebraic}. The diagonal arrows are defined
as the isomorphisms presented in \fullref{rem-prep-iso}. The
horizontal arrows are associated to the following directed
concordances between the corresponding Hamiltonian structures. The
concordance structure $\Pi$ on each of these concordances
corresponds to the linear segment between the 1--forms defining the
Hamiltonian structures. Let us work out the case of the upper left
horizontal arrow $A^+_0 \to A^0_1$ (all other cases are absolutely
similar): Here the 1--form $\Pi$ on $V \times [0,1]$ is defined as
$$\Pi = \Theta_1 + s(\Theta_0 + \e \lambda_0 - \Theta_1)\;.$$
We claim that $d\Pi$ is symplectic. Indeed, $$d\Pi = (1+\e\mu
s)d\Theta_0 + \e ds\wedge \lambda_0 + ds\wedge (\Theta_0
-\Theta_1) - (1-s)(d\Theta_0 -d\Theta_1)$$ is a small perturbation
of the non-degenerate 2--form $(1+\e\mu s)d\Theta_0 + \e ds\wedge
\lambda_0$ provided that $\Theta_1$ is sufficiently close to
$\Theta_0$, and hence it is non-degenerate as well. The claim
follows.

One readily checks that we got a commutative diagram, and hence
\fullref{lem-conc-algebraic} yields the desired isomorphism. It
preserves the grading since all the concordances involved in the
construction of this isomorphism are homotopy-like in the sense of
\fullref{rem-conc-homotopy} above.
\end{proof}

We complete this section with mentioning that the monotonicity
morphisms behave naturally with respect to diffeomorphisms. Assume
that $\cH^{\pm}_i,\; i = 0,1,$ are four Hamiltonian structures and
$$\Con(\cH^-_0,\cH^+_0)= (W_0,(\Pi_0)) \;\;\;\text{and}\;\;\;
\Con(\cH^-_1,\cH^+_1)= (W_1,(\Pi_1))$$ are two directed
concordances. Let $A \co  W_0 \to W_1$ be a diffeomorphism whose
restrictions to the boundaries $\p_{\pm}W_0$ send $\cH^{\pm}_0$ to
$\cH^{\pm}_1$ and such that $(A^*\Pi_1) = (\Pi_0)$. Then the
following diagram commutes:
\begin{equation}\label{eq-mon-naturality}
\xymatrix{ \GFH(\cH^+_0) \ar[r]^{A^+_\sharp} \ar[d]^{\mon_0} &
\GFH(\cH^+_1) \ar[d]^{\mon_1} \\
\GFH(\cH^-_0) \ar[r]^{A^-_\sharp} & \GFH(\cH^-_1) \ & }
\end{equation}
where the vertical arrows correspond to the monotonicity morphisms
defined by the concordances, and the horizontal arrows stand for
the natural morphisms induced by the restriction of $A$ to
$\p_{\pm}W_0$.

\subsection{$\GFH$ in the absence of regularity}
\label{subsec-absence} Let $\cH= (\Theta)$ be a Hamiltonian
structure, not necessarily regular. For our applications, it would
be convenient to extend the definition of generalized Floer
homology $\GFH^{(a,b)}(\cH)$ to this case. This is done as follows
(in fact, we imitate analogous construction in the context of
usual Floer homology). Assume that $a < b$ are non-critical values
of $\cA_{\cH}$. Take a sufficiently small regular perturbation
$\Theta_0$ of $\Theta$ defining a Hamiltonian structure $\cH_0 =
(\Theta_0)$ . Put
$$\GFH^{(a,b)}(\cH) = \GFH^{(a,b)}(\cH_0)\;$$
\fullref{propisomorphismon} shows that if $\Theta_0$ and
$\Theta_1$ is a pair of sufficiently small such perturbations, the
spaces $\GFH^{(a,b)}(\cH_0)$ and $\GFH^{(a,b)}(\cH_1)$ are
canonically isomorphic and hence the definition above is correct.
Similarly, one extends the monotonicity morphisms to not
necessarily regular Hamiltonian structures.

\subsection{Discussion: non-special Hamiltonian structures}\label{subsec-disc-nonspecial}
Most of the above theory works for stable but not necessarily
special Hamiltonian structures. Namely, suppose Axiom 3 in the
definition of a special structure $\cH=(V,(\Theta))$ is weakened
to simply the requirement of existence of a framing
$\lambda\in\digamma_{\sP,\sQ}$, ie, a $1$--form $\lambda$ such that
$i_R d\lambda=R$ and $\lambda(R)=1$ for a characteristic vector
field $R$ of $(\Theta))$ (see \fullref{secstable} above.)
This leads to the following modifications:
\begin{itemize}
\item The period-action equation \eqref{eq-per-act-on} does not hold, and
what is even worse, one does not necessarily have a bound on the
period of a closed orbit in terms of its action.
\item As a result we cannot prove analogues of \fullref{prop-special-stab} (openness of the stability condition)
and of \fullref{prop-orb-regular}.
\item Even if $\cH$ is regular, it is not any longer true that the sets
$\Pc^{a,b}$ for $a<b<\sC$ are finite.
\item However, assuming that $\cH$ is regular, one can choose an adjusted almost
complex structure $J$ and define the complex
$\left(C^{(a,b)}(\ocH,J),d\right)$ using the same formula
\eqref{eqboundary} as in the special case. Though the complex $
C^{(a,b)}(\ocH,J)$ can be in this case infinitely generated,
compactness results from \cite{SFT-compact} guarantee that the sum
in formula \eqref{eqboundary} contains only finite number of
terms, and hence, $\GFH(\ocH)$ is well defined in this case.
\item Similarly, one can extend to this situation
the definition of monotonicity homomorphisms and prove analogues
of Theorems \ref{thmdF} and \fullref{cor-indep}.
\end{itemize}
For certain applications (which we do not consider in this paper)
it is important to further enlarge the class of admissible
Hamiltonian structures. In particular, one can bypass Axiom 1, by
considering Floer homology with coefficients in an appropriate
Novikov ring. Without much changes one can relax Axiom 2 by
requiring that there are no contractible periodic orbits of
certain Maslov indices (cf \cite{B,Ustilovsky,Yau}. Dropping
Axiom 2 altogether leads to {\it generalized Floer homology
differential algebra} rather than complex, in the spirit of
\cite{SFT}.

\subsection{$\GFH$ for Hamiltonian functions}\label{secFloer}

\subsubsection{(Generalized) Floer homology}
In this section we analyze Floer homology of Hamiltonian
structures $\ocH_{F}$ associated to a generic time-dependent
Hamiltonian function $F(x,t)=F_t(x)$ on $M$ which equals to a
positive constant $\sC$ outside a compact set, see \fullref{ex-sympl} above.

First of all, let us recall the standard definition of
the Floer homology of $F$ (see eg \cite{BPS}). The action
functional $\cA_F$ on the space of smooth loops $\gamma\co  S^1 \to
M$ is given by
$$\cA_F(\gamma) =\int\limits_0^1-\g^*\alpha+F_t(\g(t))dt\;.$$
The values of $\cA_F$ on 1--periodic orbits of the Hamiltonian flow
generated by $F$ form the action spectrum $\Spec (F)$. In what
follows we refer to real numbers $a \in (-\infty;\sC) \setminus
\Spec (F)$ as to {\it non-critical values} of $\cA_F$.

Let $\overline{J}= \{J_t\},\; t\in S^1$, be a family of almost complex
structures compatible with $\omega = d\alpha$. Under certain
generic assumptions, the filtrated Floer homology complex
$C^{(a,b)}(F,\overline{J})$ is generated by 1--periodic orbits
$\overline \g_i$ of the Hamiltonian $F$ with
$$a < \cA_F(\overline\g_i)<b,$$
where $ a, b$ are non-critical values of $\cA_F$.

Denote by $\cS^M$ the set of free homotopy classes of loops $S^1
\to M$. The generators $\overline\g_i$ of
$C^{(a,b)}(F,\overline{J})$ are graded by their free homotopy
classes $\e \in \cS^M$ and by elements of the $\Z$--set
$\II_{\e}^M: = \JJ_{\e}^M \times \{0,1\}$, where $\JJ_{\e}^M$ is
the set of coherent symplectic trivializations of the bundle
$(TM,\omega)$ over loops representing $\e$. The details of
this construction are exactly the same as in \fullref{subsubsec-CZ} above, and we leave them to the reader. \hypertarget{Corr3}{Put}
$$\II^M = \{(\e,i)\;:\; \e \in \cS^M, i \in \II_\e^M\}\;.$$\vspace{-20pt}
\marginpar{\small\it Erratum:
  \hyperlink{Corr3R}{Text to be added here.}}

\begin{defin} \label{def-gradingset-Liouville}{\rm
The set $\II^M$ (respectively, $\II_{\e}^M$) is called {\it the
grading set} of the Liouville manifold $(M,\omega)$ (respectively,
in the class $\e$).}
\end{defin}

In the cases when $\e$ is the class of contractible
loops or when $M$ is endowed with a Lagrangian distribution we get
the canonical $\Z$--grading, as it was explained in Remarks
\ref{rem-grading-constant} and \ref{rem-grading-Lagrangian}.

The boundary operator $d\co C^{(a,b)}(F,\overline{J})\to
C^{(a,b)}(F,\overline{J})$ is defined by the formula
$$d\overline\g_i=\sum\limits_i \nu_{ij}\overline\g_j,$$ where
$\CZ(\overline\g_i)-\CZ(\overline \g_j)=1$, and $\nu_{ij}$ is the
number $(\mod\,2)$ of the trajectories of the reversed gradient
flow $-\nabla\cA_F$ on the loop space, connecting $\overline\g_i$
and $\overline\g_j$. More precisely, $\nu_{ij}$ counts components
of the moduli space $\Mc_{ij}$ of maps $u\co \R\times S^1\to M$ which
satisfy an inhomogeneous Cauchy-Riemann equation
\begin{equation}\label{eqinhomogeneous}
\begin{split}
&\frac{\p u}{\p \sigma} (\sigma,\tau)+J_\tau \frac{\p u}{\p \tau}
(\sigma,\tau) -J_{\tau} \sgrad F_\tau=0\,,\cr
&\lim\limits_{\sigma\to
-\infty}u(\sigma,\tau)=\overline\g_i(\tau),\cr
&\lim\limits_{\sigma\to
+\infty}u(\sigma,\tau)=\overline\g_j(\tau).\cr
\end{split}
\end{equation}
The homology $\FH(F)$ of the complex
$\left(C^{(a,b)}(F,\overline{J}),d\right)$ is independent of
$\overline{J}$ and called the {\it Floer homology of the
Hamiltonian } $F$.

Let us turn now to the generalized Floer homology of the
framed Hamiltonian structure $\ocH_{F}$. Recall that the
underlying Hamiltonian structure is given by the equivalence class
of 1--forms $(\Theta)$ with $\Theta= F dt-\alpha$ on $V$, and the
framing is chosen as $\lambda:= dt$. The 2--form $\Omega$ equals
$dF_t \wedge dt -\omega$. Hence the Reeb vector field equals
$R=\frac{\p}{\p t}+ \mathrm{sgrad}{F_t}$, so that all its closed
orbits are automatically non-contractible. Moreover, every orbit
representing a class from $\cS$ has period $1$.

The cut $\xi$ is the tangent bundle to the fibers $M
\times\{\mathrm{point}\}$. Extend the family of almost complex
structures $-J_t\co TM\to TM$ compatible with $-\omega = \Omega|_\xi
$ ({\bf mind the minus sign!}) to an adjusted almost complex
structure $J$ on $\bV = M \times S^1 \times \R$. This means, in
particular, that $J$ is invariant along the shifts $s \to s+c$,
where $s$ is the coordinate on $\R$, and $J\pps = R$. We assume
that $J$ is regular.

\begin{prop}\label{propFloer-classic}
The graded filtrated Floer complexes $C^{(a,b)}(\ocH_F, J)$ and\break
$C^{(a,b)}(F,\overline{J})$ coincide together with their
differentials. In particular, the (generalized) Floer homology
{\rm $\GFH(\ocH_F)$} and {\rm $\FH(F)$} are canonically
isomorphic.
\end{prop}

\begin{proof}
The complex $C^{(a,b)}(\ocH_F,J)$ is generated by periodic orbits
$\g_i(t)$ of $R$ whose free homotopy classes lie in the set $\cS$.
Each such class, say $e$, has the form $e= \e \oplus [\mathrm{point}
\times S^1]$ where $\e$ is a free homotopy class of loops $S^1 \to
M$. Write $\ell$ for the kernel of $\Omega$ and observe that the
bundle $(TV/\ell, \Omega)$ can be canonically identified with
$(TM,-\omega)$. Therefore, we have a canonical $\Z$--equivariant
bijection between the sets $\JJ_{e}$ and $\JJ_{\e}^M$ of coherent
trivializations, which yields an equivariant bijection between the
grading sets $\II_{e}$ and $\II_{\e}^M$. The orbits $\g_i$ are in
1--1 correspondence with orbits $\overline\g_i$. Namely, each such
orbit $\g_i\co S^1\to M\times S^1$ has the form
$\g_i(t)=(\overline\g_i(t), t)$. \marginpar{\small\it Erratum:
\hyperlink{Corr4R}{Phrase to be ammended.}}Moreover, we have $\deg
(\gamma_i) = \deg(\overline{\gamma_i})$ 
\hypertarget{Corr4}{since} the change of the
sign of the symplectic form does not affect Conley-Zehnder
indices. Furthermore,
$$\cA_{\cH_F}(\g_i)=\int\limits_{\g_i}F_tdt-\alpha=\int\limits_0^1-\g_i^*\alpha+F_t(\g_i(t))dt=
\cA_F(\overline\g_i)\,.$$ Hence, the two complexes have the same
grading and filtration. Now let us compare the differentials of
these complexes. Consider the complex $$(C^{(a,b)}(\ocH_F,J),
d).$$ Recall from \fullref{subsec-floer-framed} that the
coefficient at $\gamma_j$ in the expression for $d\gamma_i$ is
determined by the $0$--dimensional component of the moduli space of
$J$--holomorphic cylinders which are defined as follows. Write
$\Upsilon = \R \times S^1$ for the standard cylinder endowed with
coordinates $\sigma \in \R$ and $\tau \in S^1$. We endow
$\Upsilon$ with the standard complex structure by introducing the
complex coordinate $\sigma+i\tau$. The $\overline\p$--equation
$$\frac{\p\Phi}{\p \sigma} +J \frac{\p\Phi}{\p\tau}=0$$ for a
$J$--holomorphic cylinder $\Phi=(u,A,B)\co  \Upsilon \to M \times S^1
\times \R$ takes the form
\begin{equation}\label{eqFloer-classic}
\begin{split}
&\frac{\p A}{\p \sigma} =-\frac{\p B}{\p \tau},\cr &\frac{\p A}{\p
\tau}=\frac{\p B}{\p \sigma},\cr & \frac{\p u}{\p \sigma}
-J_{A(\sigma,\tau)}\frac{\p u}{\p \tau} + \frac{\p A}{\p
\tau}J_{A(\sigma,\tau)}\sgrad F_{A(\sigma,\tau)} + \frac{\p B}{\p
\tau}\sgrad F_{A(\sigma,\tau)}\,,\cr &\lim\limits_{\sigma\to
-\infty}u(\sigma,\tau)=\overline\g_j(\tau),\cr
&\lim\limits_{\sigma\to
+\infty}u(\sigma,\tau)=\overline\g_i(\tau),\cr
&\lim\limits_{\sigma \to \pm\infty} A(\sigma,\tau)=\tau,\cr
&\lim\limits_{\sigma \to\pm\infty}B(\sigma,\tau) = \pm\infty\;.\cr
\end{split}\end{equation}
The first two equations just say that
$B(\sigma,\tau)+iA(\sigma,\tau)$ is a holomorphic self-map of the
cylinder, which together with the asymptotic boundary conditions
imply that $A(\sigma,\tau)=\tau$ and $B(\sigma,\tau)=\sigma+c$ for
some constant $c$. Without loss of generality we can assume that
$c=0$. Thus the third equation can be rewritten as $$-\frac{\p
u}{\p \sigma} (\sigma,\tau)+J_\tau \frac{\p u}{\p \tau}
(\sigma,\tau) -J_{\tau} \sgrad F_\tau=0\;.$$ The change of
variables $(\sigma,\tau)\mapsto (-\sigma,\tau)$ transforms the
last equation together with the asymptotic boundary conditions to
equation \eqref{eqinhomogeneous} arising in the usual Floer
theory. Hence, we conclude that the Floer complexes
$C^{(a,b)}(\ocH_F,J)$ and $C^{(a,b)}(F,\overline{J})$ coincide
together with their differentials.
\end{proof}

Suppose now that we are given two Hamiltonian functions
$G,F\co M\times S^1\to\R$, $S^1=\R/\Z,$ which are constant at
infinity, and such that $G> F $. Then for any numbers $a<b $ which
are non-critical values for both $\cA_F$ and $\cA_G$, in the Floer
homology theory (see \cite{Floerr,MS}) there is defined a grading
preserving monotonicity homomorphism $\mon_{\FH}\co
\FH_*^{(a,b)}(G)\to \FH_*^{(a,b)}(F)$. On the other hand, the
Hamiltonian structures $\cH_G$ and $\cH_F$ are concordant in a
natural way, see \fullref{ex-sympl-cob} above. This
concordance defines a homomorphism $\mon_{\GFH}\co 
\GFH_*^{(a,b)}(\cH_G)\to \GFH_*^{(a,b)}(\cH_F)$, see \fullref{sec-cob-morphism} above. \fullref{propFloer-classic}
supplies us with identifications between Floer homology of $G$
(resp., of $F$) with generalized Floer homology of $G$ (resp., of
$F$). It is straightforward to check the following:

\begin{prop}\label{propFloer-monoton}
After the identifications of $\GFH$ and $\FH$ the corresponding
monotonicity morphisms coincide: $\mon_{\FH}=\mon_{\GFH}$.
\end{prop}

Moreover, the non-criticality of $a,b$ allows us to
relax the condition $G>F$ to $G\geq F$.

In our calculations below we will use the following
standard fact (see eg \cite{BPS}):
\begin{lemma}\label{lem-homot}
Let $F^{(s)},\; s \in [0,1]$, be a monotone homotopy of functions
constant at infinity: $F^{(s)} \leq F^{(t)}$ for $s <t$. Assume
that $b$ is a non-critical value of $\cA_{F^{(s)}}$ for all $s$.
Then the natural homomorphism {\rm
$\mon\co \FH^{(-\infty;b)}(F^{(1)}) \to \FH^{(-\infty;b)}(F^{(0)})$}
is an isomorphism.
\end{lemma}

\subsubsection{Symplectic homology of
domains}\label{subsubsec-SH-domains} Let $U \subset M$ be a
bounded starshaped domain with boundary $P = \p U$. Define the
action spectrum $\Spec (P) \subset \R$ of $P$ as
$$\Big{\{}-\int_\gamma \alpha\Big{\}}\;,$$
where $\gamma$ runs over all closed orbits of the naturally
oriented characteristic foliation of $P$. Let $\cX=\cX(U)$ be the
set of all non-positive Hamiltonian functions $F\co  M \to \R$
supported in a domain $U$. Define an ``anti-natural" partial order
on $\cX$ by
$$F \preceq G \;\Leftrightarrow \; F \geq G\;,$$ and introduce
symplectic homology of $U$ (see \cite{Floer-Hofer},
\cite{CieFloHof}) as
$$\SH^{(a;b)}(U) := \varinjlim_{F \in \cX}\FH^{(a,b)}(F)\;,$$
where $a < b$ lie in $(-\infty;0) \setminus \Spec (P)$. The
grading of Floer homology induces the grading of symplectic
homology by the grading set $\II^M$ of the Liouville manifold
$(M,\omega)$. With a help of the monotonicity homomorphism in
Floer homology, any inclusion $U_1 \subset U_2$ gives rise to a
natural grading-preserving morphism
$$\SH^{(-\infty;b_1)}(U_1)\to \SH^{(-\infty;b_2)}(U_2),\; b_1 \leq
b_2, $$ of symplectic homologies. The purpose of the next sections
is to extend this construction to domains in the contact manifold
$(V, \mathrm{Ker}(dt -\alpha))$.

\subsection{$\GFH$ for contact forms}\label{seccontact}

In this section we focus on the contact manifold $(V,
\mathrm{Ker}(dt-\alpha)$) where $ V = M \times S^1$, as before.

Recall (see \fullref{ex-cont}) that contact form
$\lambda= F(dt-\alpha)$ is called \emph{admissible} if $F$ equals
a constant $\sC$ at infinity and the Reeb field $R_\lambda$ of
$\lambda$ has no contractible closed orbits of period $\leq \sC$.

Denote by $\cF$ the set of all contact forms on $V$
which equal $\mathrm{const}\cdot (dt-\alpha)$ at infinity for some
positive constant, and by $\cF_{ad} \subset \cF$ the subset of all
admissible forms. With an admissible form $\lambda$ one associates
a framed Hamiltonian structure $\ocH_{\lambda} = (\lambda)$
described in \fullref{ex-cont}. Recall that the same 1--form
$\lambda$ serves as its framing. The cut $\xi$ coincides with the
contact structure on $V$, and the characteristic vector field is
the Reeb field $R_{\lambda}$.

In what follows we will write $\CH^{(a,b)}(\lambda)$ instead of
$\GFH^{(a,b)}(\ocH_{\lambda})$, and call this group the filtrated
{\it cylindrical contact homology} of the contact form $\lambda$.
We tacitly assume that $ a <b < \sC$ and that $a,b$ are
non-critical values for $\cH_{\lambda}$ (for the sake of brevity,
we say that $a,b$ are non-critical values for $\lambda$). Note
that the group $\CH^{(a,b)}(V,\lambda)$ is generated only by
periodic orbits of $R$ which belong to classes from $\cS$, unlike
the more general situation considered in
\cite{SFT},\cite{Ustilovsky}, \cite{Yau}, and \cite{B}. In
particular, multiply covered periodic orbits do not contribute to
this homology. As we already had seen in some applications
considered above (eg, see the proof of \fullref{thm-nonsq} in
\fullref{subsec-ch}) this is good enough for our purposes
because we are able to recover the remaining part of cylindrical
contact homology via a covering trick.

As in the case of Hamiltonian Floer homology, contact homology
$\CH(\lambda)$ is graded by the grading set $\II^M$ of the
Liouville manifold $(M,\omega)$ (see \fullref{def-gradingset-Liouville} above). Indeed, denote by $\ell$
the characteristic line field of the Hamiltonian structure
$\cH_{\lambda} = (\lambda)$. Then the symplectic bundle $(TV/\ell,
d\lambda)$ can be canonically identified with $(\xi, d\lambda)$.
The symplectic structure on this bundle can be homotoped to
$-\omega$ which corresponds to the case $\lambda = dt-\alpha$.
Finally, we deform $(\xi,-\omega)$ to $(TM,-\omega)$ via the
family of subbundles $(\xi_s = \Ker(dt -(1-s)\alpha), -\omega)$,
where $s \in [0;1]$. Furthermore, each class $e \in \cS$ has the
form $\e \oplus [\mathrm{point}\times S^1]$ with $\e \in S^M$. This
readily yields that $\II(\cH_{\lambda}) = \II^M$ (cf the
beginning of the proof of \fullref{propFloer-classic}
above). In the cases when either $\e$ is the class of contractible
loops, or $M=T^*X$ and we pick any class $\e$, the contact
homology are $\Z$--graded in accordance with Remarks
\ref{rem-grading-constant} and \ref{rem-grading-Lagrangian}.

If $\lambda_+ = F\lambda_-$,
where $\lambda_\pm \in \cF$ and the function $F$ satisfies $F > 1$
(resp. $F \geq 1$), then we write $\lambda_+ > \lambda_-$ (resp.
$\lambda_+\geq \lambda_-$). The forms $\lambda_\pm$ can be viewed
as sections of the symplectization $SV = V \times \R_+$ and thus
if $\lambda_+ > \lambda_-$, the part of the symplectization
bounded by the sections $\lambda_\pm$ is parameterized by a
directed concordance between the framed Hamiltonian structures
$\ocH_{\lambda_-}$ and $\ocH_{\lambda_+}$, see \fullref{ex-cont-cob} above. Hence there is defined a {\it
monotonicity homomorphism}
$$\mon(\lambda_+,\lambda_-)\co \CH^{(a,b)}(\lambda_+)\to \CH^{(a,b)}(\lambda_-).$$
In what follows we have to deal with the case $\lambda_+ \geq
\lambda_-$ (note the non-strict inequality). Let us make a
digression to settle this situation.

\begin{rem}\label{remscaling}
{\rm Let $\lambda$ be a contact form on $V$, and let $J$ be any
almost complex structure on $\bV$ adjusted to $\ocH_{\lambda}$.
Fix a real constant $c>0$. The dilation $(z,s) \to (z,cs)$ of
$\bV$ sends $J$ to an almost complex structure $J_c$ which is
adjusted to $\cH_{c\lambda}$, see \fullref{remscaling-framing}
above. Furthermore, each $T$--periodic orbit of the Reeb field of
$\lambda$ is a $cT$--periodic orbit of the Reeb field of
$c\lambda$, and hence the $c\lambda$--action differs from
$\lambda$--action by the factor $c$. Therefore, we get a canonical
isomorphism
$$c_*\co  \CH^{(a,b)}(\lambda)\to\CH^{(ca,cb)}(c\lambda)\;.$$
On the other hand, we have the following commutative diagram,
assuming $c<1$:
\begin{equation}\label{eqscaling-monotonicity}
\xymatrix{ \CH^{(0,ca)}(\lambda) \ar[d]_{\mon(\lambda,c\lambda)}
\ar[r]^{E(ca,a)} &
\CH^{(0,a)}(\lambda)\ar[d]^{\mon(\lambda,c\lambda)} \ar[dl]^{c_*}\\
\CH^{(0,ca)}(c\lambda)\ar[r]^{E(ca,a)} & \CH^{(0,a)}(c\lambda)}
\end{equation}
Indeed, the monotonicity homomorphism $\mon(\lambda,c\lambda)$ is
realized by a symplectic concordance which can be endowed, after
completion, by an almost complex structure which is
biholomorphically equivalent to a translationally invariant almost
complex structure. All index $0$ holomorphic cylinders for this
structure are trivial. A simple analysis of this situation shows
that on the level of complexes, the homomorphisms $c_*$ and
$\mon(\lambda,c\lambda)$ differ only by a varying window
homomorphism. }\end{rem}

This implies the following useful
\begin{prop}\label{propcontact-iso}
Suppose that $c<1$, $a>0$ and that the interval $[ca,a]$ does not
intersect the action spectrum for the form $\lambda$. Then the
monotonicity homomorphism
$$\CH^{(0,a)}(\lambda)\to\CH^{(0,a)}(c\lambda)$$ is an
isomorphism.
\end{prop}

This result enables us to extend the notion of
monotonicity morphism between contact homology of forms
$\lambda_+$ and $\lambda_-$ with $\lambda_+ \geq \lambda_-$. In
fact, we will need just the following particular case of this
construction. Take any real number $a$ with $0<a <
\min(\sC_-,\sC_+)$ which is non-critical for both $\lambda_-$ and
$\lambda_+$. Take $c < 1$ sufficiently close to $1$ so that $
c\lambda_- < \lambda_+$ and the interval $[ca,a]$ does not
intersect the action spectrum of $\lambda_-$. \fullref{propcontact-iso} above provides an isomorphism
$$\CH^{(0,a)}(\lambda_-)\to\CH^{(0,a)}(c\lambda_-).$$
Its composition with the (already defined) monotonicity morphism
$$\CH^{(0,a)}(\lambda_+)\to\CH^{(0,a)}(c\lambda_-)$$
gives the desired monotonicity morphism
$$\mon(\lambda_+,\lambda_-)\co  \CH^{(0,a)}(\lambda_+)\to\CH^{(0,a)}(\lambda_-)\;.$$
It does not depend on the particular choice of the constant $c$.

Given a compactly supported contactomorphism $\Phi$ and
a contact form $\lambda$, we write $\Phi_*\lambda$ for
$(\Phi^*)^{-1}\lambda$ and $\Phi_\sharp$ for the induced
isomorphism $\CH(\lambda) \to \CH(\Phi_*\lambda)$.

\begin{prop}\label{propisotopy-invariance}
Suppose we are given two contact forms $\lambda_+,\lambda_-$ such
that $\lambda_+ \geq \lambda_-$. Let $\Phi^t\co V\to V$, $t\in[0,1]$,
be a compactly supported contact isotopy with $\Phi^0=\id$ and
such that for all $t\in[0,1]$ we have $\Phi^t_*\lambda_+ \geq
\lambda_-$. Then the following diagram commutes:
\[
\xymatrix{ \CH^{(a,b)}(\lambda_+) \ar[d]^{\Phi^t_\sharp}
\ar[r]^{\iota_0} &
\CH^{(a,b)}(\lambda_-) \\
\CH^{(a,b)}(\Phi^t_*\lambda_+) \ar[ur]^{\iota_t} & }
\]
where $\iota_t=\mon(\Phi^t_*\lambda_+,\lambda_-)$.
\end{prop}

\begin{proof}
Put $\lambda_+^t =\Phi^t_*\lambda_+$. It suffices to prove the
proposition for the case when $\lambda_+^t > \lambda_-$ for all
$t$. Write $\Pi$ for the contact form $s(dt-\alpha)$ on the
symplectization $SV$ and denote $W = V \times [0;1]$. Let
$\varphi_t\co W \to SV$ be the embedding given by formula
\eqref{eqscaling-diffeo} of \fullref{ex-cont-cob} above, so
that $(W,(\varphi_t^*\Pi))$ is a directed concordance between
Hamiltonian structures $\cH_{\lambda_+^t}$ and $\cH_{\lambda_-}$.
Note that the image of $\varphi_t$ is the domain cut out from the
symplectization $SV$ by the graphs of the forms $\lambda_+^t$ and
$\lambda_-$. The isotopy $\Phi^t$ lifts to a family of
$\R_+$--equivariant symplectomorphisms $\wt{\Phi^t}$ of $SV$.
Cutting off the Hamiltonian which generates $\wt{\Phi^t}$ outside
a neighborhood of $\mathrm{graph}(\lambda_+^t)$ we get a family of
symplectomorphisms $\Psi^t$ which send $\mathrm{Image}(\varphi_0)$
onto $\mathrm{Image}(\varphi_t)$, coincide with $\wt{\Phi^t}$ near
$\mathrm{graph}(\lambda_+)$ and are equal to the identity near
$\mathrm{graph}(\lambda_-)$. Define a diffeomorphism $A\co  W \to W$ by
$$A = \varphi_t^{-1}\circ \Psi_t \circ \varphi_0\;.$$
Note that the restriction of $A$ to $V \times \{1\}$ sends the
Hamiltonian structure $\cH_{\lambda_+}$ to $\cH_{\lambda_+^t}$,
while the restriction of $A$ to $V \times \{0\}$ is the identity
map. Take now the concordances $(W,(\varphi_0^*\Pi))$ and
$(W,(\varphi_t^*\Pi))$. The proposition follows from commutative
diagram \eqref{eq-mon-naturality}.
\end{proof}

\subsection{Contact homology of domains}\label{seccont-domains}

Let $U \subset V$ be an open subset with compact closure, and let
$\lambda \in \cF_{ad}$ be an admissible form. We denote by
$\cF(U,\lambda)$ the set of contact forms which coincide with
$\lambda$ outside of $U$, and by $\cF_{ad}(U,\lambda)$ its subset
consisting of admissible forms. The set $\cF(U,\lambda)$ is
endowed with the ``anti-natural" partial order $\preceq$ defined
as follows: $$ \lambda'' \preceq \lambda' \Leftrightarrow
\lambda'' \geq \lambda'\;.$$ An increasing sequence in a partially
ordered set is called {\it dominating} if every element of the set
is $\preceq$ than some element of the sequence.

\begin{defin} {\rm An admissible form $\lambda$ on $V$ is called
$U$--{\it basic} if the partially ordered set $\cF(U,\lambda)$
admits a dominating sequence consisting of admissible forms.
Domains $U$ admitting a basic form are called {\it admissible}.}
\end{defin}

This means, in the down to earth language, that there
exists a sequence of admissible forms $\lambda_i \in
\cF_{ad}(U,\lambda)$ such that $\lambda_i/\lambda$ converges to
zero on every compact subset of $U$ as $i \to \infty$.

Let $\lambda$ be a $U$--basic form of an admissible domain $U$.
For $\e>0$ denote by $\cF_{ad}(U,\lambda,\e)$ the set of all forms
from $\cF_{ad}(U,\lambda)$ which do not have $\e$ as a critical
value. Given a dominating sequence in $\cF_{ad}(U,\lambda)$, one
can perturb it to a dominating sequence in
$\cF_{ad}(U,\lambda,\e)$ provided $\e$ is small enough. Define
{\it contact homology} of $U$ as follows:
\begin{equation}\label{eqdomains}\CH(U|\lambda)=:=
\varprojlim_{\e\to 0}\;\varinjlim_{\lambda' \in
\cF_{ad}(U,\lambda,\e)}\CH^{(0,\e)}(\lambda')\;.
\end{equation}
Let $\lambda_i \in \cF_{ad}(U,\lambda)$ be a dominating
sequence. Unveiling the definitions of direct and inverse limits
we see (after, possibly, a small perturbation of $\{\lambda_i\}$)
that every element of $\CH(U|\lambda)$ can be represented as a
sequence
$$\{x_N\}, \;\;\;x_N \in
\CH^{(0;\frac{1}{N})}(\lambda_{i_N}), \;\;\;i_N \to \infty\;,$$
such that for $K
> N$ the images of $x_N$ and $x_K$ under the natural morphisms
coincide in $\CH^{(0;\frac{1}{N})}(\lambda_{j})$ for some $j\geq
\max(i_N,i_K)$.

The grading of contact homology of forms induces the
grading of contact homology of domains by elements of the grading
set $\II^M$ of the Liouville manifold (see \fullref{def-gradingset-Liouville}).

Clearly, the notions introduced above are invariant
under the action of the group $\cG$ of all compactly supported
contactomorphisms of $V$. We sum up this in the following lemma.

\begin{lemma}\label{lmPhi-star}
Let $U $ be an admissible domain with a basic form $\lambda$. Let
$\Phi\in\cG$ be a contactomorphism. Then the domain $\wt
U=\Phi(U)$ is admissible with a basic form
$\wt\lambda:=\Phi_*\lambda$. Moreover, there is a canonical
isomorphism
$$\Phi_\sharp\co \CH(U|\lambda)\to\CH(\wt U|\wt\lambda).$$
\end{lemma}

\begin{proof}
One can define $\CH(U)$ starting from the $U\!$--basic form
$\lambda$, and define $\CH(\Phi(U))$ starting from the basic form
$\wt\lambda=\Phi_*\lambda$. A dominating sequence
$\lambda_i\in\cF_{ad}(U,\lambda)$ is mapped by $\Phi_*$ to a
dominating sequence $\wt\lambda_i\in\cF_{ad}(\Phi(U),\wt\lambda)$.
On the other hand, for any contact form $\lambda_i$ the
contactomorphism $\Phi$ induces an isomorphism
$\Phi_\sharp\co \CH^{(0,\e)}(\lambda_i)\to
\CH^{(0,\e)}(\wt\lambda_i)$. These isomorphisms commute with the
monotonicity morphisms $\mon_{ij}\co \CH^{(0,\e)}(\lambda_i)\to
\CH^{(0,\e)}(\lambda_j)$ and $\widetilde{
\mon}_{ij}\co \CH^{(0,\e)}(\wt\lambda_i)\to
\CH^{(0,\e)}(\wt\lambda_j)$, for $i<j$. Hence $\Phi_\sharp$
descends to the required isomorphism
$$\Phi_\sharp\co \CH(U|\lambda)\to\CH(\wt U|\wt\lambda)$$ between the limits.
\end{proof}

An important example of admissible domains is given by
the next lemma.

\begin{lemma}\label{lmadmissible}
Let $U \subset V$ be a fiberwise starshaped domain. Then the form
$dt-\alpha$ is $U$--basic. In particular, $U$ is admissible.
\end{lemma}

\begin{proof}
We claim that an admissible form $\lambda_0 = dt-\alpha$ is
$U$--basic. Indeed, put $P_t= \p U \cap (M \times\{t\})$. Let us
recall that the complement of $\Core(M)$ can be identified with
the symplectization $SP_t$, and the Liouville form $\alpha$ can be
written as $u\beta_t$ for the contact form $\beta_t=
\alpha|_{TP_t}$. In these coordinates we have $\lambda_0=
dt-u\beta_t$. Take a sequence $G_i$ of non-decreasing positive
functions $G_i\co \R_+ \to \R$ which are constant near $0$, equal $1$
for $u>1$ and such that $G_i$ uniformly converge to $0$ on every
compact subset of $U$. Clearly, the sequence $\lambda_i:=
G_i\lambda_0 \in \cF(U,\lambda_0)$ is dominating. To see that the
forms $\lambda_i$ are admissible note that on
$(M\setminus\mathrm{Core}(M))\times S^1$ we have
$$d\lambda_i|_{M\times \{t\}}=-d(G_iu)\wedge \beta_t-G_iu
d\beta_t\;.$$ Since $\frac{d(G_iu)}{du} > 0$ we conclude that the
restriction of $d\lambda_i$ to the fibers ${M\times \{t\}}$ is
non-degenerate. Hence the Reeb field of $\lambda_i$ is transversal
to the fibers $M\times \{t\}$, and therefore it has no
contractible closed orbits. This proves that all $\lambda_i$ are
admissible.
\end{proof}

\begin{prop}\label{propindependence} Let $U \subset V$ be
an admissible domain. With any two $U$--basic forms $\lambda$ and
$\lambda'$ one can associate an isomorphism
$$q(\lambda', \lambda)\co  \CH(U|\lambda')\to\CH(U|\lambda).$$ This
correspondence is functorial in the sense that for any three
$U$--basic forms $\lambda_1,\lambda_2$ and $\lambda_3$ we have
$$q(\lambda_1, \lambda_3) =q(\lambda_2, \lambda_3)\circ q(\lambda_1, \lambda_2)\,.$$
\end{prop}

\begin{proof}

{\bf Step 1}\qua Let $\lambda'$ and $\lambda$ be two
$U$--basic forms such that $\lambda'>\lambda$. Denote by
$\{\lambda'_i\}$ and $\{\lambda_i\}$ the corresponding dominating
sequences of admissible forms. We claim that for every $i$ there
exists $j$ such that $\lambda'_i > \lambda_j$. Indeed, there
exists a compact $ K \subset U$ such that $\lambda'_i > \lambda_1
$ on $V\setminus K $. Thus $\lambda'_i >\lambda_j $ on $V\setminus
K$ for all $j$ since $\lambda_1 \geq \lambda_j$. Choose now $j$ so
large that $\lambda'_i > \lambda_j $ on $K$. The claim follows.

It follows from the claim that the dominating admissible
sequences $\lambda'_i\in\cF_{ad}(U,\lambda')$ and
$\lambda_i\in\cF_{ad}(U,\lambda)$ can be chosen in such a way that
for all $i$ we have $\lambda'_i > \lambda_i$. By passing to double
limits in the family of monotonicity homomorphisms
$$\mon(\lambda_i',\lambda_i)\co \CH^{(0,\e)}(\lambda'_i)\to \CH^{(0,\e)}(\lambda_i)$$
we get a homomorphism $\CH(U|\lambda')\to \CH(U|\lambda)$ which
will be denoted by $q(\lambda',\lambda)$. Then, given three forms
$\lambda_1,\lambda_2,\lambda_3$ with $\lambda_i > \lambda_j $ for
$i >j$ the formula
$$q(\lambda_1, \lambda_3)=q(\lambda_2, \lambda_3)\circ q(\lambda_1, \lambda_2)\, $$ follows
from the corresponding property for the monotonicity homomorphisms
(see \fullref{thmdF} above).

\medskip
{\bf Step 2}\qua Next, we claim that
$q(\lambda_0,c\lambda_0)$ is an isomorphism for any $U$--basic form
$\lambda_0$ and any $c \in (0,1)$. Note that if
$\lambda\in\cF(U,\lambda_0)$ then $c\lambda\in\cF(U,c\lambda_0)$.
Look at commutative diagram \eqref{eqscaling-monotonicity}:
First, pass to the direct limit as $\lambda \in \cF(U,\lambda_0)$.
Then consider the inverse limit as $a \to 0$. The horizontal
arrows become the identity maps, so the vertical and the diagonal
arrows of the diagram coincide in the double limit. Note that the
vertical arrows correspond to the homomorphism
$q(\lambda_0,c\lambda_0)$. The claim follows from the fact that
the diagonal arrow $c_*$ in the diagram is an isomorphism for all
$a$ and $\lambda$.

\medskip
{\bf Step 3}\qua Now for any two forms $\lambda, \lambda'$
such that $\lambda'>\lambda$ we have for a sufficiently small
constant $c>0$:
$$\frac1c\lambda>\lambda'>\lambda>c\lambda',$$
and hence we have
$$ q(\lambda',\lambda)\circ q(\frac1c\lambda,\lambda')=q(\frac 1c\lambda,\lambda)$$
and
$$q(\lambda,c\lambda')\circ q(\lambda',\lambda) =q(\lambda',c\lambda')\,,$$
which implies that $q(\lambda',\lambda)$ is an isomorphism.
Finally for any two $U$--basic forms $\lambda_0$ and $ \lambda_1$
choose a $U$--basic form $\lambda$ such that $\lambda < \lambda_0$
and $\lambda < \lambda_1$ and put
$$q(\lambda_0,\lambda_1) = q(\lambda_1,\lambda)^{-1}\circ
q(\lambda_0,\lambda)\;.$$ It is straightforward to check that the
definition is independent of $\lambda$ and satisfies the required
functoriality condition.
\end{proof}

Therefore, from now on we can omit the $U$--basic form $\lambda$
from the notation and write
\[
\CH(U)=\CH(U|\lambda).
\]
The monotonicity homomorphism for contact homology of
forms induces an inclusion homomorphism
$\incl(U',U)\co \CH(U')\to\CH(U)$ for admissible domains $U,U' $ with
$U' \subset U$. Indeed, we can choose dominating sequences of
forms $\{\lambda_k\}$ for $U$ and $\{\lambda'_k\}$ for $U'$ so
that $\lambda'_k\geq\lambda_k$. The monotonicity maps
$\mon_k\co \CH^{(0;\e)}(\lambda'_k)\rightarrow
\CH^{(0;\e)}(\lambda_k)$ induce the inclusion homomorphism in the
double limit when $k\to\infty$ and $\e\to 0$.

On the other hand, according to \fullref{lmPhi-star}
a contactomorphism $\Phi\in \cG$ induces an isomorphism
$\Phi_\sharp\co \CH(U)\to\CH(\Phi(U))$.

It is easy to see that contact homology $\CH(U)$ of admissible
domains, inclusion morphisms $\incl(U',U)$ and induced
isomorphisms $\Phi_\sharp$ satisfy conditions
\eqref{eq-functon}-\eqref{eq-functfo} of \fullref{subsec-ch} of
the introduction. We summarize this in the following theorem in
the language of category theory. In order to avoid cumbersome
formulations we will assume that compactly supported
contactomorphisms of $V$ induce grading-preserving morphisms of
contact homology. This happens, for instance, when $M$ is simply
connected, when $\dim M\geq 4$ and $M$ is Weinstein, or when $\dim
M \geq 6$. In the latter two cases we have, maybe after a
perturbation of the Liouville structure, that $\codim(\Core
(M))\geq 2$, and hence every free homotopy class of loops in $M$
has a representative outside an arbitrarily large compact subset.

\begin{thm} \label{thm-ch-functor} Assume that either $M$ is simply connected,
or $\dim M \geq 4$ and $M$ is Weinstein, or $\dim M\geq6$.
Consider a category whose objects are admissible open domains of
$V$ and morphisms correspond to inclusions. Then the
correspondence $U \to \CH(U)$ is a $\cG$--functor of this category
to the category of $\II^M$--graded vector spaces over $\Z_2$, where
$\cG$ is the group of compactly supported contactomorphisms of $V$
and $\II^M$ is the grading set of the Liouville manifold
$(M,\omega)$.
\end{thm}

Since fiberwise starshaped domains are admissible, we
recover results stated in the beginning of \fullref{subsec-ch}
of the introduction.

In the remaining cases, eg, when $M$ is a
non-contractible Liouville surface, the same result holds true if
either we replace the group $\cG$ by its identity component, or
restrict ourselves to contact homology in the class of
contractible loops. In the general case the categorical
formulation is more clumsy and is omitted.

The next theorem and its corollary is not used in the
sequel, but we present them for the completeness of the picture.

\begin{thm}\label{propisotopy-domains}
Let $U$ be an admissible domain, and $\Phi\in\cG$ a
contactomorphism which admits an isotopy $\Phi^t$ connecting
$\Phi^1=\Phi$ with $\Phi^0=\id$ such that $\Supp\Phi^t\subset U$
for all $t\in[0,1]$. Then for any admissible sub-domain $U'\subset
U$ the following diagram commutes:
\begin{equation}\label{eqisotopy-domains}
\xymatrix{ \CH(U') \ar[d]^{\Phi^t_\sharp} \ar[r]^\incl &
\CH(U) \\
\CH(\Phi_t(U')) \ar[ur]^{\incl_t} & }
\end{equation}
where $\incl=\incl(U',U)$ and $\incl_t=\incl(\Phi^t(U'),U)$.
\end{thm}

\begin{proof}
We can choose dominating sequences of forms $\lambda_k$ for $U$
and $\lambda_k'$ for $U'$ such that for all $t\in[0,1]$ we have
$\Phi^t_*\lambda_k'\geq \lambda_k$. Then according to \fullref{propisotopy-invariance} we have a commutative diagram
\[
\xymatrix{ \CH^{(0,\e)}(\lambda_k') \ar[d]^{\Phi^t_\sharp}
\ar[r]^{\iota_0} &
\CH^{(0,\e)}(\lambda_k) \\
\CH^{(0,\e)}(\Phi^t_*\lambda_k') \ar[ur]^{\iota_t} & }
\]
where $\iota_t=\mon(\Phi^t_*\lambda_k',\lambda_k)$ is the
monotonicity homomorphism. Passing to the double limit when
$k\to\infty$ and $\e\to 0$ we get the diagram
\eqref{eqisotopy-domains}.
\end{proof}

Taking $U'=U$ in this theorem, we get the following
result.

\begin{cor} \label{cor-isotopy-invariance}
Assume that a contactomorphism $\Phi \in \cG$ is isotopic to the
identity through contactomorphisms supported in an admissible
domain $U$. Then the induced isomorphism $\Phi_\sharp\co  \CH(U) \to
\CH(U)$ is the identity map.
\end{cor}

\section{Calculations with contact homology} \label{sec-calcul}

Here we complete the proof of several ``hard" results stated in
\fullref{sec-intro}. First of all we express contact homology
of the prequantization of a symplectic domain in terms of its
filtrated symplectic homology, see \fullref{thm-calcul-split}
and \fullref{propfunctoriality} below which together form
a slightly more explicit version of \fullref{thm-chsh}. Our
approach is based on generalized Floer homology theory developed
in \fullref{secflavors}. We apply this result to calculations
with contact homology. In particular, we prove Theorems
\ref{thm-ch-ellip} and \ref{thm-ball-morph} on contact homology of
prequantizations of ellipsoids and balls in $\R^{2n}$.
Furthermore, we study contact homology of prequantizations of unit
ball bundles of closed manifolds in terms of cohomology of free
loop spaces. As a result, we prove \fullref{thmcont-elements}
on orderability of spaces of contact elements.

\subsection{Contact homology of split
domains}\label{subseq-calcul-split} Let $(M^{2n},\omega,L)$ be a
Liouville manifold with the Liouville form $\alpha = i_L\omega$.
We will work on the manifold $V = M \times S^1$, write $t
\;(\mathrm{mod} 1)$ for the coordinate on $S^1=\R/\Z$, and orient
$V$ by the volume form $(-\omega)^n \wedge dt$. Consider the
contact structure $\xi = \mathrm{Ker}(dt-\alpha)$ on $V$. Recall
that the notion of contact characteristic foliation of a
hypersurface in a contact manifold was defined in \fullref{subsec-geom-ze}.

\begin{defin}\label{def-non-resonant} {\rm A hypersurface $\Sigma$ with smooth boundary
in a contact manifold is called {\it non-resonant} if it is
transversal to the contact structure and the contact
characteristic foliation of $\Sigma$ has no closed orbits.}
\end{defin}

\begin{example}\label{exnon-resonant}
{\rm Let $P \subset \R^{2n}$ be any starshaped hypersurface. Then
the hypersurface $\wh{P}= P \times S^1$ is non-resonant if and
only if the symplectic characteristic foliation on $P$ (equipped
with the natural orientation) has no closed orbits $\gamma$ with
$\int\limits_\gamma\alpha=1$. With the notation of \fullref{subsubsec-SH-domains} this condition means that $-1 \notin
\spec(P)$.}
\end{example}

\begin{thm} \label{thm-calcul-split} Let $U \subset M$ be an open
domain with compact closure and smooth starshaped boundary $P=\p U
$. Suppose that the hypersurface $P \times S^1$ is non-resonant.
Then there exists a grading-preserving isomorphism
$$\Psi_U\co \CH(U \times S^1) \to
\SH^{(-\infty;-1)}(U)\;.$$
\end{thm}

Note that both contact and symplectic homologies in this
theorem are graded by the grading set $\II^M$ of the Liouville
manifold $(M,\omega)$.

\begin{proof}
Recall that we have the canonical decomposition $M = SP
\sqcup \Core(M)$ (see formula \eqref{eq-decomposition} of \fullref{secLiouville}). Consider the contact form $\beta =
\alpha|_{TP}$ on $P$ and fix the identification $SP = P \times
\R_+$ associated to $\beta$ (see \fullref{subsec-geom-ze}
above). We use coordinates $(x,u) \in P \times \R_+$ on $SP$ and
extend $u$ by $0$ to $\Core(M)$. With this notation we have $U =
\{u <1\}$.

\medskip
{\bf Step I\qua From contact to Floer homology}

Denote by $\cH$ the set of continuous Hamiltonians $H\co M \to
(0,+\infty)$ with the following properties:
\begin{itemize}
\item On $SP=P \times (0;+\infty)$ the Hamiltonian $H$ is a piecewise
smooth function of one variable $u$, that is $H = H(u)$;
\item $H(u) \equiv c$ for $u \in
(0,a)$ and $H(u) \equiv 1$ for $u\geq b$, where $c>0$ and $0 < a
<b <1$ are some constants depending on $H$;
\item At each smooth point $u \in (0,+\infty)$,
\begin{equation}\label{eq-Legendre}
H(u)-uH'(u)>0\;,
\end{equation}
where $H'$ stands for the derivative of $H$. Note that the
expression in the left hand side of this inequality can be
interpreted as the vertical coordinate of the intersection point
between the tangent line to $\Graph (H)$ and the vertical axis;
\item $H \equiv c$ on $\mathrm{Core}(M)= M \setminus SP$.
\end{itemize}

Given $H \in \cH$ as above, define a function $\phi_H\co 
\R_+ \to \R_+$ by $$\phi_H(u) = \frac{u}{H(u)}.$$ One readily
checks that $\phi_H$ is a homeomorphism of $\R_+$ whose inverse
can be written as $\phi_{\bar{H}}$ for some function $\bar{H}\in
\cH$. Moreover, $\phi_H$ is a diffeomorphism when $H$ is smooth.
Define a map $\Phi_H\co  V \to V$ as follows: put $\Phi_H(x,u,t) =
(x,\phi_H(u),t)$ on $P \times(0,+\infty) \times S^1$ and extend it
by $(z,t) \to (L^{-\log c}z,t)$ to $\mathrm{Core}(M) \times S^1$.
Again, $\Phi_H$ is a homeomorphism of $V$ with $\Phi_H^{-1} =
\Phi_{\bar{H}}$ and $\Phi_H$ is a diffeomorphism when $H$ is
smooth.

Put $\lambda_H:= \frac{dt-\alpha}{H}$ and $\sigma_H :=
\bar{H}dt-\alpha$. The key property of $\Phi_H$ is given by the
following

\begin{lemma} \label{lem-transform}
$\Phi_H^*\lambda_H = \sigma_H$ for every smooth $H \in \cH$.
\end{lemma}

We deduce from the lemma that for every $\epsilon \in (0,1)$
\begin{equation} \label{eq-calcul-on}
\CH^{(0;\epsilon)}(\lambda_H) = \CH^{(0;\epsilon)}(\sigma_H) =
\GFH^{(0;\epsilon)}(V, (\sigma_H))= \FH^{(0,\epsilon)}(\bar{H}),
\end{equation}
where the right equality follows from \fullref{propFloer-classic}. Furthermore, these identifications
preserve the grading since the diffeomorphism $\Phi_H$ is isotopic
to the identity. Note that $\lambda_H$ is a contact form
corresponding to the contact structure $\xi$. Hence equation
\eqref{eq-calcul-on} links together contact homology and
Hamiltonian Floer homology. Let us emphasize that the main purpose
of generalized Floer homology theory developed in \fullref{secflavors} is to justify this link. Eventually, this link
will enable us to complete the proof of the theorem.

\begin{proof}[Proof of \fullref{lem-transform}] Clearly it
suffices to verify the statement of the lemma on $SP \times S^1$.
Recall that with respect to the splitting $SP = P \times
(0,+\infty)$ the Liouville form $\alpha$ can be written as
$u\beta$ where $\beta = \alpha|_{TP}$. Put $v = \phi_H (u)$ and
observe that
$$\frac{dt-u\beta}{H(u)} = \frac{dt}{H(\phi_{\bar{H}}(v))} -
v\beta\;.$$ On the other hand,
$$\frac{\phi_{\bar{H}}(v)}{H(\phi_{\bar{H}}(v))} = v,$$
which yields \begin{equation} \label{eqtransform}
H(\phi_{\bar{H}}(v)) = \frac{\phi_{\bar{H}}(v)}{v}
=\frac{1}{\bar{H}(v)}\;.
\end{equation}
Therefore,
$$\frac{dt-u\beta}{H(u)} = \bar{H}(v)dt-v\beta,$$
as required.\end{proof}

One readily checks that the transform $H \to \bar{H}$ is
continuous in the uniform topology and anti-monotone:
\begin{equation}\label{eq-antimonotone}
H_1 \leq H_2 \Rightarrow \bar{H}_1 \geq \bar{H}_2 \;.
\end{equation}
Moreover, this transform preserves the class of piecewise linear
functions from $\cH$. We will work out this property for the
following special class of functions. For $0<a<b<1$ and $c>1$
consider a piecewise linear function $F_{a,b,c}\in \cH$ which
equals $c$ on $(0,a)$, equals $1$ on $(b,+\infty)$ and which is
linear on $[a,b]$ (see \fullref{Figon}).
\begin{figure}[ht!]
\cl{\includegraphics{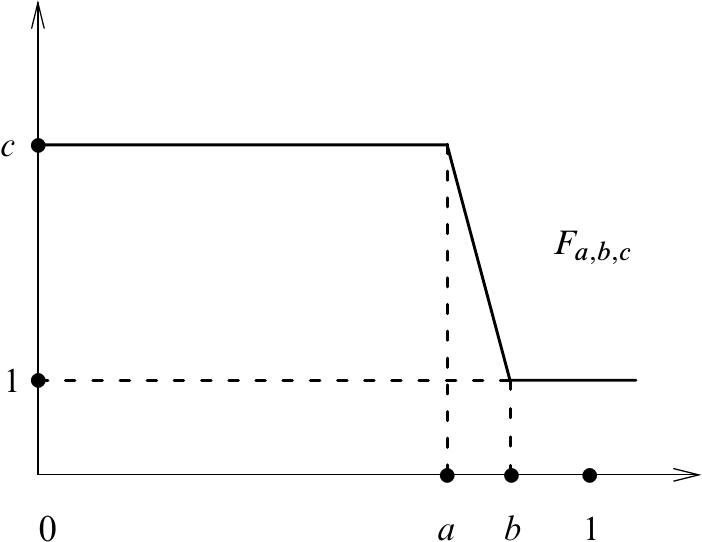}}
\caption{}\label{Figon}
\end{figure}
A straightforward calculation shows that
\begin{equation}\label{eq-Fbar}
\bar{F}_{a,b,c} = F_{a/c,b,1/c}\;.
\end{equation}
For $0 <\mu < \nu < 1$ and $\kappa < 0$ let us consider a
piecewise linear function $G_{\mu,\nu,\kappa}$ which equals
$\kappa$ on $[0,\mu]$, equals $0$ on $[\nu,+\infty]$ and which is
linear on $[\mu,\nu]$. With this notation
\begin{equation}\label{eq-Fbar-on}
\bar{F}_{a,b,c}-1 = G_{a/c,b,-1+1/c}\;,
\end{equation}
see \fullref{Figtw}.
\begin{figure}[ht!]
\cl{\includegraphics{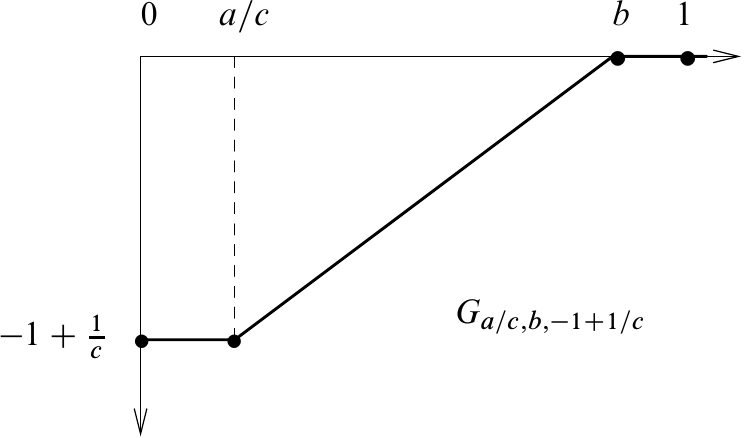}}
\caption{}\label{Figtw}
\end{figure}

Note that $G_{a/c,b,-1+1/c}$ is a non-positive Hamiltonian
supported in $U$. Such Hamiltonians appear in the definition of
symplectic homology of $U$ presented in \fullref{subsubsec-SH-domains} . A troubling point, of course, is that
the function $G_{a/c,b,-1+1/c}$ is not smooth. We are going to
take care of this right now.

\medskip
{\bf Step II\qua A smoothing procedure} 

In what follows we are going to work with Floer homology of
Hamiltonians $G_{\mu,\nu,\kappa}$, which is understood as Floer
homology of a smooth function which approximates $G_{\mu,\nu,\kappa}$
in the uniform topology. A careful choice of the smooth approximations
will enable us to perform precise homological calculations. Let us
explain the details of the approximation procedure we use.

Let $\cL := PT\R^2$ be the projectivization of the tangent
bundle of the plane $\R^2$. Define the set of {\it generalized
tangent lines} $\cT_{\mu,\nu,\kappa} \subset \cL$ to $\Graph
(G_{\mu,\nu,\kappa})$ as follows:
$$\cT_{\mu,\nu,\kappa} = \{(l,(u,G_{\mu,\nu,\kappa}(u))\}$$
so that
\begin{itemize}
\item $l$ is the usual tangent line for $u \notin\{\mu,\nu\}$;
\item $l$ lies above $\Graph (G_{\mu,\nu,\kappa})$ for $u=\nu$;
\item $l$ lies below $\Graph (G_{\mu,\nu,\kappa})$ for $u=\mu$.
\end{itemize}

The Euclidean metric on $\R^2$ gives rise to a natural
metric on $\cL$. The proof of the next elementary lemma is
straightforward.

\begin{lemma}\label{lem-approx}
For a small enough $\delta>0$ there exists a smooth function
$H=H_{\mu,\nu,\kappa}$ on $[0,+\infty)$ with the following
properties (see \fullref{Figth} ):
\begin{itemize}
\item $H(0) = \kappa$, $H'(0) = 0$;
\item $H(u)=0$ for all $u \geq 1$;
\item $H(u) < 0,H'(u)>0, H(u)-uH'(u) < 0$ for all $u \in
(0,1)$;
\item The tangent bundle $T(\Graph (H))$, considered as a
subset of $\cL$, lies in the $\delta$--neighborhood of
$\cT_{\mu,\nu,\kappa}$.
\end{itemize}
Moreover, the $\delta$--approximation $H_{\mu,\nu,\kappa}$ can be
chosen to depend continuously on the parameters
$\mu,\nu,\kappa,\delta$.
\end{lemma}
\begin{figure}[ht!]
\cl{\includegraphics{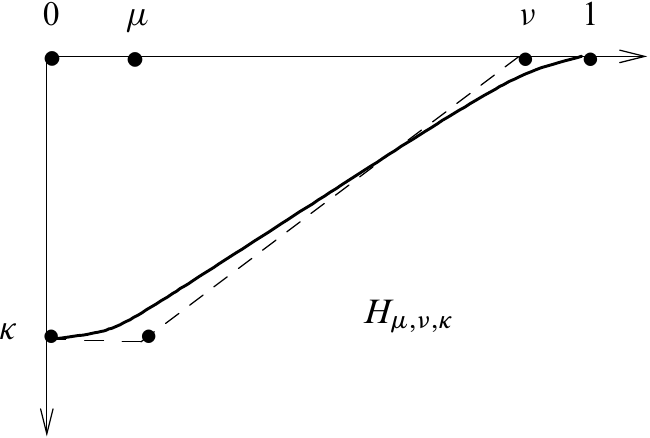}}
\caption{}\label{Figth}
\end{figure}

\medskip
{\bf Step III\qua Homological calculations} 

Take $\epsilon > 0$ small enough. Choose $a,b,c,\delta$ depending on
$\epsilon$ so that $a$ and $b$ are sufficiently close to $1$, $c$ is
large enough and $\delta$ is small enough. We assume that
\begin{equation} \label{eq-eps-del}
\epsilon > \frac{2}{c}\;.
\end{equation}
Let $H$ be a $\delta$--approximation of $G_{a/c,b,-1+1/c}$, as in
\fullref{lem-approx}. The goal of this Step is to establish the
following formula:
\begin{equation} \label{eq-calcul-tw}
\FH^{(-1;-1+\epsilon)}(H) = \FH^{(-\infty,-1+\epsilon)}(H) =
\SH^{(-\infty; -1+\epsilon)}(U)\;.
\end{equation}
To start with, let us recall (see equation \eqref{sgrad}
above) that the Hamiltonian vector field of the Hamiltonian
$H=H(u)$ is given by $H'(u)R$ on $SP= P \times (0,+\infty)$ and it
vanishes outside this domain, where $R$ stands for the Reeb vector
field of $\beta$. Hence its 1--periodic orbits are either the fixed
points lying in $M \setminus (P\times (0,1))$, or the pairs
$(\gamma,u)$ where $u \in (0,1)$ and $\gamma$ is a closed orbit of
the Reeb field $R$ of period $H'(u)$. In the latter case $-H'(u)
\in \spec(P)$. The symplectic action of an orbit corresponding to
a fixed point is simply the value of the Hamiltonian at this
point. The symplectic action of an orbit of the type $(\gamma,u)$
equals $H(u)-uH'(u)$. Recall that this quantity can be interpreted
as the vertical coordinate of the intersection point between the
tangent line to $\Graph (H)$ and the vertical axis.

Since $P \times S^1$ is assumed to be non-resonant, there exists
$\tau
>0$ such that
\begin{equation}\label{eq-forbid-int}
\spec(P) \cap (-1-\tau,-1+\tau) = \emptyset.
\end{equation}
Our next claim is that the actions of all 1--periodic orbits of $H$
are greater than $-1$. Indeed, look at \fullref{Figtw}. The
claim will readily follow from the fact that any line from
$\cT_{a/c,b,-1+1/c}$ intersects the vertical axis above $-1$.
Clearly, it suffices to check this for the line which passes
through $(b,0)$ and $(a/c,-1+1/c)$. Calculating, we get that it
intersects the vertical axis at a point $(0,y)$ with
$$y = -\frac{bc-b}{bc-a}.$$
Since $a<b$ we have $y > -1$ as required, and the claim follows.

The claim yields $$\FH^{(-1,-1+\epsilon)}(H) =
\FH^{(-\infty,-1+\epsilon)}(H),$$ which proves the left equality
in formula \eqref{eq-calcul-tw}. To show that $$
\FH^{(-\infty,-1+\epsilon)}(H) = \SH^{(-\infty, -1+\epsilon)}(U)$$
it suffices to produce a monotone homotopy $H^{(s)},\; s \geq 0,$
such that $H^{(0)} = H$, the family $H^{(s)}$ is dominating with
respect to ``anti-natural" order $\preceq$ on $\cF$, and,
crucially, for every $s$ the point $-1+\epsilon$ does not lie in
the action spectrum of $H^{(s)}$ (see \fullref{lem-homot}).

To construct $H^{(s)}$, we start with a monotone homotopy
$G^{(s)}=G_{\mu(s),\nu(s),\kappa(s)}$, $s \geq 0,$ of the function
$G^{(0)}=G_{a/c,b,-1+1/c}$ through a dominating family, see \fullref{Figfo}. Next, we take $H^{(s)}$ as a smooth
$\delta(s)$--approximation of $G^{(s)}$ provided by \fullref{lem-approx} with $\delta(s)$ small enough.

\begin{figure}[ht!]
\cl{\includegraphics[width=0.83\hsize]{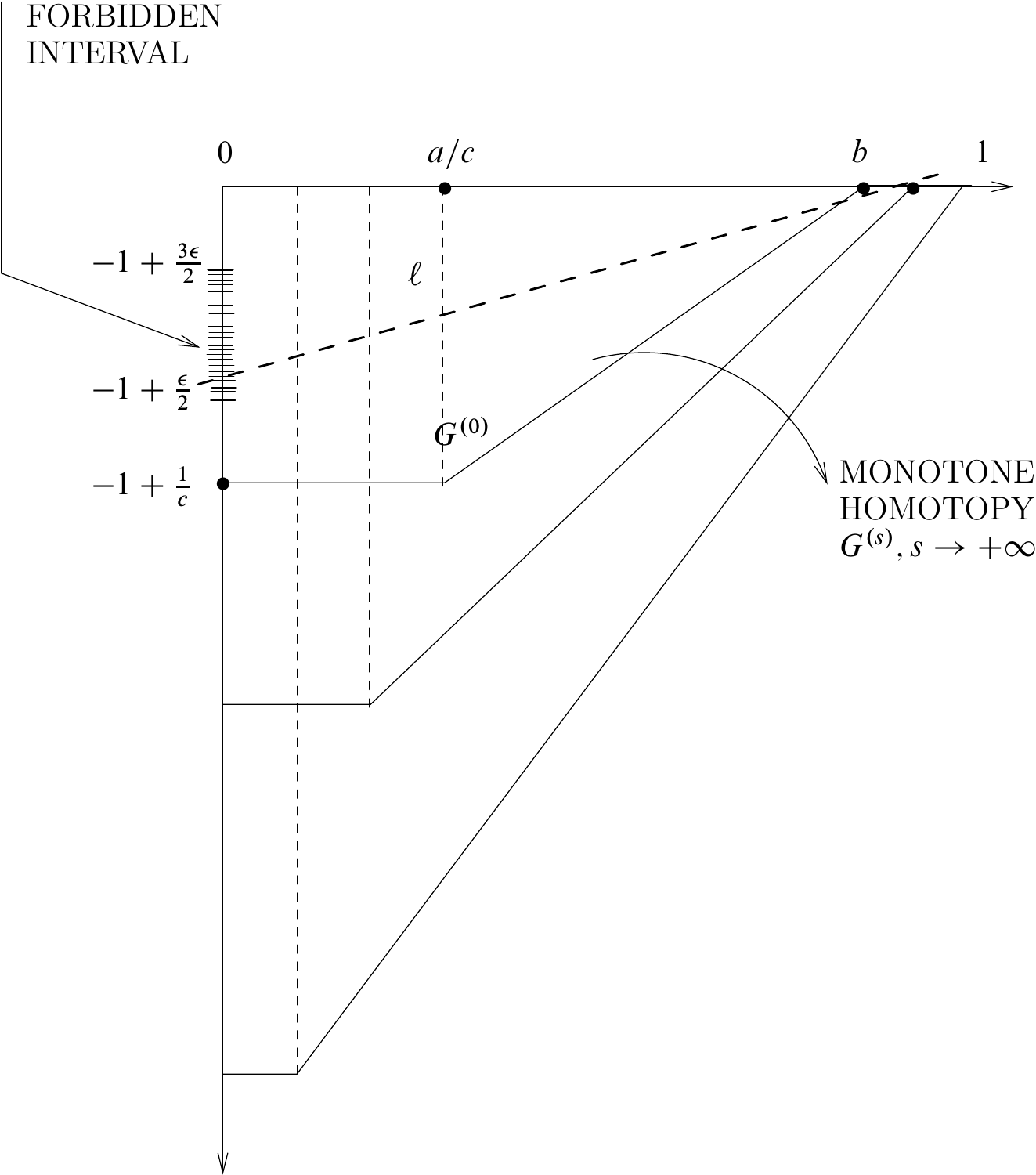}}
\caption{}\label{Figfo}
\end{figure}

We claim that for all $s$ no line $l \in
\cT_{\mu(s),\nu(s),\kappa(s)}$ whose slope lies in $-\spec(P)$ can
pass through the ``forbidden interval"
$[-1+\epsilon/2,-1+3\epsilon/2]$ of the vertical axis. Indeed, any
$l \in \cT_{\mu(s),\nu(s),\kappa(s)}$ passing through the
``forbidden interval" has the slope close to $1$. In particular,
this slope lies in $(1-\tau;1+\tau)$ for $\epsilon$ small enough
and $b$ sufficiently close to $1$. The claim follows from
condition \eqref{eq-forbid-int}. Hence we completed the proof of
formula \eqref{eq-calcul-tw}.

\medskip
{\bf The end of the proof} 

Denote $\wh U = U \times S^1$. Recall that the form $dt-\alpha$ is
$\wh U$--basic, see \fullref{lmadmissible} above. Consider a smooth
approximation of the function $F_{a,b,c}$ (see Step I) given by
$K_{a,b,c} = \overline{H_{a/c,b,-1+1/c}+1}$ (here the bar stands for
the transform described in Step I). Using formulas
\eqref{eq-calcul-on} and \eqref{eq-calcul-tw} we get that
\begin{align*}
\varinjlim_{\lambda\in\cF(\wh U,
dt-\alpha)}\!\!\!\CH^{(0,\epsilon)}(\lambda)& = \lim_{\substack{c \to
+\infty \\ a,b \to 1}} \CH^{(0,\epsilon)} (\lambda_{K_{a,b,c}})\\
& =
\lim_{\substack{c \to +\infty \\ a,b \to 1}} \FH^{(0;\epsilon)}
(H_{a/c,b,-1+1/c}+1) \cr &= \lim_{\substack{c \to +\infty \\ a,b
\to 1}}\FH^{(-1,-1+\epsilon)}(H_{a/c,b,-1+1/c}) =
\SH^{(-\infty,-1+\epsilon)}(U)\;,
\end{align*}
where the limits are understood as direct limits. Finally, we pass
to the inverse limit and see that
$$\CH(\wh U) = \varprojlim_{\e\to 0}\;
\varinjlim_{\lambda\in\cF(\wh U,
dt-\alpha)}\CH^{(0,\epsilon)}(\lambda)=$$$$\varprojlim_{\e\to 0}\;
\SH^{(-\infty,-1+\epsilon)}(U )= \SH^{(-\infty,-1)}(U)\;.$$ In the
last equality we use that
$$\SH^{(-\infty,-1+\epsilon)}(U)= \SH^{(-\infty,-1)}(U
)$$ since $-1 \notin \spec(P)$. This completes the proof.
\end{proof}
\begin{prop}\label{propfunctoriality}
The isomorphism $\Psi_U\co \CH(U \times S^1 ) =
\SH^{(-\infty,-1)}(U)\;$ constructed in \fullref{thm-calcul-split} is functorial in the following sense.
Consider two starshaped domains $U_1 \subset U_2 \subset M$ with
smooth non-resonant boundaries. Let $I \co \CH(U_1 \times S^1 )\to
\CH(U_2 \times S^1)$ and $\imath\co \SH^{(-\infty;-1)}(U_1)\to
\SH^{(-\infty,-1)}(U_2)$ be the inclusion homomorphisms in contact
and symplectic homologies. Then
\begin{equation} \label{eq-funct}
\Psi_{U_2} \circ I = \imath \circ \Psi_{U_1}\;.
\end{equation}
\end{prop}
\begin{proof}
Let us re-examine the isomorphism between $\CH(\lambda_{H})$ and
$\FH(\bar H)$ given by formula \eqref{eq-calcul-on} in Step I of
the proof of \fullref{thm-calcul-split}. Consider the
following manifolds:
\begin{itemize}
\item the symplectization $SV =V \times \R_+$ of $(V,\xi)$ equipped with the
Liouville form $\Pi = s(dt-\alpha)$;
\item the manifold $V \times \R$ endowed with the 1--form
$\Pi' = sdt - \alpha$.
\end{itemize}
Define an embedding $A \co SV\to V\times \R$ by the formula
$$F(z,t,s)=(L^{\log s}(z),t,s), \;\;\text{where}\;\;\; (z,t) \in V = M \times S^1, s\in\R_+\;.
$$
Note that $A^*\Pi' = \Pi$. Furthermore, on the subset $$P\times
\R_+ \times S^1 \times \R_+\subset SV$$ the map $A$ takes the form
$$(x,u,t,s)\mapsto (x,su,t,s).$$ Therefore, given a Hamiltonian $H \in
\cH$, the graph of the contact form
$\lambda_{H}:=\frac1H(dt-\alpha)$, viewed as a section of the
symplectization $SV \to V$, is mapped by $A$ to the hypersurface
$\{s=\bar{H}(u)\} \subset V \times \R$. Indeed, this is equivalent
to the identity
$$\bar{H}\left(\frac{u}{H(u)}\right)=\frac1{H(u)},$$
which is dual to \eqref{eqtransform}.

 Given two Hamiltonians $H_1>H_2$ from $\cH$, we observe that the
domain $Y \subset SV$ between the graphs of the forms
$\lambda_{H_2}$ and $\lambda_{H_1}$ is mapped by $A$ to the domain
$Y' \subset V \times \R$ between the graphs of Hamiltonians
$\bar{H}_2$ and $\bar{H}_1$. But the domains $Y$ and $Y'$
correspond to directed cobordisms which are used in Sections
\ref{seccontact} and \ref{secFloer} in order to define the
monotonicity maps in contact and Floer homology theories. Hence,
the isomorphism between $\CH(\lambda_{H})$ and $\FH(\bar H)$ is
functorial with respect to the monotonicity morphisms. It is
straightforward to check that this functorial correspondence
survives all the limits in the construction of $\Psi_U$.
\end{proof}

\subsection{Ellipsoids and balls} \label{subsec-ellip}

{\bf Note}\qua {\em All symplectic and contact homology in this section are
$\Z$--graded (see \fullref{rem-grading-constant} above).}

\begin{proof}[Proof of \fullref{thm-ch-ellip}] Here we
calculate the contact homology of the domain $ \wh{E}(N,R)=E(N,R)
\times S^1 \subset \C^n \times S^1$, where $E(N,R)$ is the
ellipsoid
\begin{equation} \label{eq-ellip-on}
E(N,R):= \{\pi |z_1|^2 + \frac{\pi}{N} \sum_{i=2}^n |z_i|^2 < R\},
\; N \in \N, \frac{1}{R} \notin \N\;.
\end{equation}
Note that $P:= \partial E(N,R)$ is a convex hypersurface in
$\C^n$. The action spectrum of $P$ equals $-R\N$, hence the
condition $\frac{1}{R}\notin \N$ is exactly the non-resonant
condition. Therefore one can apply \fullref{thm-calcul-split}.
Together with an obvious rescaling it yields
\begin{equation}\label{eq-calcul-ellip-tw}
\CH_*(\wh{E}(N,R)) =\SH_*^{(-\infty,-1)}(E(N,R)) =
\SH_*^{(-\infty,-\frac{1}{R})}(E(N,1)). \end{equation}
Next we calculate symplectic homology   $
\SH_*^{(-\infty,-\frac{1}{R})}(E(N,1))$ (see \cite{FHW} for a
calculation of another version of symplectic homology for
ellipsoids). Define a piecewise linear function
$$F_{a,c}\co  [0;+\infty) \to \R, \;\; a \in (0,1), c<0$$
as follows: $F_{a,c} \equiv 0$ on $[a,+\infty)$, $F_{a,c}(0) = c$
and $F_{a,c}$ is linear on $[0,a]$. We consider $F_{a,c}(u)$ as a
Hamiltonian on $\C^n$ with
$$u =\pi |z_1|^2 + \frac{\pi}{N} \sum_{i=2}^n |z_i|^2\;.$$
In order to deal with Floer homology of $F_{a,c}$ we approximate
it by a smooth function, exactly as in Step II of the proof of
\fullref{thm-calcul-split} above. We leave the details of the
approximation argument to the reader. In what follows $\cT_{a,c}$
stands for the set of generalized tangent lines to $\Graph
F_{a,c}$ which consists of the following points $(l,(u,F_{a,c}(u))
\in PT\R^2$:
\begin{itemize}
\item $l$ is the usual tangent line when $u \neq a$;
\item $l$ lies above $\Graph
(F_{a,c})$ when $u = a$.
\end{itemize}
We claim that
\begin{equation} \label{eq-calc-ellip-on}
\SH_*^{(-\infty,-\frac{1}{R})}(E(N,1)) =
\FH_*^{(-\infty,-\frac{1}{R})}(F_{a,-\frac{1}{R} - \delta})
\end{equation} provided $a$ is sufficiently close to $1$ and $\delta$
is small enough. Indeed, consider a monotone homotopy of
$F_{a,-\frac{1}{R} - \delta}$ through a dominating family of
functions $F_{a(s),c(s)},\; s \geq 0$ (see \fullref{Figfi}) and
note that no line from $\cT_{a(s),c(s)}$ which passes through a
$\frac{\delta}{2}$--neighborhood of $-1/R$ on the vertical axis has
an integral slope. The claim follows from \fullref{lem-homot}.

\begin{figure}[ht!]
\cl{\includegraphics[width= 2.7in]{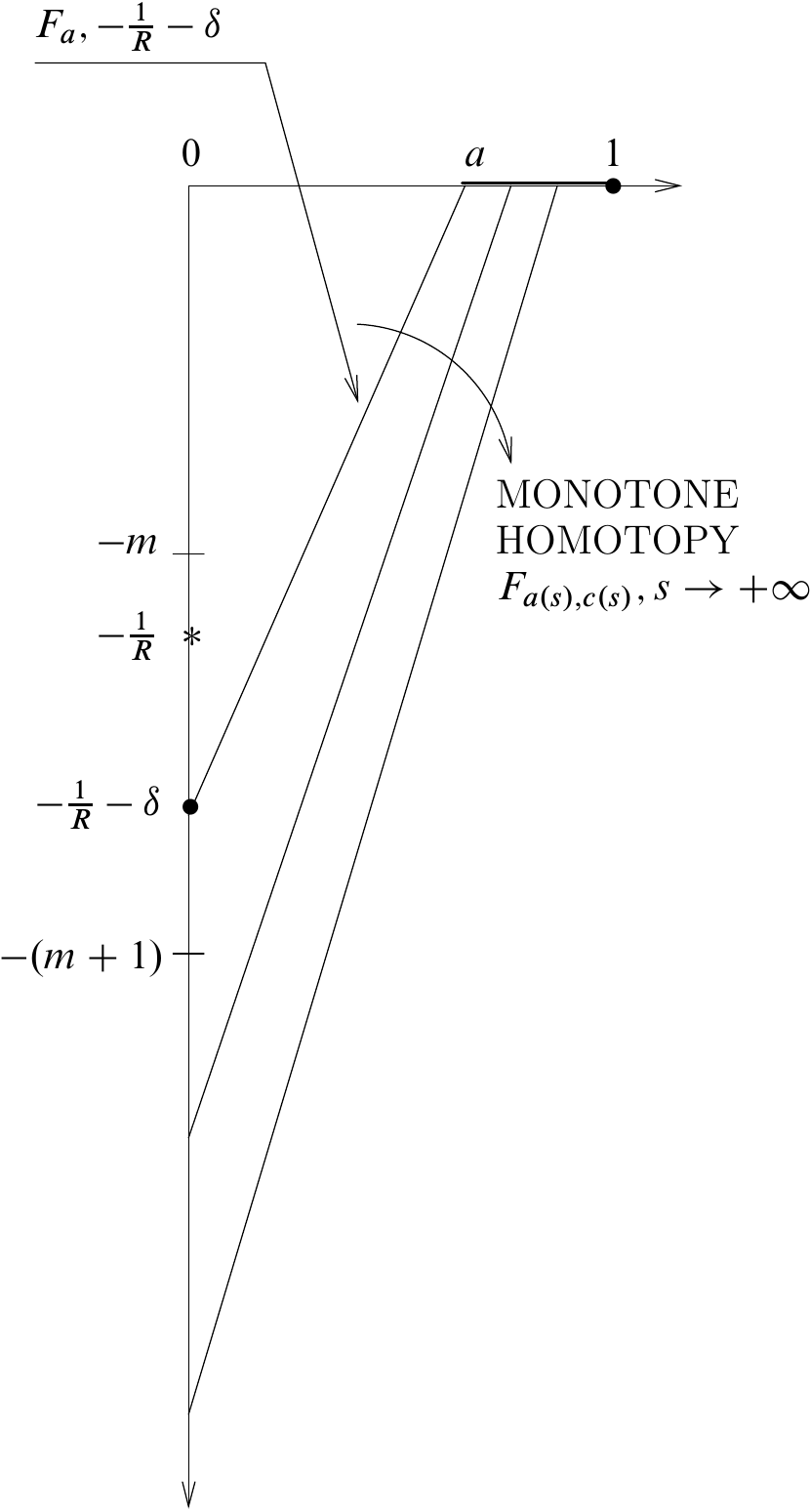}}
\caption{}\label{Figfi}
\end{figure}

Finally, note that the only 1--periodic orbit of
$F_{a,-\frac{1}{R} - \delta}$ with action $\leq -\frac{1}{R}$
corresponds to the fixed point at $u=0$. In the neighborhood of
this point,
$$F_{a,-\frac{1}{R}
- \delta} (z)= \mathrm{const} + \pi
\Big{[}\frac{1}{R}\Big{]}|z_1|^2+
\pi\Big{[}\frac{1}{NR}\Big{]}\sum_{i=2}^n |z_i|^2 + \pi H(z),$$
where $H(z) = \sum_{i=1}^n c_i|z_i|^2$ with $c_i \in (0,1)$.
Therefore, according to our convention (see \fullref{subsubsec-CZ}), the Conley-Zehnder index of our fixed point
equals
$$k(N,R) = -2\Big{
[}\frac{1}{R}\Big{]}-2(n-1)\Big{[}\frac{1}{NR}\Big{]}.$$ We
conclude that $\SH_i^{(-\infty,-1/R)}(E(N,1))$ equals $\Z_2$ for
$i = k(N,R)$ and vanishes otherwise. The theorem follows now from
equation \eqref{eq-calcul-ellip-tw}. \end{proof}

\begin{proof}[Proof of \fullref{thm-ball-morph}] Assume that
$$\frac{1}{k} < R_1 <R_2 < \frac{1}{k-1}$$ for some $k \in \N$.
Using the functoriality equation \eqref{eq-funct} we see that it
suffices to show that the inclusion $B^{2n}(R_1) \hookrightarrow
B^{2n}(R_2)$ induces an isomorphism of symplectic homologies
$$ \SH^{(-\infty;-1)}(B^{2n}(R_1)) \to
\SH^{(-\infty;-1)}(B^{2n}(R_2)).$$ Consider the family of balls
$B_s:= B^{2n}(s)$, $s \in [R_1, R_2]$. The desired statement
follows from the fact that $-1 \notin \spec(B_s)= -s\N$ for all $s
\in [R_1, R_2]$ in view of our assumptions on $R_1$ and $R_2$.
\end{proof}

\subsection{Cotangent bundles and free loop spaces}\label{seccotangent}

{\bf Note}\qua {\em All symplectic and contact homology in this section are
$\Z$--graded (see \fullref{rem-grading-Lagrangian} above). All
singular (co)homology are with $\Z_2$--coefficients.}

Let $X$ be a closed manifold, and let $T^*X$ be the
cotangent bundle of $X$ with the standard symplectic form
$\omega=d\alpha$ where $\alpha = pdq$. Consider the contact
manifold $(V = T^*X \times S^1, \xi = \mathrm{Ker}(\alpha-dt))$.
Choose a Riemannian metric on $X$ and set $U_r=\{|p|<r\}$. Thus
$U_r$ is the ball bundle associated with $T^*X$. Denote $\wh U_r =
U_r \times S^1 \subset V$. Write $\Lambda$ for the length spectrum
of the metric. Note that the hypersurface $\p \wh U_r$ is non
resonant if and only if $r^{-1} \notin \Lambda$. For the sake of
brevity, will shall call such values of $r$ non-resonant.

 Let us denote by $\cL X$ the free loop space of $X$, by
$S\co \cL X\to\R$ the action functional
$$S(\g)=\frac12\int\limits_0^1\|\frac{d\g}{dt}\|^2 dt\,,$$
and by $\cL^aX$ the set $\{S\leq a\}\subset \cL X$. We will use
below the following theorem due to J. Weber \cite{Weber} which is
based on earlier works \cite{Viterbo,Abondandolo-Schwarz,Salamon-Weber}.

\begin{thm}\label{thmloop-space}
For each non-resonant $r>0$ there exists an isomorphism
\begin{equation}
L_r\co \SH^{(-\infty, -1)}(U_r)\to H^*(\cL^{\frac1{2r^2}} X)\;.
\end{equation}
Furthermore, if $0< r<r'$ are non-resonant, the following diagram
commutes:
\[
\xymatrix{\SH^{(-\infty, -1)}(U_r) \ar[r]^{L_r} \ar[d]^{I} &
H^*(\cL^{\frac1{2r^2}} X) \ar[d]^{j^*} \\
\SH^{(-\infty, -1)}(U_{r'})\ar[r]^{L_{r'}} &
H^*(\cL^{\frac1{2(r')^2}} X) }
\]
where $I$ and $j^*$ are respective inclusion homomorphisms.
\end{thm}

\fullref{thm-calcul-split} then implies:
\begin{cor}\label{corloops}
For each non-resonant $r>0$ there exists an isomorphism
\begin{equation}
\wh L_r\co \CH(\wh U_r)\to H^*(\cL^{\frac1{2r^2}} X)
\end{equation}
such that if $r<r'$ the following diagram commutes
\[
\xymatrix{\CH(\wh U_r) \ar[r]^{L_r} \ar[d]^{I} &
H^*(\cL^{\frac1{2r^2}} X) \ar[d]^{j^*} \\
\CH(\wh U_{r'})\ar[r]^{L_{r'}} & H^*(\cL^{\frac1{2(r')^2}} X)\,.
}
\]
\end{cor}

\begin{proof}[Proof of \fullref{thmcont-elements}]

Observe that every finite covering $X \to X'$ lifts to a finite
covering $\PP_+T^*X \to \PP_+T^*X'$ which is a local
contactomorphism. Thus a positive contractible loop in
$\Cont_0(\PP_+T^*X')$ lifts to a positive contractible loop in
$\Cont_0(\PP_+T^*X)$. Therefore if $\PP_+T^*X'$ is non-orderable
then $\PP_+T^*X$ is non-orderable as well. Hence, the case when
$\pi_1(X)$ is finite can be reduced to the case when $X$ is simply
connected.

Consider natural morphisms
$$\tau_a \co  H^*(\cL X)\to H^*(\cL^aX)\;,$$
$$\tau_{E,a}\co  H^*(\cL^E X)\to H^*(\cL^aX)\;,\;\; E >a\;,$$
and denote by $\beta(a) $ and $\beta(E,a)$ the dimensions of their
images, respectively. Note that all these functions are
non-decreasing in $a$ and non-increasing in $E$.

\begin{lemma}\label{lem-dimensions-tw}
Assume that $X$ is simply connected, or that $\pi_1(X)$ has
infinite number of conjugacy classes. Then
\begin{itemize}
\item[(i)] $\beta(a)\mathop{\to}\limits_{a\to\infty}\infty\;;$
\item[(ii)] For every $a > 0$ there exists $E(a)$ such that
$\beta(E,a) = \beta(a)$ for all $E \geq E(a)$.
\end{itemize}
\end{lemma}
The proof is presented at the end of this section.

We continue the proof of \fullref{thmcont-elements}.
According to our covering argument we can assume that $X$ is
either simply connected, or $\pi_1(X)$ has infinitely many
conjugacy classes.

In view of \fullref{lem-dimensions-tw}(i) one can choose an
arbitrarily small $r$ such that $\beta(a)$ jumps when crossing
$b=\frac1{2r^2}$, that is for any sufficiently small
$\delta_1,\delta_2>0$ we have
\begin{equation} \label{eq-jump}
\beta(b-\delta_1)<\beta(b+\delta_2).
\end{equation}
We can assume that $b-\delta_1$ and $b+\delta_2$ are non-resonant.
In what follows we fix $\varepsilon > 0, r,\delta_1$ and
$\delta_2$ assuming that they are as small as needed. Put $E =
1/(2\varepsilon^2)$, $r^- = 1/\sqrt{2(b+\delta_2)}$ and $r^+ =
1/\sqrt{2(b-\delta_1)}$.

According to \fullref{thm-corresp}, a contractible positive
loop of contactomorphisms on $\PP_+T^*X$ allows us to construct a
contact squeezing of the domain $\wh U_{r^+}$ into $\wh U_{r^-}$.
Moreover, this squeezing can be done in the class of fiberwise
starshaped domains, and hence it can be performed via a contact
isotopy which fixes pointwise a small neighborhood $\wh
U_{\varepsilon}$ of the zero section. In other words, there exists
a contact isotopy $\Phi_t\co  V \to V, \; t \in [0,1],$ such that
$\Phi_0 = \id$, $\Phi_t \equiv \id$ on $\wh U_{\varepsilon}$ for
all $t$ and
$$\Phi_1 (\mathrm{Closure}(\wh U_{r^+})) \subset \wh U_{r^-}\;.$$

Note that the inclusion homomorphism
$$I \co \CH(\wh U_{\e})\to \CH(\wh U_{r^-})$$ splits as a composition
$$\CH(\wh U_{\e})\to \CH(\Phi_1(\wh U_{r^+}))\to \CH(\wh U_{r^-}),$$
and hence
$$\dim\image\left(\CH(\wh U_{\e})\to \CH(\Phi_1(\wh U_{r^+}))\right)\geq
\dim\image\left(\CH(\wh U_{\e})\to \CH(\wh U_{r^-})\right).$$ But
\begin{equation}\label{eqrankineq}
\image\left(\CH(\wh U_{\e})\to \CH(\Phi_1(\wh
U_{r^+}))\right)\simeq \image\left(\CH(\wh U_{\e})\to \CH(\wh
U_{r^+})\right).
\end{equation}
According to the commutative diagram in \fullref{corloops}
this yields that
$$\beta(E,b-\delta_1) \geq \beta(E,b+\delta_2)\;.$$
Note that $E$ can be taken arbitrarily large, and hence \fullref{lem-dimensions-tw}(ii) yields
$$\beta(b-\delta_1) \geq \beta(b+\delta_2)\;$$
in contradiction with inequality \eqref{eq-jump}. This
contradiction proves the theorem.
\end{proof}

\begin{rem}
{\rm Note that if $\pi_1(X)$ has infinitely many conjugacy classes
it is sufficient to use only $0$--dimensional cohomology. The jump
points of the function $\beta(a)$ have in this case a transparent
geometric meaning: they correspond to $\ell(\alpha)^2/2$ where
$\alpha$ runs over all non-trivial free homotopy classes of loops
on $X$ (that is non-trivial conjugacy classes in $\pi_1(X)$ ) and
$\ell(\alpha)$ is the minimal length of closed geodesic in
$\alpha$. }
\end{rem}

\begin{proof}[Proof of \fullref{lem-dimensions-tw}]
The case when $\pi_1(X)$ has infinitely many conjugacy classes is
obvious already on the level of $0$--dimensional cohomology. To
deal with the simply connected case let us recall some facts from
the theory of free loop spaces:
 
{\bf (I)}\qua The space $H^*(\cL X)$ has infinite dimension.
Let us emphasize that we are dealing here with $\Z_2$
coefficients. This is proved in \cite{FTVP}.

{\bf (II)}\qua $H^*(\cL X)$ is the inverse limit of
$H^*(\cL^E X)$ as $E \to \infty$ (see eg \cite{Weber}).

\medskip First we prove statement (i) of the lemma.
It is more convenient to work here with homology instead of
cohomology. Any homology class $\alpha\in H_*(\cL X)$ can be
represented by a cycle contained in $\cL^aX$ for some finite $a$.
Thus the number
$$A(\alpha) = \inf\{a\;|\; \alpha \in \image(\tau_a^*)\}\;,$$
where $\tau_a^*$ stands for the natural morphism $$ H_*(\cL^a X)
\to H_* (\cL X)\;,$$ is well defined for every $\alpha\in H_*(\cL
X)$. Since $\tau_a^*$ factors through $H_*(\cL^b X)$ for all $b
>a$ we have that $\alpha \in \image(\tau_b^*)$ provided $b> A(\alpha)$.
Take any $N
> 0$. Using (I), choose $N$ linearly independent cycles
$\alpha_1,...,\alpha_N \in H_*(\cL X)$. It follows that $$\beta(a)
= \dim H_*(\cL^aX) \geq N$$ provided $a>\max_i A(\alpha_i)$. Thus
$$\beta(a)\mathop{\to}\limits_{a\to\infty}\infty$$
as required.

Let us prove statement (ii) of the lemma. Put
$$\beta(\infty,a) = \lim_{E \to \infty} \beta(E,a)\;. $$
The limit exists since $\beta(E,a)$ is a non-increasing function
of $E$. Observe that $\beta(E,a) \geq \beta(a)$ for all
sufficiently large $E$, and therefore
\begin{equation}\label{eq-vsp-beta} \beta(\infty,a)\geq\beta(a).
\end{equation}
For every $c>0$, $H^*(\cL^c X)$ is a finite-dimensional space over
$\Z_2$, since $\cL^cX$ has homotopy type of a finite $CW$--complex.
Thus the family of decreasing subspaces $\image(\tau_{E,c})$
stabilizes to a subspace $F_c \subset H^*(\cL^c X)$ when $E \to
\infty$.\footnote{This is called the Mittag-Leffler condition in
homological algebra.} Clearly, $\dim F_a = \beta(\infty,a)$. Note
that $\tau_{p,q}$ sends $F_p$ {\bf onto} $F_q$ for $p > q$.
Consider the following commutative diagram:
\[ \xymatrix{ \varprojlim F_E \ar[r] \ar[d] &
F_a\ar[d] \\
\varprojlim H^*\cL^EX \ar[r] & H^*\cL^aX \ & } \] Here the
vertical arrows are natural inclusions and the horizontal arrows
are natural projections. The upper horizontal arrow is surjective,
while the inverse limit in the south-west corner equals $H^*(\cL
X)$ by (II). Thus $\beta(a) \geq \beta(\infty,a)$, and so in view
of inequality \eqref{eq-vsp-beta} we have $\beta(a) =
\beta(\infty,a)$. Since $\beta(E,a)$ is an integer, we conclude
that the function $\beta(E,a)$ equals $\beta(a)$ for a
sufficiently large $E$. This completes the proof of (ii).
\end{proof}

\section{Discussion and further directions
}\label{sec-discussion-open-problems}

Here we touch miscellaneous topics related to geometry of contact
domains and transformations. We discuss (non)-squeezing of contact
domains from the viewpoint of quantum mechanics. Then we speculate
on links between orderability and fillability of contact
manifolds. Finally, we introduce a canonical semigroup associated
to a contact manifold. This semigroup carries footprints of a
partial order even when the manifold is non-orderable.

\subsection{Contact non-squeezing at the quantum scale} \label{subsec-hbar}

We start with a brief (and standard) discussion of quasi-classical
approximation in quantum mechanics. Consider the Cauchy problem
for the Schr\"{o}dinger equation
\begin{equation}\label{schro}
i\hbar \frac{\p \psi}{\p t} = - \frac{1}{2} \hbar^2 \Delta \psi +
W(q)\psi
\end{equation}
which describes the quantum motion of a mass $1$ particle on
$\R^n$ in the presence of a potential force $F = -\frac{\p W}{\p
q}$. Here $\hbar$ is the Planck constant. Assume that the initial
complex-valued wave function is given in the form
$$\psi_0(q) = \rho_0^{\frac{1}{2}}(q) e^{iF_0(q)/\hbar}\;.$$
Here $\rho_0$ is the probability density of the particle position
in $\R^n$ and $F_0\co  \R^n \to S^1_{\hbar}$ is the phase which takes
values in the circle $S^1_{\hbar} = \R/(2\pi \hbar\Z)$.

Quasi-classical approximation gives the following recipe for an
approximate solution of the Cauchy problem which involves contact
geometry and dynamics: Consider the contact manifold $(V_{\hbar} =
\R^{2n} \times S^1_{\hbar}, \zeta)$ where $\zeta = \mathrm{Ker}(du
-pdq)$. Consider the lift $f_t\co  V_{\hbar} \to V_{\hbar}$ of the
Hamiltonian flow on $\R^{2n}$ describing the classical motion of
our particle. The flow $f_t$, which is given by the system
\[ \left\{ \begin{array}{lll}
\dot{q}=p \\ \dot{p} = -\frac{\p W}{\p q} \\ \dot{u} =
\frac{1}{2}p^2 -W(q) \end{array} \right. \] preserves the contact
structure $\zeta$.

Given a function $K\co  \R^{n} \to S^1_{\hbar}$ and a probability
density $\rho(q)dq$ on $\R^{n}$, consider a Legendrian submanifold
$$L(K) = \{u = K(q,0), p = \frac{\p K(q,0)}{\p q}\} \subset V\;$$
equipped with the density $ \sigma(\rho) := \tau^*(\rho(q)dq)$
where $\tau\co  V \to \R^{n}$ stands for the natural projection. For
$t$ small enough, the image $f_t(L(F_0))$ can be written as
$L(F_t)$ and the density $f_{t*}\sigma(\rho_0)$ can be written as
$\sigma(\rho_t)$. The main outcome of this construction is as
follows. The wave function
$$\psi_t(q) = \rho_t^{\frac{1}{2}}(q) e^{iF_t(q)/\hbar}\;$$
is an approximate solution of the Schr\"{o}dinger equation
\eqref{schro} in the following sense: it satisfies the
Schr\"{o}dinger equation up to higher order terms in $\hbar$.

The above discussion suggests that it is worth studying contact
geometry and dynamics on the manifold $(V_{\hbar},\zeta)$ keeping
track of the value of the Planck constant $\hbar$. In what follows
we reexamine our (non)-squeezing results from this viewpoint. Put
$h = 2\pi\hbar$ and for a domain $U \subset \R^{2n}$ define its
prequantization as $\wh{U}_{\hbar} = U \times S^1_{\hbar}$. Next
note that the map
$$\Phi\co  V \to V_{\hbar},\;\; (p,q,t) \to
(\sqrt{h}p,\sqrt{h}q, h(t + \frac{1}{2}p\cdot q))$$ establishes a
contactomorphism between $(V,\xi)$ and $(V_{\hbar},\zeta)$.
Furthermore, given a domain $U \subset \R^{2n}$,
$$\Phi(\wh{U}) = (\wh{\sqrt{h}\cdot U})_{\hbar}\;.$$
Therefore, after a rescaling, our squeezing and non-squeezing
results presented in Sections \ref{subsec-sq} and
\ref{subsec-capac} extend to the contact manifold
$(V_{\hbar},\zeta)$. In particular, prequantization of balls
$B^{2n}(R)$ of sub-quantum size $R < h$ in the classical phase
space gives rise to squeezable domains $\wh{B}_{\hbar}^{2n}(R)$ in
$V_{\hbar}$. One can say that {\it on the sub-quantum scale,
symplectic rigidity is lost after prequantization.} On the other
hand, rigidity of balls in the classical phase space persists
after prequantization provided the size of the ball is an integer
multiple of $h$: $\wh{B}^{2n}(mh)$ cannot be squeezed into itself.
It would be interesting to explore whether this geometric
phenomenon has some physical meaning.

Note that the balls $B^{2n}(R)$ are energy levels of the
classical harmonic oscillator, while the sequence $\{mh\}$ is, in
accordance to Planck's hypothesis, the energy spectrum of the
quantized harmonic oscillator. (Planck's hypothesis, which
appeared at the early stage of development of quantum mechanics,
turned out to be non-precise and was later on corrected in a more
advanced model. Still, it gives a satisfactory approximation to
the energy spectrum.) It is unclear whether appearance of the
sequence $\{mh\}$ in the context of (non)-squeezing is occasional
or not. A possibility, which is still open, is that
$\wh{B}^{2n}_{\hbar}(R)$ cannot be squeezed into itself {\it for
any} $R
> h$.

In \cite{deGosson1,deGosson2} de Gosson revised the classical
notion of quantum cell in the phase space $\R^{2n}$ in terms of
symplectic capacities. The following feature of de Gosson's
definition is of interest in our context: for a quantum cell $U$,
both $\uc (U)$ and $\oc (U)$ are of order $\hbar$. Interestingly
enough, the same scale appears in our (non)-squeezing problem. To
be more precise, let us introduce the following notion which is
natural when one studies transition from squeezing to
non-squeezing in families of domains. Fix a constant $D>0$ and
consider the class $\Xi_D$ of domains $U \subset \R^{2n}$ with
$\oc (U) \leq D \cdot \uc (U)$. For instance, all convex domains
lie in $\Xi_D$ for $D = 4n^2$ (see \cite{Viterbo1},\cite{CHLS}).
We say that a domain $U \in \Xi_D$ is {\it critical} if
$\wh{U}_{\hbar}$ is negligible (that is $\wh{U}_{\hbar}$ can be
contactly squeezed into $\wh{B}^{2n}_{\hbar}(r)$ for all $r >0$),
while $\wh{(2\cdot U)}_{\hbar}$ is not negligible. \fullref{thm-capacities} yields that $\uc(U) \sim \hbar$ for every
critical domain $U$.

\subsection{Fillability and orderability}\label{subsec-fillability-orderability}

\begin{conjecture}\label{conjcritical}
If a contact manifold $(P^{2n-1},\eta)$ admits a Weinstein filling
$W^{2n}$ with $H_n(W)\neq 0$ then it is orderable.
\end{conjecture}

This conjecture is supported by spaces of co-oriented contact
manifolds, see \fullref{thmcont-elements} above. If this
conjecture were true then together with \fullref{thmsubcrit}
it would leave the orderability problem for Weinstein fillable
contact manifolds open essentially only in the $1$--subcritical
case, eg, for $S^2\times S^1$ with its tight contact structure.

Note that the real projective space $P:= \R P^{2n-1}$ with the
standard contact structure, which is orderable by Givental's
theorem (see \cite{EP}), is not Weinstein fillable for any $n>2$
(though it is symplectically fillable). Indeed, suppose that $W$
is such a filling. Then
$$H^2(W,P;\Z_2)=H_{2n-2}(W;\Z_2)=0$$ for $2n-2>n$, or $n>2$.
Looking at the cohomology long exact sequence
$$\xymatrix{H^1(W,\Z_2) \ar[r]^r & H^1(P,\Z_2) \ar[r] & H^2(W,P;\Z_2)=0}$$
we see that the map $r$ is onto and hence there is a class $a\in
H^1(W;\Z_2)$ such that $b=r(a)$ is the generator of $H^1(P;\Z_2)$.
Hence, we have $b^{2n-1}\neq 0$ which contradicts to the fact that
$a^{2n-1}\in H^{2n-1}(W;\Z_2)=0$. It would be interesting to study
orderability of more general contact manifolds which are not
Weinstein fillable. It seems likely that the methods of the
current paper extend to a more general case of symplectically
fillable contact manifolds.

The case $n=2$ corresponds to the real projective $3$--space $\R
P^3$ which is contactomorphic to the space of co-oriented contact
elements of $S^2$. It is orderable in view of \fullref{thmcont-elements}.

\subsection{Contact diffeomorphisms and semigroups} \label{subsec-semigroup}
\subsubsection{Semigroups and orderability}\label{secsemi-order}
Let $(P,\eta)$ be a closed contact manifold. Denote by $G$ the
group $\mathrm{Cont}_0 (P,\eta)$ and by $\wt G$ its universal cover.
We will often identify (paths of) contactomorphisms of $P$ with
the corresponding (paths of) $\R_+$--equivariant symplectomorphisms
of the symplectization $SP$.

Consider the set of positive paths
$$f =\{f_t\},\; t \in [0,1]\;, f_0 = \id,$$
on $G$ (that is the paths generated by positive contact
Hamiltonians). We say that two such paths $f'$ and $f''$ are
equivalent if $f'_1 = f''_1$ and the paths are homotopic through
positive paths with fixed end points. Denote by
$\ocG=\ocG(P,\eta)$ the set of equivalence classes. \footnote{When
necessary, we can assume without loss of generality that our paths
are generated by 1--periodic in time contact Hamiltonians, see
Lemma 3.1.A in \cite{EP}.} The operation
$$[\{f_t\}]\circ [\{g_t\}] = [\{f_t g_t\}]\;,$$
where $\{f_t\}$ and $\{g_t\}$ are positive paths, defines a
semigroup structure on $\ocG$. Note that $\ocG$ comes along with
the natural projections $\tau\co  \ocG \to \wt{G}$ and $\pi\co  \ocG \to
G$ which are morphisms of semigroups.

The (non)-orderability of $(P,\eta)$ can be recognized in terms of
$\ocG$ and $\tau$ as follows. Denote by $\wt{\id} \in \wt{G}$ the
canonical lift of the identity. When $(P,\eta)$ is orderable,
there are no positive contractible loops on $G$ and hence
$\wt{\id}$ does not belong to the image of $\tau$. In contrast to
this, when $(P,\eta)$ is not orderable the morphism $\tau$ is
surjective. Indeed, in this case we have a positive contractible
loop $\{\phi_t\}$ of contactomorphisms with $\phi_0=\phi_1 = \id$.
Thus for every path of contactomorphisms $f= \{f_t\}$ on $G$
starting at the identity, the path $h=\{f_t\phi_t^N\}$ is positive
for $N$ large enough, and hence we get an element $[h] \in \ocG$
with $\tau([h]) = [f] \in\wt{G}$.

\subsubsection{Partial order for non-orderable manifolds}
Assume now that $(P,\eta)$ is non-orderable. In this case
$$\ocG_0:= \tau^{-1}(\wt{\id}) \subset \ocG$$
is a subsemigroup. Its elements represent connected components of
the space of positive contractible loops on $G$ based at $\id$.

\begin{prop}\label{prop-center} {\rm The center of
$\ocG(P,\eta)$ contains $\ocG_0(P,\eta)$.}
\end{prop}

The proposition immediately follows from the formula:
\begin{equation} \label{eq-permutation}
[\{f_tg_t\}] = [\{g_tg_1^{-1}f_tg_1\}]\;,
\end{equation}
where $\{f_t\}$ and $\{g_t\}$, $t \in[0,1]$, are arbitrary
positive paths with $f_0=g_0 = \id$. To prove this formula, define
new paths
\[
f'_t =\begin{cases} \id &\hbox{for}\;\; t \in [0,\frac{1}{2}];\cr
f_{2t-1}& \hbox{for}\;\; t \in [\frac{1}{2},1]\cr
\end{cases}\;\;\;\text{and}\;\;\;
g'_t =\begin{cases} g_{2t} &\hbox{for}\;\; t \in
[0;\frac{1}{2}];\cr g_1 & \hbox{for}\;\; t \in [\frac{1}{2},1]\cr
\end{cases}
\;.
\]
Note now that $$f'_tg'_t = g'_tg_1^{-1}f'_tg_1$$ for all $t$.
Furthermore, $\{f_t\}$ and $\{f'_t\}$ (respectively, $\{g_t\}$ and
$\{g'_t\}$) are homotopic with fixed end points such that all the
intermediate paths except $\{f'_t\}$ (respectively, $\{g'_t\}$)
are smooth and positive. Observe that $\{f'_t\}$ and $\{g'_t\}$
are piece-wise smooth and non-negative. By an appropriate
``smoothing of corners" argument, we correct this and get the
desired homotopy between $\{f_tg_t\}$ and $\{g_tg_1^{-1}f_tg_1\}$
through smooth positive paths, thus proving formula
\eqref{eq-permutation}.

An interesting feature of the semigroup $\ocG_0$ is given by the
following proposition.

\begin{prop}\label{prop-aperiodic}
For every $\gamma \in \ocG$ and $\theta \in \ocG_0$ there exists
$N \in \N$ such that $\theta\gamma^N = \gamma^N$.
\end{prop}

\begin{proof} Assume that $\gamma$ is represented by a
path $\{h_t\}$ generated by a contact Hamiltonian $H(z,t)>0$ and
$\theta$ is represented by a loop $f= \{f_t\}$ generated by the
contact Hamiltonian $F(z,t) > 0$. Let $f^{(s)} = \{f_{s,t}\}$, $s
\in [0,1]$, be a homotopy of $f = f^{(0)}$ to the constant loop,
and let $F_s(z,t)$ be the corresponding contact Hamiltonians.
Choose $N$ large enough so that
$$F_s(f_{t,s}y,t) + NH(y,Nt) >0$$
for all $s,t \in [0,1]$ and $y \in SP$. Consider the homotopy of
loops $\{f_{t,s}\circ h_t^N\}$ which joins $\{f_th_t^N\}$ and
$\{h_t^N\}$. The intermediate loops are generated by the
Hamiltonians
$$F_s(z,t) + NH(f_{s,t}^{-1}z,Nt)$$
which are positive due to our choice of $N$. Thus
$\theta\gamma^{N}=\gamma^N$ as required.
\end{proof}

\begin{cor}\label{cor-epigroup} For every $\theta \in
\ocG_0$ there exists $N \in \N$ such that $\theta^{N+1}=\theta^N$.
\end{cor}

The corollary implies that for each $\theta\in\ocG_0$
the sequence $\{\theta^N\}$, $N \to +\infty$, stabilizes to an
element which we denote by $\overline\theta$.

\begin{cor}\label{corunique-stable}
For any $\g,\theta\in\ocG_0$ we have
$\overline\g=\overline\theta$.
\end{cor}

Indeed, according to Propositions \ref{prop-aperiodic},
\ref{prop-center} and \fullref{cor-epigroup} we have
$$\overline\gamma=\overline\theta\overline\gamma=\overline\gamma\overline\theta=\overline\theta\,.$$

\begin{cor}\label{cor-nil}
There exists a unique ``stable" element in $\ocG_0$, denoted by
$\theta_{\st}$, which is equal to $\overline\theta$ for any
$\theta\in\ocG_0$. The element $\theta_\st$ can be characterized
by any of the following properties:
\begin{itemize}
\item $\theta_{\st}$ is the unique {\it idempotent},
$\theta_{\st}^2=\theta_{\st}$, in the semigroup $\ocG_0$;
\item $\theta_{st}$ is the {\it zero element} of $\ocG_0$, that is
for any $\gamma\in\ocG_0$ we have
$$\gamma\theta_{\st}=\theta_{st}\gamma =\theta_{st}.$$
\end{itemize}
\end{cor}

Thus we established that some power of every element in
$\ocG_0$ is equal to the zero element. Such semigroups are called
{\it nilsemigroups}.

\begin{example} \label{ex-sphere-ocG} {\rm
Consider the semigroup $\ocG_0(S^{2n-1})$ associated with the
standard contact sphere, where $n \geq 2$. The stable element
$\theta_{\st}$ can be explicitly identified in this case. Let
$\theta \in \ocG_0$ be the element represented by the positive
contractible loop $\varphi= \{\varphi_t\}$ constructed in \fullref{thm-main-consttw}. We claim that $\theta$ is an idempotent:
$\theta^2 = \theta$, and hence $\theta$ coincides with
$\theta_{\st}$. To prove this recall that $\{\varphi_t\}$ is
generated by the Hamiltonian $\Phi(z,t)$ with $\Phi(z,t) > \pi
|z|^2$ (see equation \eqref{eq-positiveloop-bound} above). On the
other hand, we have seen in the proof of \fullref{mu}(ii) that
there exists a homotopy $\varphi^{(s)} = \{\varphi_{t,s}\}$ of
$\varphi=\varphi^{(0)}$ to the constant loop $\varphi^{(1)} \equiv
\id$ such that the Hamiltonians $\Phi_s(z,t)$ generating the loops
$\varphi^{(s)}$ satisfy $\Phi_s(z,t) \geq -\pi |z|^2$. Therefore,
\begin{equation}\label{eq-idempotent}
\Phi(z,t_1)+\Phi_s(z,t_2) > 0
\end{equation}
for all $z \in \C^n \setminus \{0\}$, $t_1,t_2 \in S^1$ and $s \in
[0,1]$.

Let $\theta^* \in \ocG_0$ be the element represented by the loop
$\varphi^* = \{\varphi^{-1}_{1-t}\}$. Consider the homotopy of
$\varphi^*\varphi$ to $\varphi$ given by
$\{\varphi^{-1}_{1-t,s}\varphi_t\}$, where $s \in [0,1]$. The
Hamiltonians of the intermediate loops are given by
$$\Phi_s(\varphi_{1-t,s}z,1-t) + \Phi(\varphi_{1-t,s}z,t)\;.$$
These Hamiltonians are positive in view of \eqref{eq-idempotent},
and therefore $\theta^*\theta = \theta$.

Consider now the homotopy of $\varphi^*\varphi$ to $\varphi^*$
given by $\{\varphi^{-1}_{1-t}\varphi_{t,s}\}$, where $s \in
[0,1]$. The Hamiltonians of the intermediate loops are given by
$$\Phi(\varphi_{1-t}z,1-t) + \Phi_s(\varphi_{1-t}z,t)\;.$$
Again, these Hamiltonians are positive in view of
\eqref{eq-idempotent}, and therefore $\theta^*\theta = \theta^*$.
We conclude that $\theta^*=\theta$ and hence
$$\theta^2 = \theta^*\theta = \theta\;,$$
and the claim is proved. }
\end{example}

\begin{rem}\label{remunique-stable-element}
{\rm No complete description of the semigroup $\ocG_0$ is
available even in the simplest case when $P$ is the
three-dimensional sphere $S^3$ equipped with the standard contact
structure. In fact, we do not know a single example when $\ocG_0$
consists of more than just the zero element $\theta_\st$.}
\end{rem}

The next proposition describes in more detail the structure of the
semigroup $\ocG $. Set
$$\Gamma = \{\gamma \in \ocG\;:\;\gamma \theta_{\st} = \g\}\;.$$

\begin{prop}\label{propstructure-semigroup}
\begin{itemize}
\item[(i)] $\Gamma$ is a subgroup of $\ocG$ with the neutral element
$\theta_{\st}$;
\item[(ii)] the restriction of the natural projection
$\tau\co  \ocG \to \wt{G}$ to $\Gamma$ is a group isomorphism;
\item[(iii)] for every $x \in \ocG$ there exists $N \in \N$ such
that $x^N \in \Gamma$ (and hence $\ocG$ is an {\rm epigroup} in
the language of semigroup theory);
\item[(iv)] $\Gamma$ is {\rm a two-sided ideal} in $\ocG$, that is
for every $x \in \ocG$ and $\g \in \Gamma$ we have $x\g,\g x \in
\Gamma$.
\end{itemize}
\end{prop}

\begin{proof} Clearly, $\Gamma$ is closed under the
semigroup multiplication. By the definition of $\Gamma$, the
element $\theta_{\st}$ is its neutral element. Recall that the
morphism $\tau\co  \ocG \to \wt{G}$ is onto. Thus, given any $x \in
\Gamma$, we can choose $y \in \ocG$ with $\tau(y) = \tau(x)^{-1}$.
This yields $\tau(xy) = \tau(yx)= \wt{\id}$, and hence
$xy,yx\in\ocG_0$. Put $y'=y\theta_{\st}$. Note that $y'\in \Gamma$
since $\theta_{st}^2 =\theta_{st}$. Recalling that $\theta_{st}$
is the zero element of $\ocG_0$, we get $$xy' =
xy\theta_{\st}=\theta_{\st}, \;\hbox{and}\; y'x = y\theta_{st}x=
yx\theta_{st}=\theta_{st}\,.$$ Thus $y'$ is the inverse of $x$ in
$\Gamma$, and this proves Property (i). Taking into account that
$x\theta_{st} \in \Gamma$ and $\tau(x\theta_{st}) = \tau(x)$ for
every $x \in \ocG$, we get that $\tau|_{\Gamma}$ is onto since
$\tau$ is onto. Let us show that the kernel of the homorphism
$\tau|_\Gamma$ is trivial. Indeed, take any
$\theta\in\left(\tau|_\Gamma\right)^{-1}(\wt\id)=\Gamma\cap\ocG_0$.
Then we have $\theta\theta_{\st}=\theta$ by the definition of
$\Gamma$ and $\theta\theta_\st=\theta_\st$ since $\theta_{\st}$ is
the zero element of $\ocG_0$. Hence, $\theta=\theta_\st$, and we
get Property (ii). Property (iii) follows from \fullref{prop-aperiodic}. Property (iv) is a straightforward
consequence of the definitions.
\end{proof}

Let us introduce a new relation on the semigroup $\ocG$. Given
$x,y\in\ocG$ we say that $x \psucceq y$ whenever either $x=y$ or
$x = zy$ for some $z \in \ocG$.

\begin{proposition}\label{propsemigroup-order}
The relation $\psucceq$ is
\begin{itemize}
\item[(i)]
reflexive and transitive;
\item[(ii)] trivial on $\Gamma$, ie, $x\psucceq y$ and $y\psucceq x$ for
any $x,y\in\Gamma$;
\item[(iii)] satisfies the condition $x\psucceq y$ for any $x\in\Gamma$ and
$y\in\ocG\setminus\Gamma$;
\item[(iv)] defines a genuine partial order on $\ocG \setminus \Gamma$.
\end{itemize}
\end{proposition}
\begin{proof}
Property (i) is obvious. Property (ii) follows from the fact that
$\Gamma$ is a subgroup of $\ocG$. Take any $x \in \Gamma$ and
$y\in\ocG\setminus\Gamma$. Recall that $y\theta_{st} \in \Gamma$.
Using Properties (i) and (ii) we get $$x \psucceq y\theta_{st}
\psucceq y\;,$$ which proves Property (iii). Let us check now that
the relation $\psucceq$ is anti-symmetric on
$\ocG\setminus\Gamma$. Suppose that $x \psucceq y$ and $y \psucceq
x$ for some $x,y \in \ocG \setminus \Gamma$. Assume on the
contrary that $x \neq y$. Then $x = ay$ and $y = bx$ for some $a,b
\in \ocG$. Thus $x = abx$, and so $\tau(ab) = \wt{\id}$, which
implies that $\theta:=ab \in\ocG_0$. Multiplying the equation $x
=\theta x$ by successive powers of $\theta$ we have
$$x =\theta x=... =\theta^Nx=\theta_{\st}x$$
provided $N$ is large enough. We get a contradiction with our
assumption $x \notin \Gamma$. Thus the remaining possibility is
$x=y$, and we proved Property (iv).
\end{proof}

Since $\Gamma$ is a two-sided ideal of $\ocG$, we can consider
{\it the Rees quotient semigroup} $\ocG /\Gamma $. By definition,
this is a new semigroup on the set $\cR=(\ocG \setminus \Gamma)
\cup \{0\}$ in which the product of elements $x,y \in {\cR}$
concides with their product $xy$ in $\ocG$ if $x,y,xy \in \ocG
\setminus \Gamma$, and equals $0$ otherwise. Property
\ref{propstructure-semigroup}(iii) yields that for every element
$r$ of $\cal{R}$ there exists $N \in \N$ such that $r^N = 0$.
\fullref{propsemigroup-order} implies that the relation
$\psucceq$ on $\ocG \setminus \Gamma$ defines a genuine partial
order on $\cal{R}$, denoted in the same way, where we declare that
$0 \psucceq r$ for every $r \in \cal{R}$.
\begin{prop}\label{prop-biinv}
The partial order $\psucceq$ on $\cal{R}$ is bi-invariant: if $x
\psucceq y$ for some $x,y \in \cal{R}$ then $xc \psucceq yc$ and
$cx \psucceq cy$ for every $c \in \cal{R}$.
\end{prop}

Before proving the proposition let us make the following
general remark. The group $G$ acts on $\ocG$ by conjugations, that
is the action of an element $g \in G$ is given by
$$A_g ([\{f_t\}]):= [\{g^{-1}f_tg\}]\;.$$
This action preserves point-wise the subsemigroup $\ocG_0$ and
leaves invariant the subgroup $\Gamma$. Formula
\eqref{eq-permutation} above yields
\begin{equation}\label{eq-permut-on}
xy = y\cdot A_{\pi(y)}x = A_{\pi(x)}^{-1}y \cdot x
\end{equation}
for all $x,y \in \ocG$.

\begin{proof}[Proof of \fullref{prop-biinv}] Assume that
$x,y \in \cal{R}$ with $x \psucceq y$. Suppose that $x \neq 0$
(the case $x=0$ is trivial). By definition, $x = zy$ for some $z
\in \ocG$. Given any $c \in \cR$ we have $xc = zyc$ and hence $xc
\psucceq yc$. Further, $x = yw$ with $w = A_{\pi(y)}z$ in view of
formula \eqref{eq-permut-on}. Applying this formula again, we get
$$cx = cyw = A_{\pi(cy)}^{-1}w \cdot cy\;$$
and hence $cx \psucceq cy$. \end{proof}

\subsubsection{Semigroup $\cR$ and contact squeezing}

A priori, the set $\ocG \setminus \Gamma$ might be empty even for
non-orderable contact manifolds (though we are not aware of such
an example). Interestingly enough, in order to prove that it is
non-empty for the standard contact sphere (see \fullref{prop-semigr-nonempty} below), we use the contact
non-squeezing theorem. This suggests that non-triviality of the
partially ordered semigroup $\cR$ is a ``hard" aspect of the
``soft" non-orderability phenomenon.

Let us begin with a general setting. Recall that whenever
$(P,\eta)$ is an ideal contact boundary of a Liouville manifold
$(M,\omega,L)$, \fullref{lmpath-domain} establishes a
correspondence between elements $\g \in \ocG(P,\eta)$ and
fiberwise starshaped domains $U(\g) \subset M \times S^1$ up to
contact isotopy. Generalizing the definition given in \fullref{subsec-capac}, we say that a subset of $M \times S^1$ is {\it
negligible} if it can be contactly squeezed into an arbitrarily
small neighborhood of $\Core(M) \times S^1$.

\begin{prop}\label{lmGamma-squeezing}
For any $\gamma \in\Gamma$ the domain $U(\gamma)$ is negligible.
\end{prop}

\begin{proof} Take any positive time-independent contact Hamiltonian
$E$ on $SP$. The claim can be equivalently expressed by saying
that for any $C>0$ the element $\gamma$ can be generated by a
1--periodic Hamiltonian which is $\geq C\cdot E(z)$. But by the
definition of $\Gamma$ we have $\gamma\theta_\st=\gamma$. The
claim follows because $\theta_\st^N = \theta_{st}$ for all $N \in
\N$ and hence $\gamma$ can be generated by a contact Hamiltonian
which is $\geq C \cdot E(z)$ with an arbitrarily large $C$.
\end{proof}

Let us apply this result to the case when $(P,\eta)$ is
the standard contact sphere $S^{2n-1}$. The following proposition
shows that, despite the fact that the sphere is not orderable for
$n\geq 2$, the semigroup
$$\cR(S^{2n-1})=\left(\ocG \setminus\Gamma\right)\cup\{0\}$$
contains a subset of cardinality continuum. Let
$\epsilon_c\in\ocG$, $c>0$, be the element represented by the path
$$z \to e^{2\pi ic t }z$$
with $t \in [0,1]$.

\begin{prop}
\label{prop-semigr-nonempty} $\epsilon_c\in \ocG \setminus \Gamma$
for each $c\in(0,1)$.
\end{prop}

\begin{proof}[Proof of \fullref{prop-semigr-nonempty}] 
It follows from
\fullref{thm-nonsq} that if $c\in(0,1)$ then the domain
$U(\epsilon_c) = \wh{B}^{2n}(\frac1c)$ cannot be squeezed into
$\wh B^{2n}(1)$. Then \fullref{lmGamma-squeezing} yields
$\epsilon_c \in \ocG \setminus \Gamma$. This completes the proof.
\end{proof}

Furthermore, the partial order $\psucceq$ on $\cR $ is
non-trivial since $\epsilon_a \psucceq \epsilon_b$ for $1>a>b>0$.

\subsection{Quantum product in contact homology and other useful
tools}\label{secquantum} Let $(P,\eta)$ be an ideal contact
boundary of a Liouville manifold $(M,\omega,L)$. Under suitable
assumptions on the first Chern class of $(M,\omega)$ (see \fullref{secflavors} above) we have a well defined contact homology
theory for fiberwise starshaped domains in $M \times S^1$. Denote
by $\CH(\gamma)$ the contact homology of the domain $U(\gamma)$.
It would be interesting to relate the correspondence $\gamma \to
\CH(\gamma)$ with the multiplication in $\ocG$. For instance, it
sounds likely that there exists a natural pairing
$$\CH(\gamma_1)\otimes\CH(\gamma_2) \to \CH(\gamma_1\gamma_2)$$
whose construction imitates the pair of pants product in Floer
homology. This pairing is a potentially useful tool for studying
the semigroup $\ocG(P,\eta)$.

It would be interesting to elaborate this idea in the case when
$(P,\eta)$ is a prequantization space of a closed symplectic
manifold $(N,\sigma)$. Write $p\co  P \to N$ for the natural
projection, denote by $e_t$ the natural circle action on $P$, and
write $\beta$ for the connection 1--form on $P$ whose kernel equals
$\eta$.

For instance, think of the Hopf fibration $ S^{2n-1} \to \C
P^{n-1}$ as of the prequantization of the complex projective space
equipped with a suitable multiple of the Fubini-Study form. Here
the sphere is considered as the ideal contact boundary of $\C^n$,
and $e_t(z) = e^{2\pi i t}z$.

Returning to the general setting, take any path of Hamiltonian
diffeomorphisms $\{h_t\}$ on $N$ representing an element
$\tilde{h}$ in the universal cover $\wt{\mathrm{Ham}(N,\sigma)}$ of
the group of Hamiltonian diffeomorphisms of $(N,\sigma)$. We can
lift it to a path of contactomorphisms of $(P,\eta)$ in the
following way. Let us normalize the Hamiltonian $H_t$ on $N$
generating $\{h_t\}$ by requiring that $H_t$ has zero mean with
respect to the symplectic volume on $N$. Look at the projection
$SP \to P \to N$. Take the (time-dependent) contact Hamiltonian on
$SP$ whose restriction to $\mathrm{graph}(\beta) \subset SP$
coincides with the lift of $H_t$. This Hamiltonian generates a
flow, say $f_t$, which is a lift of the original Hamiltonian flow
$h_t$. Note now that the path $\{e_{ct} f_t\}$, $t \in [0,1]$, is
positive for sufficiently large $c$ and hence we get a well
defined element $\gamma_c \in \ocG$. Thus we can associate to
$\tilde{h} \in \wt{\mathrm{Ham}(N,\sigma)}$ a 1--parametric family of
vector spaces $\CH(\gamma_c)$, where $c$ is sufficiently large. It
would be interesting to establish a link between this family of
vector spaces and Hamiltonian Floer homology of $\tilde{h}$.

It seems also useful to develop an equivariant version of the
contact homology theory. For instance, $\Z/2$--equivariant contact
homology for domains in $\R^{2n}\times S^1$ may help recovering
Givental's result about the orderability of the standard contact
structure on $\R P^{2n-1}$.

\section[Appendix A: Proof of \ref{lem-somewhereinjective}]{Appendix A: Proof of \fullref{lem-somewhereinjective} }

The proof is divided into several steps.

{\bf Step 1}\qua Note that $F$ is a proper map. Denote by $\Sigma$ the image of
$F$. A point $p \in \Sigma$ is regular if $\Sigma$ is a smooth
submanifold of $W$ in a neighborhood of $p$, otherwise $p$ is
called singular. The Micallef-White theorem \cite{MiW} states that
locally a non-parametrized $J$--holomorphic curve is
$C^1$--diffeomorphic to a usual holomorphic curve. Thus singular
points form a discrete subset $\mathrm{Sing}(\Sigma)$ in
$\overline{V}$, and near a singular point $p$ the set
$\Sigma\setminus\{p\}$ is a union of a finite number of punctured
discs. This enables us to perform normalization of $\Sigma$ (eg
see Chirka \cite{Chirka}). We get a Riemann surface
$\tilde{\Sigma}$ together with a proper holomorphic map $\pi\co 
\tilde{\Sigma} \to \Sigma$ such that $\pi$ is a diffeomorphism
outside a discrete subset $S \subset \tilde{\Sigma}$ and $\pi(S)$
is exactly the set of singular points of $\Sigma$. Furthermore
there exists a proper lift $\tilde{F}\co  \Upsilon \to
\tilde{\Sigma}$ such that $F = \pi \circ \tilde{F}$.

We claim that the map $\tilde{F}$ is a diffeomorphism. To see
this, denote by $k$ and $g$ the number of ends and the genus of
$\tilde{\Sigma}$, respectively. Since $\gamma_+ \neq \gamma_-$ we
have $k \geq 2$. Let $d$ be the degree of $\tilde F$ and $r$ the
total number of singular points counted with the multiplicites. By
the Riemann-Hurwitz formula, where we use that $\Upsilon$ is a
cylinder, we have
$$d(2g-2+k) + r = 0\;.$$
Since $d\geq1, g\geq 0,k\geq 2$ and $r \geq 0$ we obtain that
$g=r=0$ and $k=2$. In particular, $\tilde{\Sigma}$ is a cylinder
and, since $r=0$, the map $\tilde{F}$ is a non-ramified covering.
Taking into account that the orbits $\gamma_+$ and $\gamma_-$ are
simple, we conclude that $d=1$, and hence $\tilde F$ is a
diffeomorphism. The claim is proved.

In view of this discussion we can assume that $F$ itself is a
diffeomorphism of $\Upsilon \setminus \mathrm{Sing}(F)$ onto $\Sigma
\setminus \mathrm{Sing}(\Sigma)$. Note that the Cauchy-Riemann
equation \eqref{eqdbar} implies that the vertical component
$\varphi$ of a $J$--holomorphic map is related to its $V$--component
$f$ by the equation
\begin{equation} \label{eq-vsp-varphi}
d\varphi = f^*\lambda \circ i
\end{equation}
where $i$ is the complex structure on $\Upsilon$ and $\lambda$ is
the 1--form defining the framing of $V$.
\medskip

{\bf Step 2}\qua Let us denote by $\mathrm{Crit}(f) \subset \Upsilon$ the set of
critical points of $df$. In view of the Cauchy-Riemann equation,
the image $F(\mathrm{Crit}(f))$ corresponds to the tangency points
of $\Sigma$ with the $J$--holomorphic characteristic foliation of
$\overline{V}$. It follows from the generalized similarity
principle due to Hofer and Zehnder \cite{HZ} that $\mathrm{Crit}(f)$
is a discrete subset of $\Upsilon$. This yields that
\begin{equation}\label{eq-vsp-regtw}
\mathrm{Closure} (f(\mathrm{Crit}(f))) \subset f(\mathrm{Crit}(f)) \cup
\g_- \cup \g_+ \;.
\end{equation}

{\bf Step 3}\qua Outside the set $\mathrm{Crit}(f)$ the map $f$ is non-tangent to
the characteristic foliation of $V$, and in particular to the
periodic orbits at infinity $\gamma_- \cup \g_+$. Thus the
(clearly, open) set $\Upsilon' = \Upsilon \setminus f^{-1}(\g_-
\cup \g_+)$ is dense in $\Upsilon$.

\medskip
{\bf Step 4}\qua We claim that the injectivity points lying in $\Upsilon'$ form
an open subset. Indeed, let $z \in \Upsilon'$ be an injectivity
point. Assume on the contrary that there exists a sequence
$\{z_j\} \to z$ such that $z_j$ is not an injectivity point. Note
that $d_{z_j}f \neq 0$ for large $j$ since $d_zf \neq 0$. Thus,
after passing to a subsequence, we can assume that there exists
another sequence $\{w_j\}$ such that $w_j \neq z_j$ for all $j$
and $f(w_j) = f(z_j)$. Since $f$ is injective in a small
neighborhood of $z$ the sequence $\{w_j\}$ does not have $z$ as a
limit point. If it has another limit point, say $w$, we have
$f(z)=f(w)$ which contradicts to the fact that $z$ is an
injectivity point. Thus $\{w_j\}$ contains an unbounded
subsequence. Since $f(z) = \lim f(w_j)$ we conclude that $f(z) \in
\g_- \cup \g_+$ in contradiction with the definition of the set
$\Upsilon'$.
\medskip

{\bf Step 5}\qua Put $$\Upsilon_0 =\Upsilon \setminus f^{-1}\Big{(}
f(\mathrm{Crit}(f)) \cup \g_- \cup \g_+ \Big{)} = \Upsilon'
\setminus f^{-1}\Big{(} f(\mathrm{Crit}(f))\Big{)}\;.$$ This set is
open in view of inclusion \eqref{eq-vsp-regtw}. Let us check that
it is dense in $\Upsilon'$ (and hence dense in $\Upsilon$).
Indeed, assume on the contrary that this is not the case. Then
there exists an open subset, say $P \subset \Upsilon'$, with $f(P)
\subset f(\mathrm{Crit}(f))$. But $f(\mathrm{Crit}(f))$ is a countable
subset of $V$ (see Step 2) so $df$ vanishes identically on $P$.
Equation \eqref{eq-vsp-varphi} yields that $d\phi = 0$ on $P$ as
well. Thus $dF=0$ on $P$, and by the unique continuation argument
the map $F$ is constant, a contradiction.
\medskip

{\bf Step 6}\qua It suffices to show that injectivity points are dense in
$\Upsilon_0$. Take any point $z \in \Upsilon_0$. There exist
neighborhoods $D$ of $z$ and $U$ of $f(z)$ with the following
properties.
\begin{itemize}
\item $f(D) \subset U$;
\item $U \cap (\gamma_- \cup \gamma_+) = \emptyset$;
\item $f|_D$ is an embedding.
\end{itemize}

Let as analyze injective points in $D$. Notice first that in view
of the boundary conditions for $f$ there exists a sufficiently
large closed annulus $A \subset \Upsilon$ containing $D$ and such
that $U \cap f(\Upsilon \setminus A) =\emptyset$. Furthermore,
every point in $f^{-1}(f(z))$ is non-critical for $df$ due to our
definition of $\Upsilon_0$. Thus, enlarging, if necessary, $A$ and
shrinking $D$ and $U$, we can achieve that $f^{-1}(U) \cap A$ is a
union of a finite number of pairwise disjoint open discs
$B_1,...,B_N$ such that $f|_{B_j}$ is an embedding. One of these
discs, say $B_1$, contains $D$. Put $B = B_j$ with $j \geq 2$. We
claim that $D \setminus f^{-1}(f(D) \cap f(B))$ is dense in $D$.
Indeed, otherwise there exist open subsets $D' \subset D$ and $B'
\subset B$ such that $f(D')= f(B')$, and hence there is a a {\it
holomorphic} diffeomorphism $\psi \co  D' \to B'$ such that $f|_{B'}
\circ \psi = f|_{D'}$. Equation \eqref{eq-vsp-varphi} yields
$\varphi|_{B'} \circ \psi = \varphi|_{D'}+c$ where $c \in \R$ is a
constant. If $c = 0$ we have $F|_{D'} = F|_{B'}$ which contradicts
Step 1 of the proof. Thus $c \neq 0$. Recall that the adjusted
almost complex structure $J$ on $\overline{V} = V \times \R$ is
invariant under the shift $T_c\co  (y,s) \to (y,s+c)$ along the
$\R$--direction. Thus our conclusion above can be reformulated as
follows: the cylinder $\Sigma=F(\Upsilon)$ and $T_c(\Sigma)$
intersect over an open subset. By the unique continuation
argument, we have that $\Sigma = T_c(\Sigma)$. Thus $\Sigma =
T_{Nc}(\Sigma)$ for all $N \in \Z$. Taking $N \to \pm\infty$ and
looking at the boundary conditions, we get that $\g_- =\g_+$,
which contradicts our assumption that the orbits at infinity are
distinct.

Thus we proved that every point $z_0 \in \Upsilon_0$ has a
neighborhood where injective points are dense. This completes the
proof. \qed

\section{Appendix B: The Olshanskii criterion and
non-orderability of $S^3$}\label{app-olsh}

In this Appendix we prove that the cone $C$ in the Lie algebra of
$PU(2,1)$ induced from the cone of non-negative contact
Hamiltonians under the natural monomorphism $PU(2,1) \to
\Cont(S^3)$ does not generate a genuine partial order. For the
proof we apply Olshanskii's criterion \cite{Ol}. The text below is
based in parts on the notes which were kindly supplied to us by G
Olshanskii. These notes contain a translation of this criterion,
which is formulated in \cite{Ol} in the abstract language of Lie
theory, into the language of matrices.  Below we present an
exposition of Olshanskii's criterion in the context of $PU(2,1)$.

Consider the form $q(u,w)$ on $\C^3$ given by
$$q(u,w) = u_1 \bar{w}_1 + u_2 \bar{w}_2 - u_3\bar{w}_3\;.$$
By definition, $U(2,1)$ is the group of all $\C$--linear
transformations preserving $q$. In particular, this group
preserves the set $$\{u \in \C^3\;:\; q(u,\bar{u}) < 0\} \;.$$
Passing to the projectivization and setting $z_1=u_1/u_3$ and $z_2
= u_2/u_3$, we get that the group $PU(2,1) := U(2,1)/S^1$ acts by
biholomorphic automorphisms of the ball
$$B^4: = \{z \in C^2\;:\; |z_1|^2 + |z_2|^2 < 1\}\;.$$
These automorphisms extend to diffeomorphisms of the boundary
sphere $S^3$. They obviously preserve the field of complex tangent
lines to $S^3$, and, therefore,  act by contactomorphisms of the standard contact $S^3$.

We identify the real Lie algebra ${\gothic g}$ of $PU(2,1)$ with
$su(2,1)$. It will be useful to describe $\gothic g$ explicitly .
It consists of those traceless complex matrices $A$ of the order
$3 \times 3$ which satisfy $A^*I +IA=0$, where
\[I =\left(
\begin{array}{ccc}
1 & 0 & 0\\
0 & 1 & 0\\
0 & 0 & -1 \end{array} \right)\;,\] and $A^*$ stands for the
complex conjugate of the transposed matrix to $A$. Introduce the
following basis in $\gothic g$:
\[E_1 =\left(
\begin{array}{ccc}
i & 0 & 0\\
0 & 0 & 0\\
0 & 0 & -i \end{array} \right)\;,\; E_2 =\left(
\begin{array}{ccc}
0 & 0 & 0\\
0 & i & 0\\
0 & 0 & -i \end{array} \right)\;,\]
\[F =\left(
\begin{array}{ccc}
0 & 1 & 0\\
-1 & 0 & 0\\
0 & 0 & 0 \end{array} \right)\;,\; {\widetilde F} =\left(
\begin{array}{ccc}
0 & i & 0\\
i & 0 & 0\\
0 & 0 & 0 \end{array} \right)\;,\]
\[G_1 =\left(
\begin{array}{ccc}
0 & 0 & 1\\
0 & 0 & 0\\
1 & 0 & 0 \end{array} \right)\;,\; {\widetilde G}_1 =\left(
\begin{array}{ccc}
0 & 0 & i\\
0 & 0 & 0\\
-i & 0 & 0 \end{array} \right)\;,\]
\[G_2 =\left(
\begin{array}{ccc}
0 & 0 & 0\\
0 & 0 & 1\\
0 & 1 & 0 \end{array} \right)\;,\; {\widetilde G}_2 =\left(
\begin{array}{ccc}
0 & 0 & 0\\
0 & 0 & i\\
0 & -i & 0 \end{array} \right)\;.\]

Consider a Cartan subalgebra ${\gothic h} = \Span_{\R}(E_1,E_2)$.
Our next task is to describe explicitly the cone $C \cap {\gothic
h}$. Let us introduce symplectic polar coordinates
$$\rho_k = \pi|z_k|^2, \varphi_k = \frac{1}{2\pi}{\mathrm Arg}(z_k),\;k=1,2,$$
in $\C^2$. With this notation the standard contact form on $S^3$
can be written as $$\beta= \rho_1 d\varphi_1 + \rho_2d\varphi_2\;.$$
Note that the element $E= a E_1 + b E_2 \in {\gothic h}$
generates the flow
\[u \mapsto \left(
\begin{array}{ccc}
e^{iat} & 0 & 0\\
0 & e^{ibt}  & 0\\
0 & 0 & e^{-i(a+b)t} \end{array} \right)u\] on $\C^3$ and,
therefore, the flow
\[z \mapsto \left(
\begin{array}{cc}
e^{i(2a+b)t} & 0 \\
0 & e^{i(a+2b)t} \end{array} \right)z\;\] on $\C^2$. When
restricted to $S^3$, the latter flow  is generated by the vector
field
$$X= \frac{1}{2\pi} (2a+b)\frac{\partial}{\partial \varphi_1}+
\frac{1}{2\pi} (a+2b)\frac{\partial}{\partial \varphi_2}\;.$$ By
definition, $E \in C$ if and only if $\beta(X) \geq 0$ everywhere
on $S^3$, which  means that both coefficients $2a+b$ and $2b+a$
are non-negative. Summing up, we get
\begin{equation} \label{eq-ap-poscone}
C \cap {\gothic h} = \{(a,b)\;:\; 2a+b \geq 0,\; a+2b \geq 0\},
\end{equation}
where  $(a,b)$ are  coordinates which correspond to the basis $E_1,E_2$ of
$\gothic h$.

Denote by ${\gothic g}_{\C}$ the complexified Lie algebra
${\gothic g} \otimes \C$. In order to avoid the conflict of
notation, we denote by $j$ the natural complex structure on
${\gothic g}_{\C}$, that is we set $jx:= x \otimes i $
Following the notation from \cite{Ol}, we put
$${\gothic h}_{Re} = j{\gothic h}\subset {\gothic g}_{\C}\;.$$

\begin{thm}\label{thm-ol}{\rm\cite{Ol}}\qua  There exists a cone
$c_0 \subset {\gothic h}_{Re}$ with the following property: Given
any convex cone $K \subset {\gothic g}$ which is invariant under
the adjoint action, it defines a genuine partial order on the
universal cover of the Lie group if and only if $$jK \cap {\gothic
h}_{Re} \subset {\pm} c_0\;.$$
\end{thm}

In view of \eqref{eq-ap-poscone} we have to apply this
theorem to the cone
\begin{equation} \label{eq-ap-poscone-on}
c = \{(a,b)\;:\; 2a+b \geq 0,\; a+2b \geq 0\}\;,
\end{equation}
where the $(a,b)$ coordinates on ${\gothic h}_{Re}$ correspond to
the basis $jE_1,jE_2$.

Next, we present Olshanskii's algorithm which enables one to
describe explicitly the cone $c_0$ appearing in \fullref{thm-ol} in terms of the structure theory of the Lie algebra
$\gothic g$.

\medskip
{\bf Structure theory of $su(2,1)$} 

We start with the explicit form of operators $ad(E_1)$ and
$ad(E_2)$. In order to describe these operators it will be convenient
to introduce the matrix
\[J= \left(
\begin{array}{cc}
0 & -1 \\
1 &  0 \end{array} \right)\] and to consider the direct sum
decomposition
$${\gothic g} = \Span_{\R} (E_1,E_2)\oplus \Span_{\R}(F,{\wt F})
\oplus \Span_{\R}(G_1, {\wt G}_1) \oplus \Span_{\R}(G_2, {\wt
G}_2)\;.$$ One readily calculates that $ad(E_1)$ is given by the
matrix
$$0 \oplus J \oplus 2J \oplus J\;,$$ and
$ad(E_2)$ is given by the matrix
$$0 \oplus (-J) \oplus J \oplus 2J\;.$$
Up to a positive multiple, the Killing form $Q$ restricted to
${\gothic h}_{Re}$ is given by the matrix
\[ \left  (
\begin{array}{cc}
2 & 1 \\
1 &  2 \end{array} \right)\;\] in the basis $(jE_1,jE_2)$. The
space ${\gothic h}_{Re}$ can be identified with the dual to
${\gothic h}$ with the help of the Killing form. Hence, we can
assume that the roots of $({\gothic g}_{\C}, {\gothic h}_{\C})$
lie in ${\gothic h}_{Re}$. Diagonalizing operators $ad(E_1)$ and
$ad(E_2)$ we calculate that in the above basis the roots have the
form
$$\pm \gamma = \pm  (1,-1), \;\pm \alpha_1 = \pm (1,0),\; \pm \alpha_2 = \pm (0,1)\;.$$ The
root system is of the type $A_2$, and one can declare
$\gamma,\alpha_1,\alpha_2$ to be positive roots.

Next, we have the Cartan decomposition ${\gothic g} = {\gothic t}
\oplus {\gothic p}$ into the direct sum of the subspaces formed by
skew-Hermitian and Hermitian matrices respectively. Thus ${\gothic
t}= \Span_{\R}(E_1,E_2,F,{\wt F})$ and ${\gothic p}=
\Span_{\R}(G_1,{\wt G}_1,G_2, {\wt G}_2)$. Recall that a root is
called non-compact if the corresponding eigenspace is contained in
${\gothic p} \otimes \C$. It is easy to verify that the positive
non-compact roots are given by $\alpha_1$ and $\alpha_2$. Denote
by $c_{min}$ the cone in ${\gothic h}_{Re}$ generated by these
roots:
\begin{equation} \label{eq-ap-poscone-tw}
c_{min} = \{(a,b)\;:\; a \geq 0,\; b \geq 0\}\;,
\end{equation}
where the $(a,b)$--coordinates on ${\gothic h}_{Re}$ correspond to
the basis $jE_1,jE_2$.

Note that the center $\zeta$ of the algebra ${\gothic t}$ is
generated by the  vector $E_1+E_2$. Furthermore, the roots of the
algebra ${\gothic t}_{\C}$ with respect to the Cartan subalgebra
${\gothic h}_{\C}$ are given by $\pm \gamma$. Let  $W$ be  the
corresponding Weyl group which acts on ${\gothic h}_{Re}$ by the
permutation of coordinates in the basis $(jE_1,jE_2)$.

{\bf Construction of the cone $c_0$} 

We work with the basis $(jE_1,jE_2)$ of ${\gothic h}_{Re}$.

\medskip{\bf Step 1}\qua Choose a maximal subsystem of positive non-compact
pairwise orthogonal roots. Since $\alpha_1$ and $\alpha_2$ form
the angle $\pi/3$, we choose only one of them, say $\alpha_1$.
Denote by $H_1$ the vector of the form $s\alpha_1$ with $s>0$ so
that
$$Q(\alpha_1,H_1) = 2\;.$$ We calculate that $H_1 = (1,0)$.

\medskip{\bf Step 2}\qua Choose the vector $Z \in j\zeta$ such that
$Q(\alpha_1,Z)=Q(\alpha_2,Z)=2$. We calculate that $Z =
(2/3,2/3)$. Put $$H_0 = Z-H_1 = \frac{1}{3} \cdot (-1,2)\;.$$

\medskip{\bf Step 3}\qua Consider the cone $c_1$ spanned by $WH_0 \cup
c_{min}$. Note that $c_1$ is spanned by vectors
$$\frac{1}{3} \cdot (-1,2) \;\;\text{and}\;\;\frac{1}{3} \cdot
(2,-1)\;.$$ The cone $c_0$ from \fullref{thm-ol} is defined as
the cone whose dual with respect to the  form $Q$ equals
$c_1$:
$$c_1 = \{x \in {\gothic h}_{Re}\;:\; Q(x,y) \geq 0 \;\;\forall y
\in c_0\}\;.$$ One readily calculates that $c_0 = c_{min}$, where
$c_{min}$ is given by inequalities \eqref{eq-ap-poscone-tw}.

\medskip{\bf The final step}\qua  Recalling that the cone $c = jC
\cap {\gothic h}_{Re}$ is given by inequalities
\eqref{eq-ap-poscone-on} we get that $c_0$ is strictly contained in
$c$. Thus, by \fullref{thm-ol} the cone $C$ does not generate
a genuine partial order. \qed

\begin{rem} \label{rem-aptw-on} {\rm
Olshanskii's method does not give rise to an explicit description
of a positive contractible loop of contactomorphisms which lies in
$PU(2,1) \subset \Cont(S^3)$. However,  imitating our construction in \fullref{subsec-main-consttw}
above, we  present here  such a loop. Let us identify $PU(2,1)$ with the group of   complex automorphisms of
$B^4$.  We define a distinguished contactomorphism (cf \fullref{subsec-dist-cont} above) by the formula
$$b(z_1,z_2) = \left(\frac{\cosh \alpha \cdot z_1 +1}{z_1 + \cosh
\alpha}\;, \frac{ \sinh \alpha \cdot z_2}{z_1 + \cosh
\alpha}\right)\;.$$ Put $ e_tz:= e^{2\pi i t }z$ and $$
f_t(z_1,z_2) = (e^{2\pi i t} z_1, e^{-2\pi i t}z_2),$$ as in
\fullref{subsec-main-consttw} above. A lengthy but
straightforward calculation shows that the (obviously
contractible) loop $e_{-t}f_{3t}be_{t}b^{-1}$  (cf \fullref{thm-main-consttw})  is positive,  provided that $\alpha>0$ is small
enough.

In order to reveal the geometry of the distinguished
contactomorphism $b$, let us perform the Cayley transform
\cite[Chapter 2.3]{Rudin}
$$\Phi\co  B^4 \to \Omega,\;\; (z_1,z_2) \mapsto
\frac{i}{1-z_1}(z_1+1,z_2)\;,$$ where
$$\Omega= \{(w_1,w_2)\in \C^2\;:\; \mathrm{Im}\;w_1 > |w_2|^2\}\;$$
is a multi-dimensional analogue of the upper half-plane. This
transform is a biholomorphism between
$\mathrm{Closure}(B^4)\setminus \{(1,0)\}$ and
$\mathrm{Closure}(\Omega)$. Let $\eta$ be the field of complex
tangent lines to $\partial \Omega$. Put $w_1=x+iy, w_2 = u +iv$.
One readily shows that in the coordinates $(x,u,v)$ the contact
manifold $(\partial \Omega ,\eta)$ is simply $\R^3$ equipped with
the standard contact structure
$$\mathrm{Ker}(dx + 2(udv-vdu))\;.$$
The map $\hat{b} := \Phi b\Phi^{-1}\co  \Omega \to \Omega$ turns out
to be a {\it non-isotropic dilation} $$(w_1,w_2) \mapsto (s^2 w_1,
sw_2),$$ where $s \to +\infty$ as $\alpha \to 0$. Of course, the
restriction of $\hat{b}$ to the boundary has the friendly form
$$(x,u,v) \to (s^2x,su,sv)\;.$$}
\end{rem}

\begin{rem}\label{rem-aptw-tw}{\rm It is unclear to us
whether the explicit analytic formula for the positive
contractible loop in $PU(2,1) \subset \Cont(S^3)$ presented above
helps to simplify our calculations in \fullref{sectw-stab},
namely to "extend" this loop to a nonnegative contractible loop in
stabilizations, and to calculate the sharp lower bound for the
invariant $\mu(\Delta)$ defined by formula \eqref{eq-mu}. }
\end{rem}

\begin{rem}\label{rem-aptw-th}{\rm Oshanskii's paper
\cite{Ol} provides information on the orderability of very general
finite-dimensional simply-connected Lie groups in terms of the
structure theory. This subject was further developed in
\cite{Neeb,Hilgert}. It would be interesting to explore the
notions introduced in \fullref{subsec-semigroup} in this
context.}
\end{rem}

\medskip
{\bf Acknowledgements} 

We cordially thank Ms  M~Hercberg for preparing LaTex figures for
this paper. We are indebted to E~Giroux for pointing out that our
\fullref{thm-nogo} follows from the Olshanskii criterion. We are
grateful to G~Olshanskii for providing illuminating notes which
helped us a lot to understand his paper \cite{Ol} and for useful
comments on the first draft of Appendix B. We thank J~McCleary for
providing us with reference \cite{FTVP}, R~Cohen, S Goberstein,
M~Sapir and S~Weinberger for useful consultations, D~McDuff,
F~Schlenk and the anonymous referee for useful comments on the
manuscript.

 Our numerous meetings at Stanford University, American Institute of
Mathematics (Palo Alto), Institute of Advanced Studies (Princeton) and
Tel Aviv University were indispensable for completing this project. We
express our gratitude to these institutions. The third named author is
grateful to the Center for Dynamics and Geometry at Pennstate for the
opportunity to lecture on preliminary results of the present paper in
Spring 2005, and in particular to D~Burago and S~Tabachnikov for
their warm hospitality.

In its early stages, this project was supported in part by United
States--Israel Binational Science Foundation grant number 1999086.

The authors acknowledge the support of the following grants:

{\sl Y Eliashberg}\qua NSF grants
DMS--0204603 and DMS--0244663\\ {\sl S\,S Kim}\qua  NSF grant
DMS--972992, WISE fellowship from the University of Southern
California and FCT/SFRH/BPD (Portugal)\\ {\sl L Polterovich}\qua
Michael Bruno Memorial Award

\bibliographystyle{gtart}
\bibliography{link}

{\small\parskip 0pt\vskip11pt minus 5pt\relax
{\sl \def\\{\futurelet\next\nocommawithnl}\def\nocommawithnl
  {\ifx\next\newline\else\unskip,\space\ignorespaces\fi}
  \theaddress\par}
{\rightskip0pt plus .4\hsize
{\def\tempab{}\tt\def~{\lower3.5pt\hbox{\char'176}}\def\_{\char'137}%
\ifx\theemail\tempab\else
  \vskip5pt minus 3pt\theemail\par\fi
  \ifx\theurl\tempab\else
  \vskip5pt minus 3pt\theurl\par\fi}
  \vskip11pt minus 5pt
            Proposed:\qua\theproposer\hfill
            Received:\qua\receiveddate\break
            Seconded:\qua\theseconders\hfill
        \ifx\reviseddate\tempab
              Accepted:\qua\accepteddate\else Revised:\qua
  \reviseddate\fi\break}}

\newpage

\title[Erratum to ``Geometry of contact transformations and 
       domains'']{Erratum to ``Geometry of contact transformations\\and domains: 
       orderability versus squeezing''}

\volumenumber{13}
\issuenumber{2}
\publicationyear{2009}
\papernumber{29}
\startpage{1175}
\endpage{1176}

\lognumber{1291}

\doi{}
\MR{}
\Zbl{}

\received{1 February 2009}
\accepted{1 February 2009}
\published{1 February 2009}
\publishedonline{1 February 2009}
\count0=1175

\def\theabstract{}
\makeatletter
\let\@primclass\relax
\let\@secclass\relax
\makeatother

\maketitle
\section*{}\vglue -80pt
\addcontentsline{toc}{section}{Erratum}
\hypertarget{Err}{The} 
purpose of this erratum is to correct a number of inconsistencies
in our paper (above).  These are related to the grading of
generalized Floer homology and do not affect formulations and
proofs of the main results of the paper. The source of the mistakes
in the grading is as follows.  Consider the framed Hamiltonian
structure
$$(V=M \times
S^1, \Theta = -\alpha + H_t dt,\lambda= dt)$$ associated to a
Hamiltonian function $H:M \times S^1 \to \R$ on an exact symplectic
manifold $(M,\omega=d\alpha)$, as in Example 4.2.2. The symplectic
form $\Omega = d\Theta$ on the cut of this Hamiltonian structure
equals the form $-\omega$ (mind the minus sign). In the proof of
Proposition 4.37 we claim that  ``the change of the sign of the
symplectic form does not affect Conley--Zehnder indices", which is
wrong for the definition of Conley--Zehnder indices as presented in
Section 4.4.1. In order to correct this, the following changes
should be made in paper:

\medskip
\centerline{\sc List of corrections}

\medskip
\hyperlink{Corr1}{Page $1677$, line $14$:} replace $2k$ by  $-2k$.

\hyperlink{Corr2}{Page $1677$, line $-8$:} replace $g_tR_*^{tT}g_0^{-1}$ by
$(g_tR_*^{tT}g_0^{-1})^{-1}$.

\hyperlink{Corr3}{Page  $1698$, line $-8$:} \hypertarget{Corr3R}{Add} the following text: Let us
emphasize that in the usual Floer homology the CZ--index of the orbit
is defined through the linearization of the Hamiltonian flow along
the orbit. Furthermore,  the action of an integer $k \in \Z$ on the
set of coherent symplectic trivializations corresponds to the
twisting by a loop of symplectic matrices with the Maslov index
$+2k$.  Thus, informally speaking, the gradings of closed orbits in
the generalized Floer homology for stable Hamiltonian structures and
for the usual Hamiltonian Floer homology are defined in the
``opposite way".

\hyperlink{Corr4}{Page 1699, line $-1$:} Replace
``\hypertarget{Corr4R}{since} the change of the sign of the symplectic
form does not affect Conley--Zehnder indices" to ``due to our
conventions on the Conley--Zehnder index for generalized and usual
Floer homologies".

With these corrections, the rest of the paper remains
unchanged.

\medskip
Let us emphasize also that in \fullref{thmloop-space} the isomorphism
$L_r$ between symplectic homology of a tube in $T^*X$ and cohomology
of the corresponding sublevel set of the energy functional on the
loop space $\cL X$ does not preserve grading. It follows from Salamon
and Weber \cite{sw}
that if the manifold $X$ is orientable and the degrees of
non-contractible orbits are determined by the vertical Lagrangian
subbundle of $T(T^*X)$, the isomorphism $L_r$ sends
$\SH_k^{(-\infty, -1)}(U_r)$ to $H^{n-k}(\cL^{\frac1{2r^2}} X)$
where $n = \dim X$.

Let us also point out that our definition of the Conley--Zehnder
index for paths of symplectic matrices in $\R^{2n}$ differs from the
one defined by Robbin and Salamon in \cite{rs}: the sum of two indices
equals $n$.

\end{document}